\documentclass[aoas, reqno]{imsart}
\RequirePackage[OT1]{fontenc}
\RequirePackage{amsmath}
\RequirePackage{natbib}
\RequirePackage[colorlinks,citecolor=blue,urlcolor=blue]{hyperref}
\RequirePackage{hypernat}
\usepackage{algorithm}
\usepackage{algorithmicx}
\usepackage{algpseudocode}
\usepackage{natbib}
\usepackage{graphicx}
\usepackage{epstopdf}
\usepackage{rotating}
\usepackage{subfigure}
\usepackage{longtable}
\usepackage{engord}
\usepackage{url}
\usepackage{enumerate}
\usepackage{color}
\usepackage{layout}
\usepackage{theorem}
\usepackage{amssymb}
\usepackage{float}
\usepackage{hhline}
\usepackage{version}
\usepackage{fancyheadings}
\usepackage{appendix}
\usepackage{epstopdf}
\usepackage{chemarrow}
\usepackage{amsmath}
\usepackage{array}
\usepackage{soul}
\usepackage{hyperref}
\usepackage{tabularx}
\usepackage{multirow}

\newcommand {\ctp}{\cite}

\newtheorem{thm}{Theorem}

\newtheorem{remark}{Remark}

\begin{document}

\begin{frontmatter}
\title{Estimating the Rate Constant \\ from Biosensor Data via \\ an Adaptive Variational Bayesian Approach}
\runtitle{Estimation of Rate Constant by Biosensor Data}

\begin{aug}
\author{\fnms{Ye} \snm{Zhang}\thanksref{t1}\ead[label=e1]{ye.zhang@oru.se}},
\author{\fnms{Zhigang} \snm{Yao}\thanksref{t2}\ead[label=e2]{zhigang.yao@nus.edu.sg}},
\author{\fnms{Patrik} \snm{Forss\'{e}n}\thanksref{t3}\ead[label=e3]{patrik.forssen@kau.se}}
\and
\author{\fnms{Torgny} \snm{Fornstedt}\thanksref{t3}\ead[label=e4]{torgny.fornstedt@kau.se}}

\thankstext{t1}{The work of the first author is supported by the Swedish Knowledge Foundation (No. 20170059) and the Alexander von Humboldt foundation.}
\thankstext{t2}{The work of the second author is supported by the MOE grants Tier 1 R-155-000-196-114 and Tier 2 R-155-000-184-112 at the National University of Singapore.}
\thankstext{t3}{The work of the last two authors is supported by the Swedish Research Council (No. 2015-04627) and the Swedish Knowledge Foundation (No. 20170059).}
\thankstext{t4}{Correspondence should be addressed to: zhigang.yao@nus.edu.sg.}
\runauthor{Y. Zhang, Z. Yao, P. Forss\'{e}n and T. Fornstedt}

\affiliation{Chemnitz University of Technology,  National University of Singapore \and Karlstad University }

\address{Y. Zhang\\
Faculty of Mathematics\\
Chemnitz University of Technology\\
09107 Chemnitz, Germany\\
\printead{e1}
}

\address{Z. Yao\\
corresponding author\\
Department of Statistics and Applied Probability\\
National University of Singapore\\
21 Lower Kent Ridge Road\\
Singapore 117546\\
\printead{e2}
}

\address{P. Forss\'{e}n and T. Fornstedt\\
Department of Engineering and Chemical Sciences\\
Karlstad University \\
65188  Karlstad, Sweden\\
\printead{e3}\\
\printead{e4}
}

\end{aug}

\begin{abstract}
The means to obtain the rate constants of a chemical reaction is a fundamental open problem in both science and the industry. Traditional techniques for finding rate constants require either chemical modifications of the reactants or indirect measurements. The rate constant map method is a modern technique to study binding equilibrium and kinetics in chemical reactions. Finding a rate constant map from biosensor data is an ill-posed inverse problem that is usually solved by regularization. In this work, rather than finding a deterministic regularized rate constant map that does not provide uncertainty quantification of the solution, we develop an adaptive variational Bayesian approach to estimate the distribution of the rate constant map, from which some intrinsic properties of a chemical reaction can be explored, including information about rate constants. Our new approach is more realistic than the existing approaches used for biosensors and allows us to estimate the dynamics of the interactions, which are usually hidden in a deterministic approximate solution. We verify the performance of the new proposed method by numerical simulations, and compare it with the Markov chain Monte Carlo algorithm. The results illustrate that the variational method can reliably capture the posterior distribution in a computationally efficient way. Finally, the developed method is also tested on the real biosensor data (parathyroid hormone), where we provide two novel analysis tools~-- the thresholding contour map and the high order moment map~-- to estimate the number of interactions as well as their rate constants.
\end{abstract}


\end{frontmatter}

\section{Introduction}
\label{intro}

In the modern world, biosensors have made a significant impact in many fields, such as antibody-antigen interactions, immunology, virology, and the pharmaceutical industry \citep{Sanvicens-2011,Nicholls-2015}. Hence, during the last few decades, there has been an accelerated technological development of biosensor instruments, e.g., surface plasmon resonance, quartz crystal microbalance, etc. A simplified biosensor system is presented in Figure \ref{Biosensor}, and its physical mechanism is briefly discussed in Appendix A. To design an appropriate biosensor instrument that is biocompatible or specifically functionalized, scientists must know the physical chemistry of biomolecular/cell surface interactions \citep{Telesca-2012}. The reliable analysis of biomolecular interactions is crucial in both science and the industry, e.g., it is required to fulfill modern drug quality assurance criteria. In order to understand the interactions, scientists collect biosensor data that measures the analyte biomolecules of several different concentrations on a sensor chip with immobilized ligand molecules that form complexes with analytes. This kind of biosensor data is usually called a sensorgram, where the systems response, proportional to total complex concentration, is measured over time for different analyte injections. In this paper, we focus on this particular type of data and aim to obtain information about the interactions, i.e., their numbers and the corresponding rate constants. The biosensor data is traditionally processed using a simple model fitting procedure assuming just one, or perhaps two, distinct interactions. These traditional approaches might not reflect the true, complicated, and heterogeneous molecular interactions for large active pharmaceutical ingredients, and can thus lead to wrong mechanistic conclusions. Therefore, a more advanced analysis of the biosensor data is necessary to avoid the costly and time consuming procedure of repeating the same ligand-analyte binding experiment over and over using chips with different ligand binding or different buffers, etc., until one of the tests results in ``good enough'' data that works with the standard rate constant estimation \citep{Gray-1985}.

\begin{figure}[!b]
\centering
\includegraphics[width=4in]{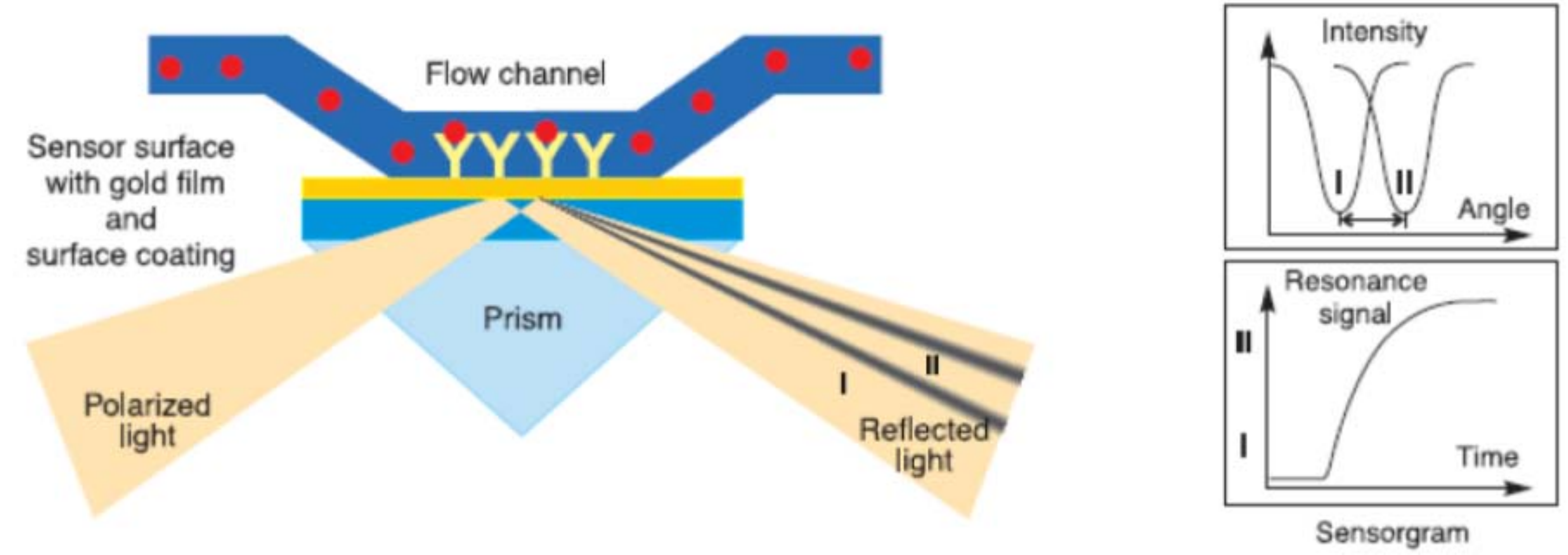}
\caption[A simplified biosensor system]{A simplified biosensor system.}
\label{Biosensor}
\end{figure}

\subsection{Existing work in Statistics and Computational Chemistry}

The problem of estimating the chemical rate constants of interactions through statistical approaches has been investigated for long time \citep{Box1965,Atherton1975,Ziegel1980,Komorowski2011,Gupta2014,Pischel2017}. To our best knowledge, most of the existing statistical approaches require a priori knowledge of the number of interactions in a chemical reaction, i.e., they first assume the existence of parallel reaction (due to potential interactions), and then estimate the rate constants for corresponding interactions. Therefore, in practice, the problem is usually solved by a two-step strategy. The first step is to determine the number of interactions by dimensionality reduction technique (e.g., principal component analysis, regression tree \citep{Loh2013}), and then fit the data to the parameters (i.e., rate constants correspond to these reactions). It is clear that the first step partially determines the quality of the results in the second step. However, it can be observed that the second step procedure also potentially has a positive influence on the estimation in the first step, particularly when the data is under-sampled and affected by large noise. Therefore, the separation of the procedure into two steps does not seem optimal. The question that now arises is whether the number of interactions and the values of rate constants corresponding to each interaction in a single step can be estimated from the biosensor data, therefore making the modeling more robust. These considerations form the motivations behind our work.

Another line of research, inspired by Adsorption Energy Distribution (AED) calculations \citep{Guiochon-1995}, which are successfully used for steady-state data usually from small molecular systems, as well as a couple of numerical solvers for dynamic biosensor data, such as the Interaction Distribution (ID) (or Interaction Map) \citep{Svitel-2003,Altschuh-2012} and Adaptive Interaction Distribution Algorithm (AIDA) \citep{AIDA}, have recently gained much attention. With these solvers, one can see the number of interactions in the system, along with their associated rate constants, as peaks on a surface. The mathematical basis behind these solvers is to use a single step model to estimate the number of rate constants and their values simultaneously. They utilize regularization methods to impose a sense of well-posedness for the model and solve the resulting deterministic problems by optimization. However, the result of these methods strongly depends on the choice of system parameters, e.g., the type of regularization penalty and the value of the regularization parameter. For our application, as well as other problems where no ground truth is available, the problem of selecting the regularization penalty term and the regularization parameter remains an open issue \citep{FanLv}, as most tuning parameters such as the generalized cross-validation (GCV) \citep{Tibshirani}, $C_p$ style statistics \citep{Efron}, and the consistent estimator of degrees of freedom of the LASSO in the $C_p$, AIC, and BIC criteria \citep{Zou2007}, are based on an asymptotic property that is practically infeasible.

The shortage of existing methods motivates us to combine statistical and deterministic methods to develop a more robust single step approach via a statistical strategy. Note that this high level idea has been used in \cite{Gorshkova-2008} to some extent. However, from a mathematical viewpoint, their method can still be classified in the branch of deterministic methods, since it estimates the distribution of rate constant directly by solving a special regularized optimization problem. Though the authors have taken into account the contribution of a prior expectation of rate constants on the regularization term in their optimization formulation, the mathematical model should be deterministic if one considers the distribution of rate constant, which is also termed as the rate constant map, as a deterministic object (i.e. an element in the space of functions). Furthermore, the justification of rate constants as random variables is missing, and no uncertainty quantification is discussed under their framework. Meanwhile, the difficulty of choosing an appropriate regularization parameter remains. This motivates us to regard the rate constant map as a random variable, which incorporates the regularization in a completely different way, and overcomes the difficulty of selecting regularization parameter, hence may carry the intrinsic information from the uncertainty quantification of the rate constant map.

\subsection{Motivation}\label{Motiva}

\begin{figure}[!htb]
\centering
\subfigure[]{
\includegraphics[clip, trim=0 0.02in 0 0.02in, width=2.2in]{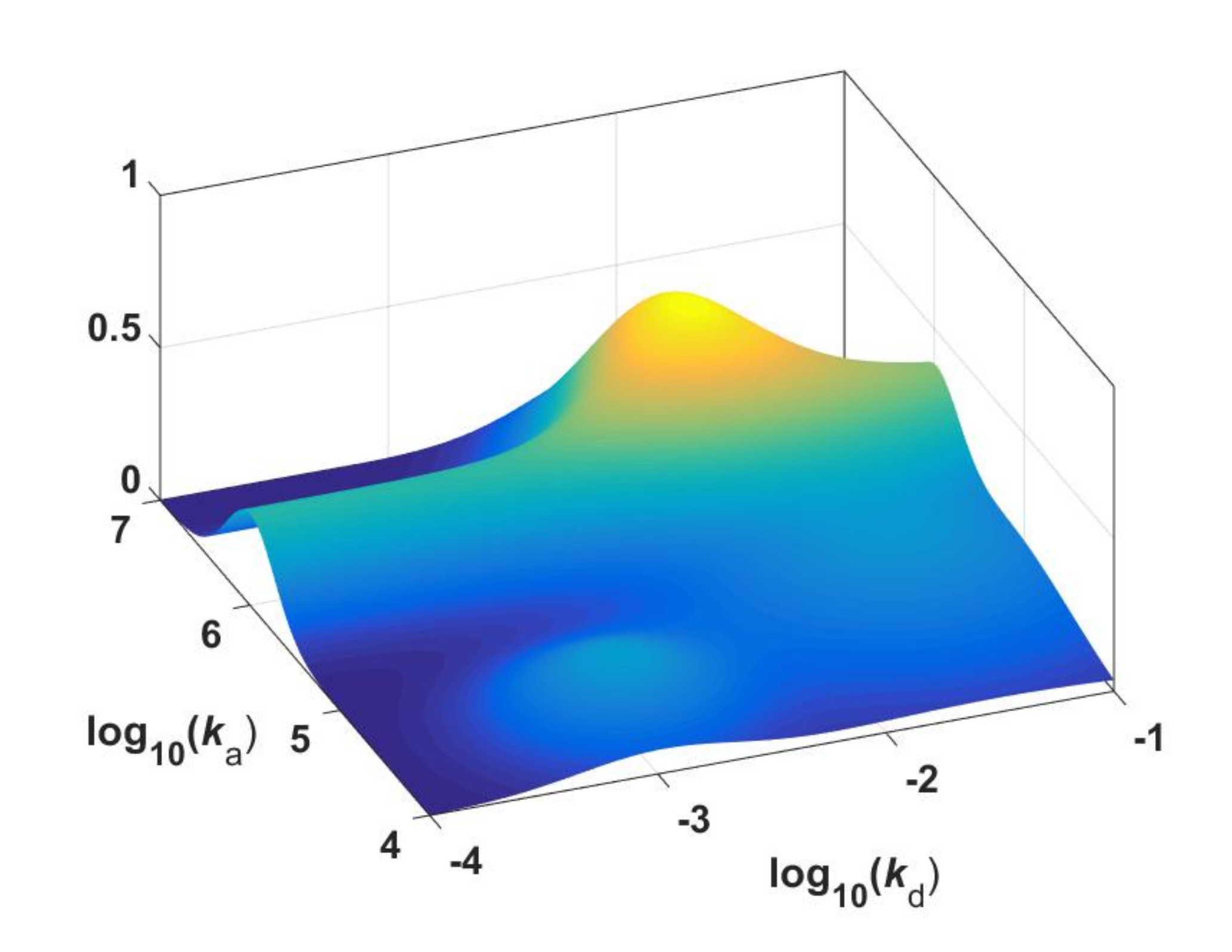}}
\subfigure[]{
\includegraphics[clip, trim=0 0.02in 0 0.02in, width=2.2in]{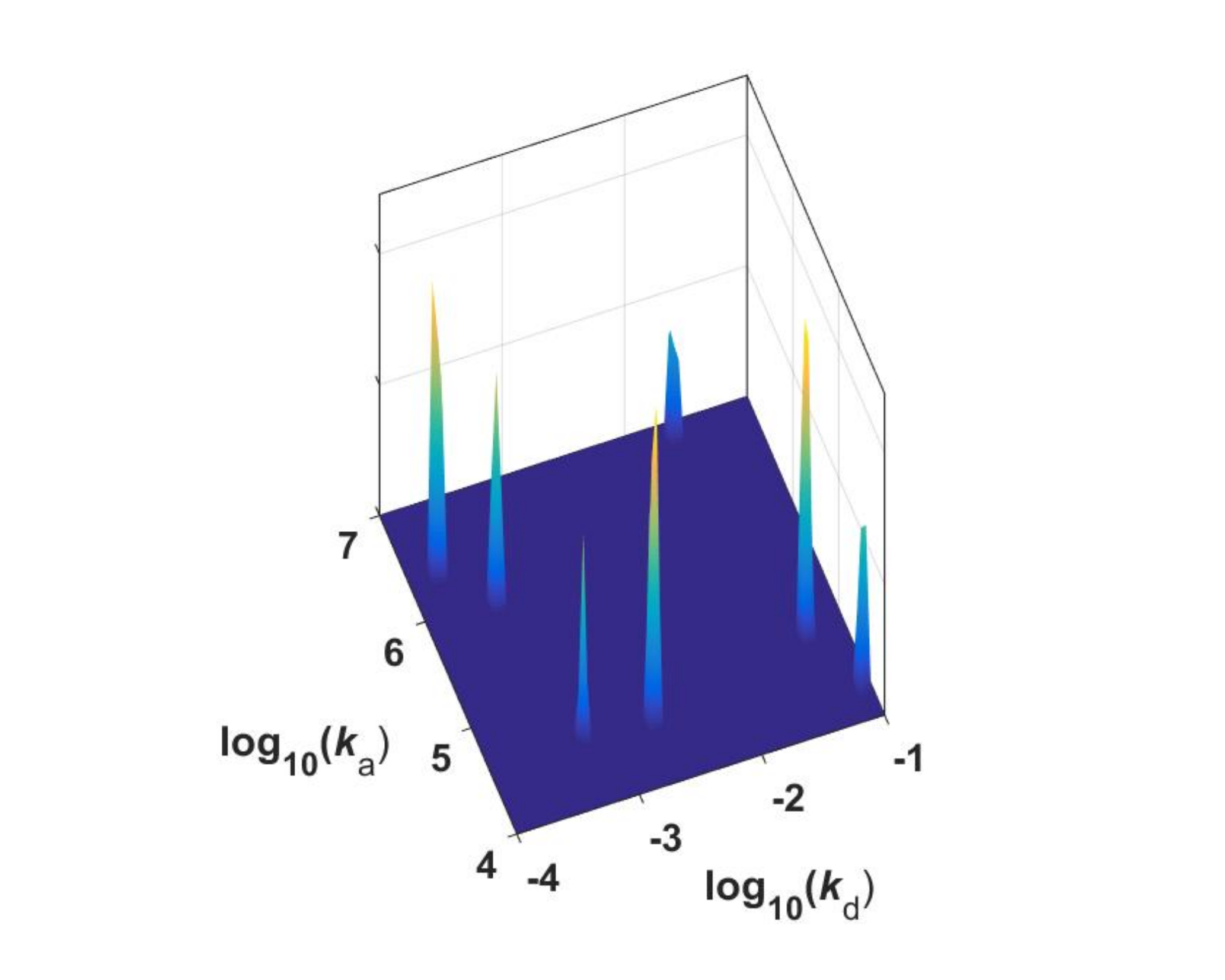}}
\caption[Exact and expected rate constants]{(a) The Tikhonov method. The estimated rate constant map is too smooth to enable the deduction of the number of interactions and the value of association and dissociation constants (the position of peaks in the rate constant map). (b) The sparse regularization method ($\ell_1$ regularization). The information of the interaction is completely hidden in the estimated rate constant map.}
\label{deterministic}
\end{figure}

As discussed above, due to the diversity of data structure and the appearance of noise in many biosensor problems, the existing methods do not work efficiently. The central problem is that most of them are not capable of inferring the interaction information correctly. For instance, in the assumption that two interactions exist in our biosensor system, the exact value of rate constants in our data structure called the dissociation constant $k_d$ (measures the rate at which a molecule complex diassociates) and association constant $k_a$ (rate at which molecules form complexes), see Appendix A for details, are displayed in (a) of Figure \ref{IntroExact}. This means that two parallel interactions exist and their corresponding rate constants are $(\log_{10}(k_d), \log_{10}(k_a))=(-3,4.5)$ and $(-2, 6.5)$. Now, we solve the problem using the conventional rate constant map method, which will be introduced later in Section 2 and Appendix A. This method provides a map of distribution of rate constants (called a rate constant map), from which we can figure out the number of interactions and the rate constants of these interactions. The desired estimated rate constant map should exhibit the information of the interactions. One example of this rate constant map is displayed in (b) of Figure \ref{IntroExact}, where people can easily derive the interaction numbers and rate constants from the peaks of the rate constant map. However, in many cases, the regularized deterministic approaches cannot offer the rate constant map that is correctly shaped. For instance, in our parathyroid hormone application, the Tikhonov regularization provides an oversmoothing effect on the rate constant map, while, the sparse $\ell_1$ regularization offers a rather unsatisfactory reconstruction of the rate constant map, see Figure \ref{deterministic} for details. This motivates us to build new statistical models and investigate the methods that can efficiently estimate the rate constant map with uncertainty quantification, from which we can explore the intrinsic heterogeneity of the biosensor system~-- the interaction information including the number of interactions and the corresponding rate constants. The philosophy of our approach is to assume that the rate constants are no longer static. By utilizing the non-static rate constants, where the rate constant map $f$ is modeled as a high dimensional random variable, we will be able to investigate the distribution of $f$ and provide estimates of the variability. Based on the estimated distribution of $f$, we can produce different kinds of maps (e.g., (b) in Figure \ref{IntroExact} or (a) and (b) in Figure \ref{deterministic}) which provides a multi-scale perspective of the rate constant map. It turns out that the random rate constant map $f$ can be used to find the interaction information, e.g., the number of the interactions and rate constants for each interaction, which will be called the intrinsic quality of a chemical system.

\begin{figure}[!t]
\centering
\subfigure[]{
\includegraphics[clip, trim=0 2.2in 0 2.2in, width=2in]{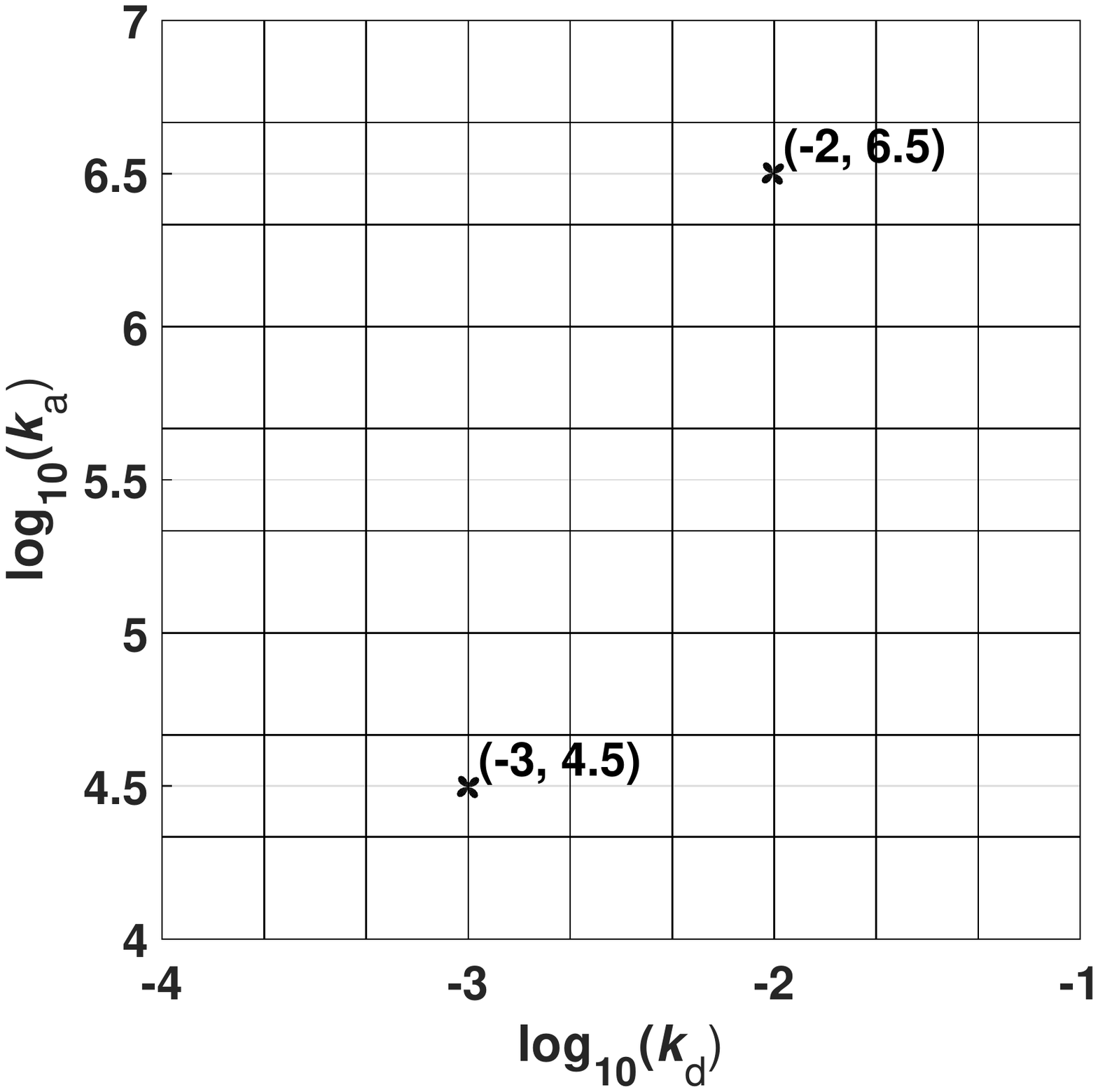}}
\subfigure[]{
\includegraphics[clip, trim=0 0.02in 0 0.02in, width=2.2in]{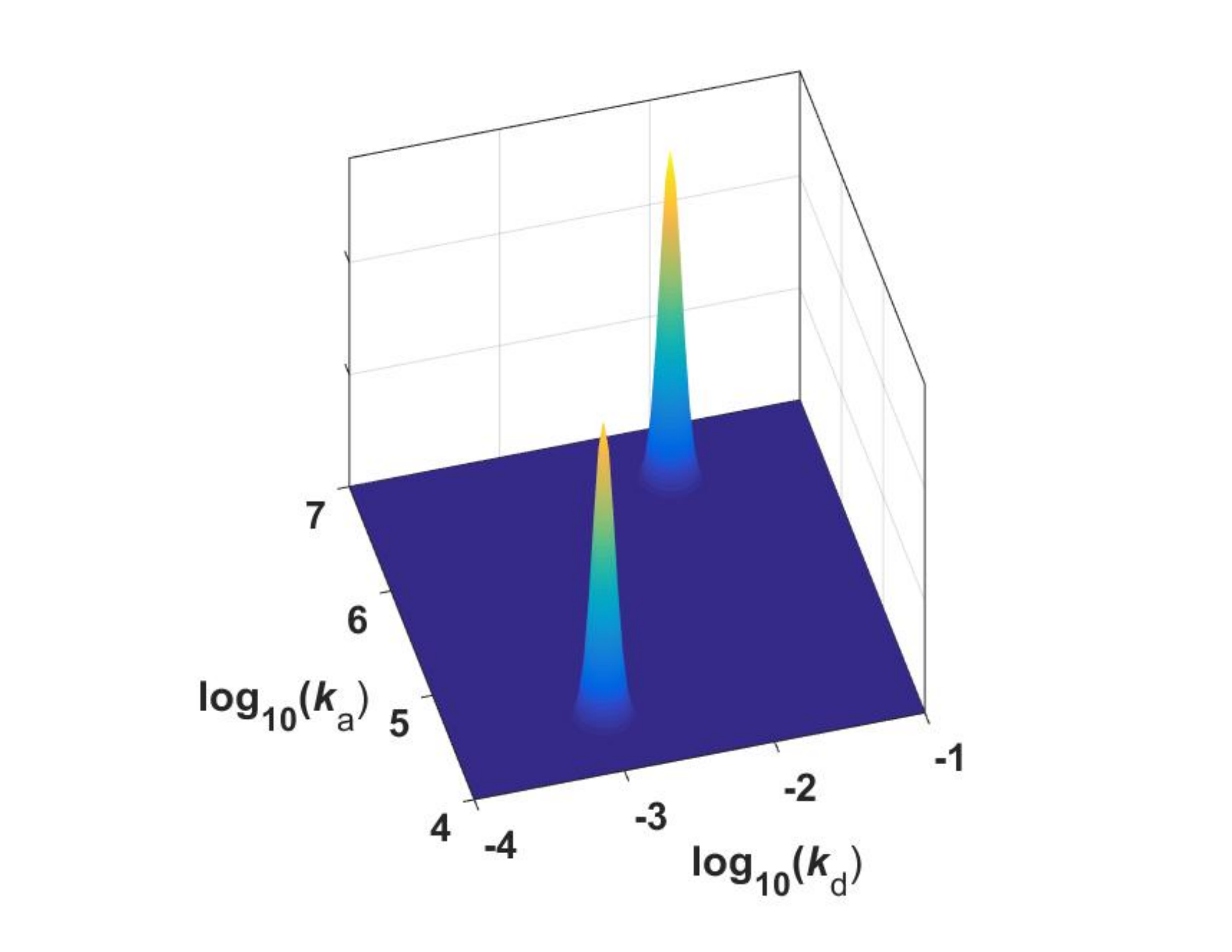}}
\caption[Exact and expected rate constants]{(a) The exact value of association and dissociation constants. (b) One example of the expected rate constant map $f(k_a, k_d)$.}
\label{IntroExact}
\end{figure}

In comparison with classical deterministic regularization methods, Bayesian approaches have distinct features for inverse problems (in our case, the inverse problem is the estimation of interaction information from indirectly measured biosensor data) that make them attractive for statistical inference \citep{Bernardo1994,Kennedy2001,Evans2002,Kaipio2007}. First, the Bayesian posterior distribution provides uncertainty quantification by assigning probabilities to the space of all possible inverse solutions that are consistent with the observed data. This uncertainty quantification is lacking in all the deterministic regularization methods, which only yield one single solution as a point estimate. Second, the use of prior distribution allows more flexibility in modeling. For instance, sparsity-inducing priors such as the spike-and-slab prior can be used if sparse solutions are desirable for high dimensional models. Constraints on the parameter space, or the space of solutions in our context, can be imposed directly through constraining the support of the prior distribution. Third, the tuning parameters that are used in the frequentist regularization methods become hyper-parameters in the Bayesian framework. They are usually assigned some hyper-priors and can hence be made fully adaptive to the observed data. On the other hand, since the posterior distribution does not have an analytically tractable closed form in most cases, Bayesian inference typically requires the numerical sampling of the posterior distributions by Markov chain Monte Carlo (MCMC), which can be computationally intensive and may suffer from slow mixing due to high dimensionality or strong correlations between the parameters in their joint posterior distribution. As an alternative approach, the variational Bayesian (VB) inference \citep{Bishop2006,Blei-2017} has gained popularity for its computational efficiency over posterior sampling methods in such complex Bayesian models. VB approximates the true posterior distribution with a family of computationally tractable distributions by minimizing the KL divergence between them. Although it has been empirically observed that VB occasionally underestimates the posterior variances, the gain in computational efficiency from VB is substantial in many applications; see, e.g., \ctp{Parisi1988}, \ctp{Peterson-1987}, \ctp{Jordan-1999}, \ctp{Jordan-2008}, \ctp{Jin-2010} and references therein.

In this paper, we build a Bayesian framework for uncertainty quantification in the biosensor problem and propose a new VB approach for computing the posterior distribution. Our VB approach is novel and motivated by three physical facts that are unique to the biosensor problem, as well as other similar problems, modeled by the integral equations (\ref{IntegralEq}). First, in our real world biosensor problem, the noise structure can be heterogeneous among different sensorgrams, which results in extremely high dimensionality of the parameter spaces that makes any posterior sampling algorithms such as the MCMC practically infeasible. Therefore, we develop a VB approach to circumvent the computational problem. Second, the solution of the rate constant map must be nonnegative for meaningful physical interpretation. This fact requires a constrained prior support and eventually a new VB algorithm that differs from the existing ones. Finally, the dimensionality of parameter space in most existing statistical models should be fixed, which requires the a priori knowledge of the distribution of rate constants in the biosensor system that is actually unknown before we estimate them. To overcome this difficulty, based on the newly developed VB algorithm we introduce an adaptive discretization technique that automatically adjusts the dimensionality of solution space during the evolution of our algorithm, according to the uncertainty quantification of the obtained temporary qualities.

\subsection{Objective of the paper}
In Section 2, based on the rate constant map theory, we establish a single step model connecting the biosensor data and the rate constant map, which is used for estimating the number of interactions and the rate constants simultaneously. Section 3 presents our newly developed method~-- an Adaptive Variational Bayesian Approach (AVBA)~-- for solving the proposed single step model. In Section 4, we illustrate the approach using an artificial problem, and compare its performance with the MCMC algorithm. The proposed method is tested on real biosensor data (parathyroid hormone) in Section 5, where we will develop two AVBA-based analysis tools that represent the intrinsic property of the interaction, and can be used for estimating the number of interactions as well as their rate constants. Finally, concluding remarks are given in Section 6.

\section{A mathematical model for estimating the rate constant map}
\label{modeling}

As mentioned before, in this work, we use the rate constant map $f(k_a,k_d)$ (i.e., the distribution of association constant $k_{a}$ and dissociation constant $k_{d}$) to explore the interaction information, i.e., their numbers and (active) rate constants. Based on the \emph{Rate Constant Map Theory}, see details in the Appendix A, the measured sensorgram $R_{obs}$, which is dependent on the analyte concentration $C$ and time $t$, and the rate constant map $f$ are related through the following integral equation
\begin{equation}
R_{obs}(t;C) = \int_{\Omega} K(t,C;k_a,k_d) f(k_a,k_d) d k_a d k_d, \quad  (k_a,k_d)\in \Omega,
\label{IntegralEq}
\end{equation}
where $\Omega\subset \mathbb{R}^2$ is the interested domain of rate constants, and the kernel function $K(\cdot)$ is defined as
\begin{equation}\label{Rij}
K(t,C;k_a,k_d) = \left\{
\begin{array}{>{\displaystyle}l>{\displaystyle}l}
& 0, \qquad t\leq t_0 + \Delta t, \\
& \frac{k_{a} C}{k_{d}+ k_{a} C} \left( 1- e^{-(k_{d}+ k_{a} C)(t-t_0)} \right), \\
& \qquad\qquad\qquad\qquad\qquad t_0+ \Delta t < t \leq t_0+ t_{inj}+ \Delta t, \\ &
\frac{k_{a} C}{k_{d}+ k_{a} C} \left( 1- e^{-(k_{d}+ k_{a} C) t_{inj}} \right) e^{-k_{d}(t-t_0-t_{inj})} , \\
& \qquad\qquad\qquad\qquad\qquad\qquad\qquad  t> t_0+ t_{inj}+ \Delta t,
\end{array}\right.
\end{equation}
where the physical meaning of parameters $t_0, t_{inj}, \Delta t$ can be found in Appendix A. All of them are given constants.

Since the solution to the integral equation (\ref{IntegralEq}) does not have an analytically closed form, we have to solve (\ref{IntegralEq}) numerically. Moreover, in order to represent the rate constant map on computer, we need to discretize the function $f(k_a,k_d)$ by a vector $\mathbf{c}$, associated with appropriate coordinates. By slightly abusing the notation, we shall call $\mathbf{c}$ the rate constant map as well. In this work, we adapt the finite element technique for discretization. The relation between $f$ and its finite element approximation is provided in the supplementary material. We discretize our bounded domain $\Omega \subset \mathbb{R}^2$ by mesh $\mathcal{T}$ using non-overlapping triangles $\{\bigtriangleup_\mu\}^M_{\mu=1}$ with the standard shape regularity assumption. Introduce the finite element space $V_n$ of piecewise linear elements on the triangulation $\mathcal{T}$ with the number of notes $n$. Then, the projection of $f(k_a,k_d)$ in $V_n$ (denote as $f_n$) can be decomposed as
\begin{equation}
f_n(k_a, k_d) = \sum^n_{l=1} \mathbf{c}_l \cdot \varphi_l(k_a, k_d),
\label{decomposition_f}
\end{equation}
where $\{\mathbf{c}_l\}$ are the coordinates of $f_n$ under basis $\{\varphi_l\}$ of the linear space $V_n$.

\begin{remark}\label{RemarkAdaptive}
The dimensionality $n$ of approximate solution space cannot be fixed, since we have no a priori information about the distribution of peaks of the rate constant map $f$, and cannot provide a good mesh $\mathcal{T}_n$ which nodes hit or close to the positions of the peaks of $f$. A naive way to capture the locations of the peaks of $f$ is to use an intensive grid with large $n$ so that every possible peak is located in a small neighborhood of a node of mesh. However, such an idea cannot be implemented in practice due to computational complexity. Alternatively, in this work, we start with a small $n$ paired with a coarse mesh $\mathcal{T}_n$, and then increase $n$ with a finer local mesh that only covers the region with possible peaks of $f$ until an approximate solution with satisfactory accuracy is obtained. Hence, instead of a fixed model (\ref{FiniteModel}) below, our algorithm will deal with a sequence of models of type (\ref{FiniteModel}) with dynamic dimensionality $n$.
\end{remark}

Furthermore, denote by $\mathbf{R}^j=[R_{obs}(t_1;C_j), R_{obs}(t_2;C_j), ..., R_{obs}(t_{N_T};C_j)]^T$ the measured
sensorgram with the injection $C_j$ and the time grid $\{t_i\}^{N_T}_{i=1}$. Then, using the decomposition (\ref{decomposition_f}) and considering the inaccuracy of the measurement data, the finite element approximation of integral equation (\ref{IntegralEq}) can be written as a system of algebraic equations:
\begin{equation}\label{FiniteModel0}
\mathbf{R}^j = \mathbf{K}^j  \mathbf{c} + \mathbf{\epsilon}^j, \quad j=1, ..., N_C,
\end{equation}
where $\mathbf{c}=[\mathbf{c}_1, \mathbf{c}_2, ..., \mathbf{c}_n]^T$, $\mathbf{\epsilon}^j\in \mathbb{R}^{N_T}$ is the additive error for $j$-th sensorgram, and matrix $\mathbf{K}^j$ is defined by $[\mathbf{K}^j]_{i,l} = \int_{\Omega} K(t_i,C_j;k_a,k_d) \varphi_l(k_a, k_d) d k_a d k_d.$

Denote $\mathbf{R}=[\mathbf{R}^1; ...; \mathbf{R}^{N_C}]$, $\mathbf{K}=[\mathbf{K}^1; ...; \mathbf{K}^{N_C}]$ and $\mathbf{\epsilon}=[\mathbf{\epsilon}^1; ...; \mathbf{\epsilon}^{N_C}]$, where for column vectors $\mathbf{a}$ and $\mathbf{b}$, $[\mathbf{a};\mathbf{b}]=[\mathbf{a}^T,\mathbf{b}^T]^T$. Then, the system (\ref{FiniteModel0}) can be rewritten into a compact form
\begin{equation}\label{FiniteModel}
\mathbf{R} = \mathbf{K}  \mathbf{c} + \mathbf{\epsilon}.
\end{equation}

The deterministic approaches solve (\ref{FiniteModel}) by the iterative regularization (e.g., the Landweber method \ctp{engl1996regularization}), or the variational regularization with appropriate regularization term (e.g., the Tikhonov method), to obtain a regularized approximate solution $\mathbf{c}_{reg}$ of (\ref{FiniteModel}). Using the obtained $\mathbf{c}_{reg}$, an approximate rate constant map $f_n$ can be reconstructed by (\ref{decomposition_f}), which is a function well defined in the whole domain $\Omega$. The active association and dissociation constants ($k_{a}$ and $k_{d}$) for the given kinetic reaction process of the biomolecular system should be the peaks in the estimated rate constant map $f_n$, see \ctp{Svitel-2003} and \ctp{Altschuh-2012}. In this work, we refer to these peaks as the Intrinsic Property of the Interaction (IPoI). Now, the problem at hand is how to determine the IPoI in a chemical reaction through the single step model (\ref{IntegralEq}) (or its finite dimensional analogue (\ref{FiniteModel})). This problem is seriously ill-posed since (i) one can always construct a large number of (deterministic) rate constant map $f$ satisfying the equation (\ref{IntegralEq}). In the reduced finite dimensional model (\ref{FiniteModel}), various numerical experiments have shown that $\mbox{rank}(\mathbf{K}_{N_T N_C \times n})\ll \min\{N_T N_C, n\}$. (ii) Note that the integral operator in (\ref{IntegralEq}) is a compact operator, which is affected by an instability phenomenon where a small amount of noise in the measurement data can lead to enormous errors in the estimates \citep{Tikhonov-1977}. Moreover, solutions to the deterministic models, e.g., the regularized solutions in \ctp{Svitel-2003}, \ctp{Altschuh-2012}, \ctp{Gorshkova-2008} and \ctp{AIDA}, do not represent the pure information of IPoI even though the experimental data is accurate enough. This phenomenon implies that the deterministic regularized solutions contain redundant information about a chemical reaction besides the IPoI. This drawback, as well as the problem of choosing a regularization term in the deterministic models aforementioned, motivates us to model the quantity $\mathbf{c}$ in (\ref{FiniteModel}) (or $f$ in (\ref{IntegralEq})) as a random vector (or function), from which we can study the uncertainty of the rate constant. As demonstrated in Section 5, the uncertainty of the quantity $\mathbf{c}$ can be used to interpret the difference between the deterministic solution of (\ref{FiniteModel}) and the IPoI.

\section{An adaptive variational Bayesian approach}

\subsection{Bayesian framework}

A detailed understanding of the dynamics of the rate constant map $\mathbf{c}$ given data $\mathbf{R}$ and the model (\ref{FiniteModel}) is crucial in the context of our problem.

For simplicity, a Gaussian distribution of noise on each sensorgram is customarily assumed, i.e., $\mathbf{\epsilon}^j$ ($j=1, ..., N_C$) are independent and identically distributed (i.i.d.) additive Gaussian random vectors with mean zero and variance $\sigma^2_j I_{N_T}$, where $I_{N_T}$ denotes the identity matrix of size $N_T$. The assumption of normality is preferred due to the fact that Gaussian noise is present at the sensorgrams, and sensorgram noise is usually much smaller than the signal. Although correlated sensorgram noise is more realistic than ``homogenerous'' sensorgrams, this could complicate the problem. Then, the likelihood $p(\mathbf{R}|\mathbf{c})$ is given by
\begin{equation}
p(\mathbf{R}|\mathbf{c}) \propto \exp \left\{ -\frac{1}{2} (\mathbf{R} - \mathbf{K}  \mathbf{c})^T \Sigma^{-1} (\mathbf{R} - \mathbf{K}  \mathbf{c}) \right\} \prod\limits^{N_C}\limits_{j=1} \left\{ \sigma^{-N_T}_j \right\} ,
\end{equation}
where $\varpropto$ denotes up to a multiplicative normalizing constant, and the covariance matrix of noise $\Sigma=\mbox{diag}(\sigma^2_1 I_{N_T}, ..., \sigma^2_{N_C} I_{N_T})$.

A natural mechanism for regularization is in the form of {\it prior} information. A versatile prior distribution of $\mathbf{c}$ is the Markov random field prior; see~\ctp{Babacan-2008} for details. Define $\mathbf{1}_{\Theta}$ as the indicator function on the set $\Theta$, i.e. $\mathbf{1}_{\Theta}(\mathbf{c})=1$ for $\mathbf{c}\in \Theta$, while $\mathbf{1}_{\Theta}(\mathbf{c})=0$ for $\mathbf{c}\not\in \Theta$. Denote by $\mathbf{c}\geq0$ the nonnegativity by components. For computational tractability, we use the following conjugate prior distribution
\begin{equation}\label{prior}
p(\mathbf{c}) \propto \sigma^{-n}_{\mathbf{c}} \exp \left\{ -\frac{1}{2 \sigma^2_{\mathbf{c}}} (\mathbf{L}  \mathbf{c})^T (\mathbf{L}  \mathbf{c}) \right\} \mathbf{1}_{\mathbf{c}\geq0},
\end{equation}
where the matrix $\mathbf{L}$ encapsulates the structure of the interactions between neighboring components of the solution $\mathbf{c}$. $\mathbf{L}$ is given and assumed to be full rank. Note that our choice of \eqref{prior} also reflects the physical meaning of the rate constant.  If the hyper-parameters, $\{\sigma_j\}^{N_C}_{j=1}$ and $\sigma^2_{\mathbf{c}}$, are given and the nonnegative constraint of the solution $\mathbf{c}$ is ignored, the prior of $\mathbf{c}$ in (\ref{prior}) for the \textit{maximum a posteriori} estimate is exactly the same as the $L_2$-type penalty term $\|\mathbf{L}  \mathbf{c}\|^2$ in the classical variational regularization theory \citep{Stuart-2010}, where $\|\cdot\|$ denotes the standard Euclidean norm of a vector or the Frobenius norm of a matrix. For the classical Tikhonov regularization, $\mathbf{L}$ stands for the discrete Laplacian.

To complete the model specification, we impose conjugate priors on the parameters $\{\sigma^2_j\}^{N_C}_{j=1}$ and $\sigma^2_{\mathbf{c}}$. In other words, for each of these parameters, we impose an independent inverse-gamma prior with the density $IG(x;\alpha,\beta) = \frac{\beta^\alpha}{\Gamma(\alpha)} x^{-\alpha -1} \exp \left( - \frac{\beta}{x} \right)$, where $\alpha$ is the shape parameter, $\beta$ is the scale parameter, and $\Gamma (\cdot)$ denotes the standard Gamma function. Then, plugging-in the density functions in (\ref{FiniteModel}), the posterior distribution of  $\mathbf{c}, \sigma^2_1, \sigma^2_2, ..., \sigma^2_{N_C}, \sigma^2_{\mathbf{c}}$ given the data $\mathbf{R}$ is
\begin{equation}\label{ppdfEnd}
\begin{array}{>{\displaystyle}l>{\displaystyle}l}
& p(\mathbf{c}, \sigma^2_1, \sigma^2_2, ..., \sigma^2_{N_C}, \sigma^2_{\mathbf{c}}| \mathbf{R})  \varpropto p( \mathbf{R} | \mathbf{c} , \sigma^2_{\mathbf{c}} ) \cdot p( \mathbf{c} | \sigma^2_1, \sigma^2_2, ..., \sigma^2_{N_C}) \cdot p(\sigma^2_{\mathbf{c}} ) \cdot \prod\limits^{N_C}\limits_{j=1} p(\sigma^2_j)\\
& \varpropto \sigma^{-n-2\alpha_{\mathbf{c}} -2}_{\mathbf{c}} \cdot \exp \left( - \frac{\beta_{\mathbf{c}}}{\sigma^2_{\mathbf{c}}} \right) \prod\limits^{N_C}\limits_{j=1} \left\{ \sigma^{-N_T -2\alpha_j -2}_j \exp \left( - \frac{\beta_j}{\sigma^2_j} \right) \right\} \\ &
\qquad \cdot \exp \left\{ -\frac{1}{2} (\mathbf{R} - \mathbf{K}  \mathbf{c})^T \Sigma^{-1} (\mathbf{R} - \mathbf{K}  \mathbf{c}) - \frac{1}{2 \sigma^2_{\mathbf{c}}} (\mathbf{L}  \mathbf{c})^T (\mathbf{L}  \mathbf{c}) \right\} \cdot \mathbf{1}_{\mathbf{c}\geq0},
\end{array}
\end{equation}
where the values of hyperparameters $(\alpha_{\mathbf{c}}, \beta_{\mathbf{c}})$ and $(\alpha_j,\beta_j)$ are given based on simulation results of artificial problems. The posterior distribution $p(\mathbf{c}, \sigma^2_1, \sigma^2_2,$ $..., \sigma^2_{N_C}, \sigma^2_{\mathbf{c}}| \mathbf{R})$ in (\ref{ppdfEnd}) is the full Bayesian solution to our finite model (\ref{FiniteModel}), and it encapsulates all the information including IPoI. However, due to the high dimensionality of $\mathbf{c}$, \eqref{ppdfEnd} is also a distribution that lives in a space of very high dimension, e.g., $10^3\sim 10^5$, for our real data. Hence, it is advisable to develop tools for exploring this very high dimensional posterior distribution. In the next section, we develop an approximate inference method based on the mean field variational approximation for exploring the posterior in \eqref{ppdfEnd}.

\subsection{Variational approximation algorithm}

Due to the presence of several variance parameters, the posterior distribution (\ref{ppdfEnd}) does not have an explicit closed form. One way to explore this high dimensional posterior is to perform posterior sampling based on MCMC algorithms. Although the MCMC has the advantage of being asymptotically exact \citep{Robert-2004}, its convergence is also often known to be difficult to diagnose \citep{Brooks-2002,Bishop2006}, and it suffers from slow mixing in the presence of high dimensional parameters. As the dimensionality of the problem in the biosensor system is extremely large, we shall take an alternative route and focus on the variational approximation approach. The essential idea of variational methods consists of first transforming the problem into an equivalent optimization problem, and then obtaining an approximate distribution $q(\cdot)$ to the true posterior $p(\cdot)$ by solving the optimization problem inexactly. The idea behind VB approximation is to find a family of tractable distributions $\{q^k(\cdot)\}_k$ to approximate the posterior density by minimizing the Kullback-Leibler (KL) divergence~\citep{Kullback1959}, while still capturing distinct features of the posterior distribution (\ref{ppdfEnd}) in a computationally efficient way.

The KL divergence is a non-symmetric measure of the difference between two probability distributions, and is defined as
\begin{equation}\label{KL}
\begin{array}{>{\displaystyle}l>{\displaystyle}l}
&\mathcal{D} \left( q(\mathbf{c}, \sigma^2_1, ..., \sigma^2_{N_C}, \sigma^2_{\mathbf{c}}) | p(\mathbf{c}, \sigma^2_1, ..., \sigma^2_{N_C}, \sigma^2_{\mathbf{c}}| \mathbf{R}) \right) = \\
&  \int q(\mathbf{c}, \sigma^2_1, ..., \sigma^2_{N_C}, \sigma^2_{\mathbf{c}}) \ln \left( \frac{q(\mathbf{c}, \sigma^2_1, ..., \sigma^2_{N_C}, \sigma^2_{\mathbf{c}})}{p(\mathbf{c}, \sigma^2_1, ..., \sigma^2_{N_C}, \sigma^2_{\mathbf{c}}| \mathbf{R})} \right) d \mathbf{c} d \sigma^2_1 \cdot\cdot\cdot d \sigma^2_{N_C} d \sigma^2_{\mathbf{c}} = \\
& \int q(\mathbf{c}, \sigma^2_1, ..., \sigma^2_{N_C}, \sigma^2_{\mathbf{c}}) \ln \left( \frac{q(\mathbf{c}, \sigma^2_1, ..., \sigma^2_{N_C}, \sigma^2_{\mathbf{c}})}{p(\mathbf{c}, \sigma^2_1, ..., \sigma^2_{N_C}, \sigma^2_{\mathbf{c}}, \mathbf{R})} \right) d \mathbf{c} d \sigma^2_1 \cdot\cdot\cdot d \sigma^2_{N_C} d \sigma^2_{\mathbf{c}} \\
& \qquad  + \ln(p(\mathbf{R})) := {\rm ELBO} \left( q(\mathbf{c}, \sigma^2_1, ..., \sigma^2_{N_C}, \sigma^2_{\mathbf{c}}) | p(\bullet,\mathbf{R}) \right) + \ln(p(\mathbf{R})),
\end{array}
\end{equation}
where $p(\mathbf{R})= \int p(\mathbf{c}, \sigma^2_1, ..., \sigma^2_{N_C}, \sigma^2_{\mathbf{c}}, \mathbf{R}) d \mathbf{c} d \sigma^2_1 \cdot\cdot\cdot d \sigma^2_{N_C} d \sigma^2_{\mathbf{c}}$ is a normalizing constant, and the global distribution $p(\mathbf{c}, \sigma^2_1, ..., \sigma^2_{N_C}, \sigma^2_{\mathbf{c}}, \mathbf{R})$ has the form
\begin{equation}\label{globalDistribution}
p(\mathbf{c}, \sigma^2_1, ..., \sigma^2_{N_C}, \sigma^2_{\mathbf{c}}, \mathbf{R}) = p(\mathbf{c}|  \sigma^2_{\mathbf{c}}) \cdot p( \mathbf{R}| \mathbf{c}, \sigma^2_1, ..., \sigma^2_{N_C}) \cdot p(\sigma^2_{\mathbf{c}}) \cdot \prod\limits^{N_C}\limits_{j=1}  p(\sigma^2_j).
\end{equation}

Note that minimizing the KL distance in (\ref{KL}) is equivalent to minimizing the first term, which is termed as the \emph{evidence lower bound} and denoted by ELBO \citep{Blei-2017}. In this way, we have successfully transformed the sampling problem into an equivalent optimization problem of finding a simpler distribution $q(\mathbf{c}, \sigma^2_1, ..., \sigma^2_{N_C}, \sigma^2_{\mathbf{c}})$ by minimizing ${\rm ELBO}$.

If we impose no constraint on the approximation $q(\mathbf{c}, \sigma^2_1, ..., \sigma^2_{N_C}, \sigma^2_{\mathbf{c}})$, minimizing ${\rm ELBO}$ is numerically intractable. The intractability is largely due to the correlations among the parameters. To enable the computational tractability, we adapt the idea of \emph{mean-field variational family} \citep{Wand2011,Blei-2017}, and impose an independence condition among the parameter components $\mathbf{c}, \sigma^2_1, ..., \sigma^2_{N_C}$ and $\sigma^2_{\mathbf{c}}$ as
\begin{equation}\label{independence}
q(\mathbf{c}, \sigma^2_1, ..., \sigma^2_{N_C}, \sigma^2_{\mathbf{c}}) = q(\mathbf{c}) q(\sigma^2_1) \cdot\cdot\cdot q(\sigma^2_{N_C}) q(\sigma^2_{\mathbf{c}}).
\end{equation}
Under assumption (\ref{independence}), we can find an effective approximate posterior density, denoted by $q(\mathbf{c}, \sigma^2_1, ...,$ $\sigma^2_{N_C}, \sigma^2_{\mathbf{c}})$ hereafter, using an alternating direction iterative algorithm, i.e., at $k$-th step, we solve the following optimization problems:
\begin{equation}\label{Algorithm1}
\left\{\begin{array}{>{\displaystyle}l>{\displaystyle}l}
& \min_{q(\mathbf{c})} {\rm ELBO}\left( q(\mathbf{c})q^{k-1}(\sigma^2_1) \cdot\cdot\cdot q^{k-1}(\sigma^2_{N_C}) q^{k-1}(\sigma^2_{\mathbf{c}}) | p(\bullet,\mathbf{R}) \right),  \\
& \min_{q(\sigma^{2}_j)} {\rm ELBO} \Big( q^k(\mathbf{c}) q^{k}(\sigma^2_1) \cdot\cdot\cdot q^{k}(\sigma^2_{j-1}) q^{k-1}(\sigma^2_{j+1}) \cdot\cdot\cdot q^{k-1}(\sigma^2_{N_C}) q^{k-1}(\sigma^2_{\mathbf{c}}) \\
& \qquad\qquad\qquad\qquad  | p(\bullet,\mathbf{R}) \Big), \qquad j=1, ..., N_C, \\
& \min_{q(\sigma^{2}_{\mathbf{c}})} {\rm ELBO}\left( q^k(\mathbf{c})q^k(\sigma^2_1) \cdot\cdot\cdot q^k(\sigma^2_{N_C}) q(\sigma^2_{\mathbf{c}}) | p(\bullet,\mathbf{R}) \right).
\end{array}\right.
\end{equation}

Here and later on, the (superscript and subscript) index $k$ denotes the quantity of interest at $k$-th iteration of the algorithm. Using condition (\ref{independence}), one can derive the explicit formulas for the minimizers of the optimization problems in (\ref{Algorithm1}) by examining the optimality system; i.e., the following theorem holds.

\begin{thm}\label{thm1}
Assume that the conditional independence condition (\ref{independence}) holds. Then, the minimizers of the optimization problems in Algorithm 1 at each iteration have the following explicit formulas:
\begin{equation}\label{minimizer_c}
\left\{\begin{array}{>{\displaystyle}r>{\displaystyle}l}
q^k(\mathbf{c}) &= \mathcal{N}_+(\mathbf{c}_k,\Sigma^{\mathbf{c}}_{k}),
\mathbf{c}_k = \Sigma^{\mathbf{c}}_{k} \mathbf{K}^T \Sigma^{-1}_{k} \mathbf{R},
\Sigma^{\mathbf{c}}_{k} =  \left( \mathbf{K}^T \Sigma^{-1}_{k} \mathbf{K} + \sigma^{-2}_{\mathbf{c},k} \mathbf{L}^T \mathbf{L} \right)^{-1}, \\
\Sigma^{-1}_{k} &= \mathbb{E}_{q^{k-1}(\sigma^2_1) \cdot\cdot\cdot q^{k-1}(\sigma^2_{N_C}) q^{k-1}(\sigma^2_{\mathbf{c}})} \left[ \Sigma^{-1} \right], \\
\sigma^{-2}_{\mathbf{c},k} &= \mathbb{E}_{q^{k-1}(\sigma^2_1) \cdot\cdot\cdot q^{k-1}(\sigma^2_{N_C}) q^{k-1}(\sigma^2_{\mathbf{c}})} \left[ \sigma^{-2}_{\mathbf{c}} \right], \\
q^k(\sigma^2_j) &= IG\left(\sigma^2_j; \alpha_j+\frac{N_T}{2}, \beta_j + \frac{1}{2} \mathbb{E}_{q^k(\mathbf{c})} \left[ (\mathbf{R}^j - \mathbf{K}^j  \mathbf{c})^T (\mathbf{R}^j - \mathbf{K}^j  \mathbf{c}) \right] \right), \\
q^k(\sigma^2_{\mathbf{c}}) &= IG\left(\sigma^2_{\mathbf{c}}; \alpha_{\mathbf{c}} + \frac{n}{2}, \beta_{\mathbf{c}} + \frac{1}{2} \mathbb{E}_{q^k(\mathbf{c})} \left[ (\mathbf{L} \mathbf{c})^T (\mathbf{L} \mathbf{c}) \right] \right),
\end{array}\right.
\end{equation}
where the density function for truncated normal distribution $\mathcal{N}_+(\mathbf{c}^*,\Sigma^{\mathbf{c}}_{*})$ is
\begin{equation}\label{ModifiedGaussian}
\frac{1}{(1-HN(0)) \sqrt{(2\pi)^n |\Sigma^{\mathbf{c}}_{*}|}} \exp \left\{ -\frac{1}{2} (\mathbf{c} - \mathbf{c}^*)^T (\Sigma^{\mathbf{c}}_{*})^{-1} (\mathbf{c} - \mathbf{c}^*) \right\} \mathbf{1}_{\mathbf{c}\geq0},
\end{equation}
$HN(0)=\mathbb{P}(\xi\leq0)$ with $\xi\sim \mathcal{N}(\mathbf{c}^*,\Sigma^{\mathbf{c}}_{*})$ and $\mathbb{E}_q [\bullet]$ is the expectation with respect to the density $q$.
\end{thm}

In practice, the value of $HN(0)$ can be estimated numerically. 
The proof of Theorem \ref{thm1} follows a standard argument in theory of variational inference, and we provide a sketch in the supplementary material \citep{ZhangYaoSupplementary}.

Now, we discuss the convergence issue of scheme (\ref{Algorithm1}).

\begin{thm}
The sequence $\left\{ q^k(\mathbf{c}) q^k(\sigma^2_1) \cdot\cdot\cdot q^k(\sigma^2_{N_C}) q^k(\sigma^2_{\mathbf{c}}) \right\}$ generated by scheme (\ref{Algorithm1}) converges, upon a subsequence, to a stationary point $q^*(\mathbf{c}) q^*(\sigma^2_1) \cdot\cdot\cdot q^*(\sigma^2_{N_C}) q^*(\sigma^2_{\mathbf{c}})$ of the KL distance functional $\mathcal{D}$, which satisfies
\begin{equation}\label{result}
\left\{\begin{array}{>{\displaystyle}r>{\displaystyle}l}
q^*(\mathbf{c}) &= \mathcal{N}_+(\mathbf{c}^*,\Sigma^{\mathbf{c}}_{*}), \\
q^*(\sigma^2_j) &= IG\left(\sigma^2_j; \alpha_j+\frac{N_T}{2}, \beta_j + \frac{1}{2} \mathbb{E}_{q^*(\mathbf{c})} \left[ (\mathbf{R}^j - \mathbf{K}^j  \mathbf{c})^T (\mathbf{R}^j - \mathbf{K}^j  \mathbf{c}) \right] \right), \\
q^*(\sigma^2_{\mathbf{c}}) &= IG\left(\sigma^2_{\mathbf{c}}; \alpha_{\mathbf{c}} + \frac{n}{2}, \beta_{\mathbf{c}} + \frac{1}{2} \mathbb{E}_{q^*(\mathbf{c})} \left[ (\mathbf{L} \mathbf{c})^T (\mathbf{L} \mathbf{c}) \right] \right), \\
\mathbf{c}^* &= \Sigma^{\mathbf{c}}_{*} \mathbf{K}^T \Sigma^{-1}_{*} \mathbf{R}, ~
\Sigma^{\mathbf{c}}_{*} =  \left( \mathbf{K}^T \Sigma^{-1}_{*} \mathbf{K} + \sigma^{-2}_{\mathbf{c},*} \mathbf{L}^T \mathbf{L} \right)^{-1}, \\
\Sigma^{-1}_{*} &= \mathbb{E}_{q^*(\sigma^2_1) \cdot\cdot\cdot q^*(\sigma^2_{N_C}) q^*(\sigma^2_{\mathbf{c}})} \left[ \Sigma^{-1} \right], ~
\sigma^{-2}_{\mathbf{c},*} = \mathbb{E}_{q^*(\sigma^2_1) \cdot\cdot\cdot q^*(\sigma^2_{N_C}) q^*(\sigma^2_{\mathbf{c}})} \left[ \sigma^{-2}_{\mathbf{c}} \right].
\end{array}\right.
\end{equation}
\end{thm}

We refer readers to Appendix B for the proof of Theorem 2.  To that end, let us consider the stopping principle of our scheme. Various stopping criteria exist for an iterative algorithm, e.g., ELBO is usually calculated (whenever possible) as a stopping criterion in the variational inference. 
For our problem, a natural and simple stopping criterion is to utilize the mean and covariance of the interested quantity~-- the rate constant map $\mathbf{c}$. Hence, we define the accuracy of our algorithm as
\begin{equation}\label{accuracyAlgorithm}
\Delta_1(q^k(\mathbf{c})) = \frac{\|\mathbf{c}_{k}- \mathbf{c}_{k-1}\|}{\|\mathbf{c}_{k-1}\|}, \quad \Delta_2(q^k(\mathbf{c})) = \frac{\|\Sigma^{\mathbf{c}}_{*,k}- \Sigma^{\mathbf{c}}_{*,k-1}\|}{\|\Sigma^{\mathbf{c}}_{*,k-1}\|}.
\end{equation}

Finally, by coupling the stopping criteria and Theorem 1, we offer a variational Bayesian algorithm for reconstructing the rate constant map in the biosensor system in Algorithm 1.

\begin{algorithm}[H]
\label{Alorithm2}
\caption{A variational Bayesian algorithm for solving (\ref{FiniteModel}).}
\begin{algorithmic}[1]
\Require Model parameters $\{\alpha_j, \beta_j\}^{N_C}_{j=1}$ and $\alpha_{\mathbf{c}}$, $\beta_{\mathbf{c}}$. Initial guesses of $q^0(\mathbf{c})$. Tolerance $\varepsilon$.

\Ensure The approximate posterior distribution of $p(\mathbf{c})$ is $q^k(\mathbf{c})$.

\State $k \gets 1$

\While{$\Delta_i(q^k(\mathbf{c}))> \varepsilon$}

\State Update the distributions of hyperparameters $\sigma^2_{j,k}$ and $\sigma^2_{\mathbf{c},k}$ by formula (\ref{result}), replacing $\mathbf{c}$ by $\mathbf{c}_k$.

\State Update the distribution of solution $\mathbf{c}_k$ by formula (\ref{minimizer_c}).

\State $k \gets k+1$
\EndWhile
\end{algorithmic}
\end{algorithm}

\subsection{An adaptive strategy to improve the quality of the rate constant map}

As mentioned in Remark \ref{RemarkAdaptive}, we start with a low dimensionality $n$ of solution space where the rate constant map $f$ is only estimated on a set of coarse distributed points in $\Omega$. Therefore, the resolution of the rate constant map is quite low. In order to improve the quality of the estimated rate constant map $f_n$, we adopt the idea of the oriented adaptive discretization, from numerical technique in finite element methods \citep{Chen-2009,AIDA}. The idea, combined with our approach, is to first solve the problem (\ref{IntegralEq}) through the variational Bayesian approach (Algorithm 1), obtaining the solution (the distribution of the random vector $\mathbf{c}$) to the current triangulation. The quality is then estimated using the solution, and is used to mark a set of triangles to be refined. Triangles are refined in a way that maintains two of the most important properties of the triangulations: shape regularity and conformity.

Given an estimated distribution $\mathcal{N}_+ (\mathbf{c}^*,\Sigma^{\mathbf{c}}_{*})$, let $\widehat{\mathbf{c}}, \underline{\mathbf{c}}$ and $\overline{\mathbf{c}}$ be the sample mean, the lower and upper endpoints of the 95\% confidence interval by samples from Algorithm 1. Define by $V^{mean}_\mu [\widehat{\mathbf{c}}]$ the variation of the mean vector $\widehat{\mathbf{c}}$ over the elements $\bigtriangleup_\mu$, i.e. $V^{mean}_\mu [\widehat{\mathbf{c}}]= \frac{1}{3} \sum_{\imath\neq \jmath, \imath, \jmath\in\{1, 2, 3\}} \left| \widehat{\mathbf{c}}_{\mu_\imath} -  \widehat{\mathbf{c}}_{\mu_\jmath} \right|$. Similarly, we can define the variation of $\underline{\mathbf{c}}$ and $\overline{\mathbf{c}}$, denoted by $V^{lower}_\mu [\underline{\mathbf{c}}]$ and $V^{upper}_\mu [\overline{\mathbf{c}}]$ respectively. Finally, define the refinement indicator $V_\mu$ by
\begin{equation}\label{Indicator}
V_\mu [\widehat{\mathbf{c}}, \underline{\mathbf{c}}, \overline{\mathbf{c}}] = \max \left\{ V^{mean}_\mu [\widehat{\mathbf{c}}] ,  V^{lower}_\mu [\underline{\mathbf{c}}], V^{upper}_\mu [\overline{\mathbf{c}}] \right\}.
\end{equation}
Then, the refinement should be done in the triangles of all points in the finite element mesh $\mathcal{T}$ where the function $V_\mu$ achieves its maximum; i.e., refine the mesh in such triangles of $\Omega$ where
\begin{equation}
V_\mu [\widehat{\mathbf{c}}, \underline{\mathbf{c}}, \overline{\mathbf{c}}] \geq \tau \max_{\bigtriangleup_\mu\in \mathcal{T}} V_\mu [\widehat{\mathbf{c}}, \underline{\mathbf{c}}, \overline{\mathbf{c}}],
\label{refinement1}
\end{equation}
where $\tau\in (0,1)$ are numbers which should be chosen computationally.

Finally, we use the longest edge refinement rule for dividing the marked triangles, so that the mesh obtained by this dividing rule still conforms and has a shape regular~\citep{Rivara-1984}. Note that apart from the marked triangles, the additional triangles are also refined to recover the conformity of triangulations. In this work, we also control the number of elements added to ensure the overall optimality of the refinement procedure.

Let us present our main algorithm for solving the two-dimensional Fredholm integral equation of the type (\ref{IntegralEq}) in Algorithm 2.

\begin{algorithm}[!htb]
\label{Alorithm3}
\caption{An adaptive variational Bayesian approach (AVBA) for reconstructing the rate constant map in the biosensor system.}
\begin{algorithmic}[1]
\Require Initial mesh: $\mathcal{T}_1$. Tolerance number: $\tau$. Tolerance error: $\varepsilon$. Sample number: $N$. The maximum iteration number $K_{max}$.

\Ensure The estimated rate constant mean map: $\widehat{\mathbf{c}}=\widehat{\mathbf{c}}_k$. The lower and upper rate constant maps: $\underline{\mathbf{c}}=\underline{\mathbf{c}}_k$ and $\overline{\mathbf{c}}=\overline{\mathbf{c}}_k$.

\State $k \gets 1$

\While{$\Delta_i(q^k(\mathbf{c}))> \varepsilon$ \textbf{and} $k < K_{max}$}

\State Obtain the distribution $\mathcal{N}_+(\mathbf{c}_k,\Sigma^{\mathbf{c}}_{k})$ of random vector $\mathbf{c}$ of the finite element solution of problem (\ref{IntegralEq}) on the mesh $\mathcal{T}_k$ by Algorithm 1.

\State Compute the sample mean $\widehat{\mathbf{c}}_k$, the lower confidence limit $\underline{\mathbf{c}}_k$ and the upper confidence limit $\overline{\mathbf{c}}_k$ using samples for $\mathcal{N}_+(\mathbf{c}_k,\Sigma^{\mathbf{c}}_{k})$ from Algorithm 1.

\State Refine the mesh $\mathcal{T}_k$ at all points where $V_\mu \geq \tau \max_{\bigtriangleup_\mu\in \mathcal{T}} V_\mu$.

\State Construct a new mesh $\mathcal{T}_{k+1}$.

\State $k \gets k+1$
\EndWhile
\end{algorithmic}
\end{algorithm}

\section{Simulation study}

\begin{figure}[!b]
\centering
\subfigure[]{
\includegraphics[clip, trim=0 3.9in 0 1in, width=2.2in]{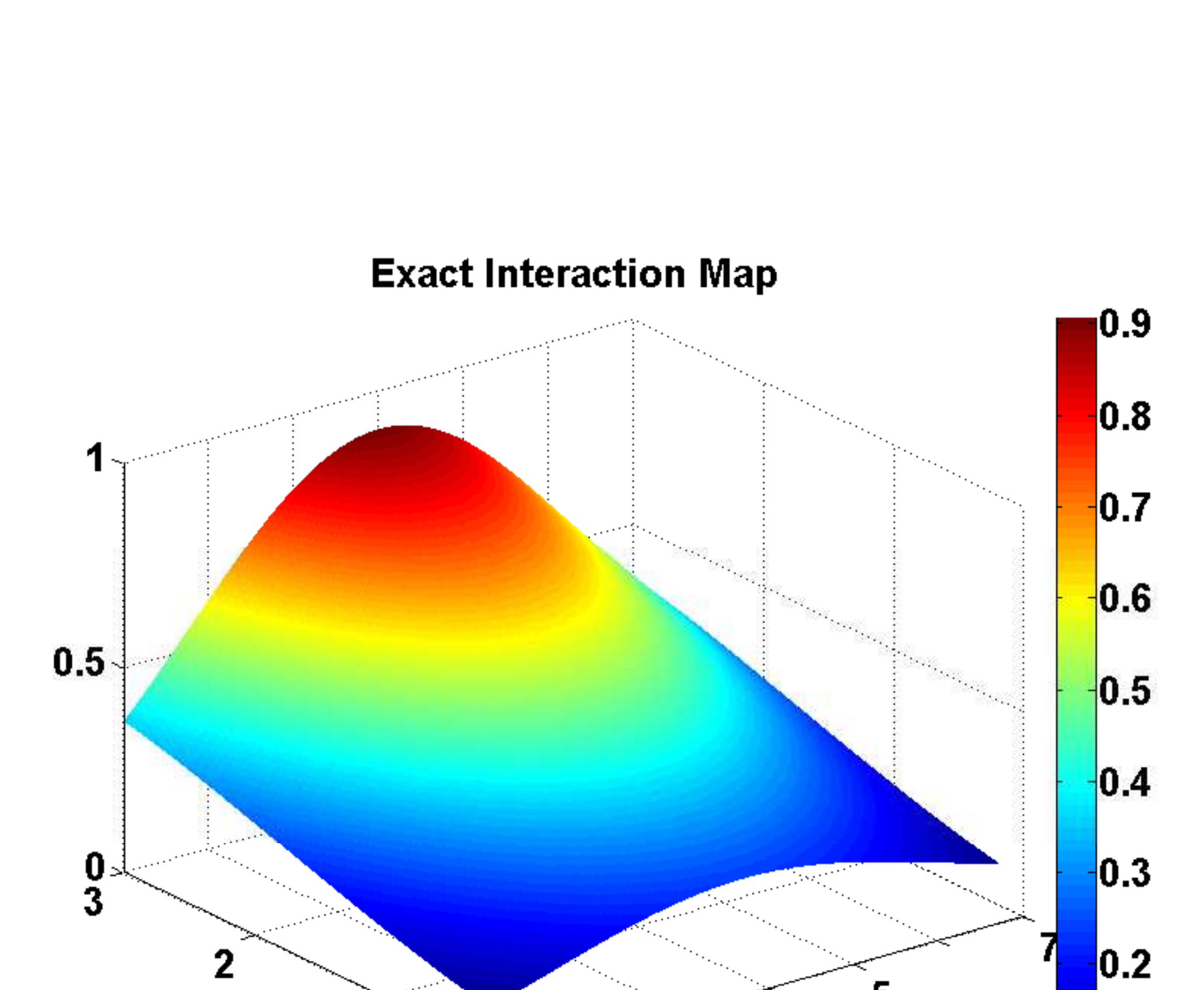}}
\subfigure[]{
\includegraphics[clip, trim=0 3.5in 0 1in, width=2.2in]{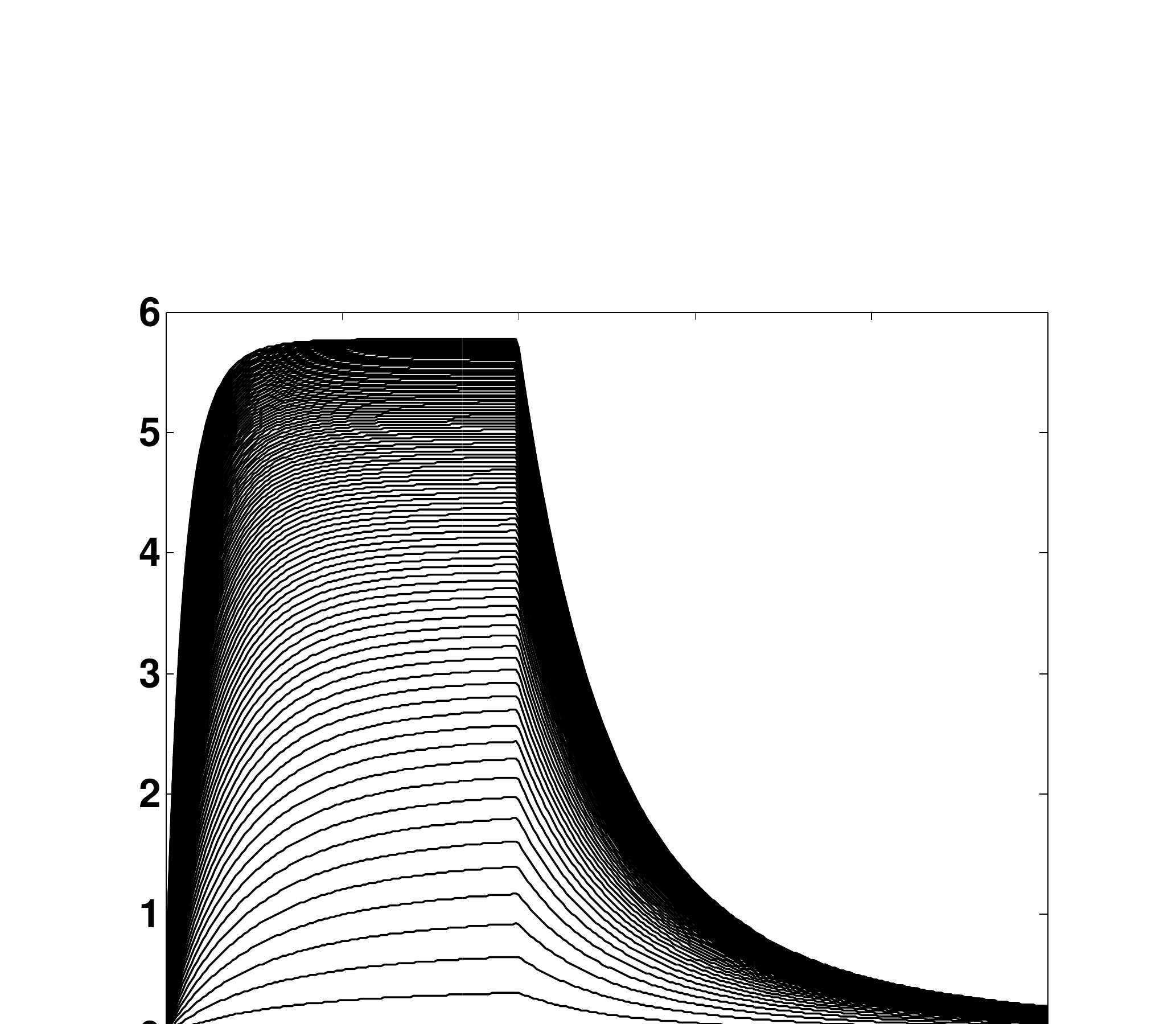}}
\caption[Exact rate constant map and accurate measurement]{(a) The exact rate constant map for model problem. (b) The synthetic time dependent noise-free data.}
\label{ExactSolution}
\end{figure}

This section presents the numerical results for the inverse problem (\ref{IntegralEq}), of which the adaptive variational Bayesian approach (AVBA), i.e., Algorithm 2, is used to illustrate its accuracy and flexibility. The simulation consists of three steps. First, a simulated response function $R_{obs}$ is generated by computer according to (\ref{IntegralEq}) for a given rate constant map $\bar{f}(k_a,k_d)=\exp\left( -0.1*((k_a-4)^2+(k_d-4)^2) \right)$ in the domain $\Omega=[1,7]\times[0,3]$~-- see (a) in Figure \ref{ExactSolution}. Denote by $\left\{ \bar{R}_{p,q} \right\}^{N_T, N_C}_{p,q=1}$ (see (b) in Figure \ref{ExactSolution}) the collected exact data at a uniform time grid $\left\{ t_p \right\}^{N_T}_{p=1}$ ($N_T=300$, $T=[0,4]$, $t_0=0$ and $t_{inj}=2$), at a uniform analyte concentration grid $\left\{ C_q \right\}^{N_C}_{q=1}$ ($N_C=100$ and $\Xi=[0.001,2]$). The synthetic noisy data is generated by
\begin{equation*}
R^\delta_{obs,p,q} = \bar{R}_{p,q} + \max\left\{ 0, \max_{p} \bar{R}_{p,q} \delta \omega_q \right\}, \quad p=1, ..., N_T,~ q=1, ..., N_C,
\end{equation*}
where $\delta$ indicates the noise level, and $\omega_q$ is a Gaussian random variable with zero mean and standard deviation $\sqrt{q}$.  At the last step of the simulation, the observation data $R^\delta_{obs}$ is processed through our algorithm, and the retrieved rate constant map $f(k_a,k_d)$ is compared with the input map. To assess the accuracy of the approximate solution, we define the $L_2$-norm relative error for an estimated random vector $\mathbf{c}$ with distribution $\mathcal{N}_+ (\widehat{\mathbf{c}}^*,\Sigma^{\mathbf{c}}_{*})$ as
\begin{equation*}\label{L2Err}
L2Err = \frac{\|\bar{f}(k_a,k_d) - \hat{f}(\widehat{\mathbf{c}}^*)\|_{L^2(\Omega)}}{\|\bar{f}(k_a,k_d)\|_{L^2(\Omega)}} , \quad \hat{f}(\widehat{\mathbf{c}}^*)= \sum^n_{l=1} \widehat{\mathbf{c}}^*_l \cdot \varphi_l(k_a, k_d).
\end{equation*}

A demonstration of Algorithm 2 can be found in the supplementary material. The fundamental assumption of our approach, according to (\ref{independence}), is the conditional independence of random variables $\{\mathbf{c}, \sigma^2_1, ..., \sigma^2_{N_C}, \sigma^2_{\mathbf{c}}\}$ given the data $\mathbf{R}$. To examine the assumption (\ref{independence}), we compute the Pearson correlation coefficients $\rho_{\mathbf{c},\sigma^2_j} (j=1, ..., N_C)$, $\rho_{\mathbf{c},\sigma^2_{\mathbf{c}}}$ and $\rho_{\sigma^2_i,\sigma^2_{\mathbf{c}}} (j=1, ..., N_C)$ between the vectors $\mathbf{c}$, $\sigma^2_{\mathbf{c}}$ and $\sigma^2_j$ respectively, from the MCMC samples. The correlation coefficient between the vector $\mathbf{c}$ and scalars are computed component-wise, i.e., with $\mathbf{c}_i$. The results, for the example in the supplementary material, are shown in Figure \ref{correlation}, where the abscissa $i$ (or $j$) denotes the $i$th (or $j$th) component. Overall, the correlation coefficients between $\sigma^2_{\mathbf{c}}$ and $\sigma^2_i$ are very small, with a maximum norm $|\rho_{\sigma^2_i,\sigma^2_{\mathbf{c}}}|_{\infty}:= \max\limits_{1\leq i \leq N_C} |\rho_{\sigma^2_i,\sigma^2_{\mathbf{c}}}|$, smaller than 0.004 for all three noise levels. The correlation coefficients between $\mathbf{c}$ and $\sigma^2_i$ are also small with a maximum norm $|\rho_{\mathbf{c},\sigma^2_j}|_{\infty}:= \max\limits_{1\leq i \leq n, 1\leq j \leq N_C} |\rho_{\mathbf{c}_i,\sigma^2_j}|$ smaller than 0.08 for all three noise levels. The correlation coefficients between $\mathbf{c}$ and $\sigma^2_{\mathbf{c}}$ are slightly larger with the maximum norm $|\rho_{\mathbf{c},\sigma^2_{\mathbf{c}}}|_{\infty}:= \max\limits_{1\leq i \leq n} |\rho_{\mathbf{c}_i,\sigma^2_{\mathbf{c}}}|$ close to 0.14. Hence, we can conclude that the correlation between $\mathbf{c}$, $\sigma^2_{\mathbf{c}}$ and $\sigma^2_j$ is relatively weak.

\begin{figure}[H]
\centering
\subfigure[]{
\includegraphics[clip, trim=0 2.8in 0 1in, width=1.4in]{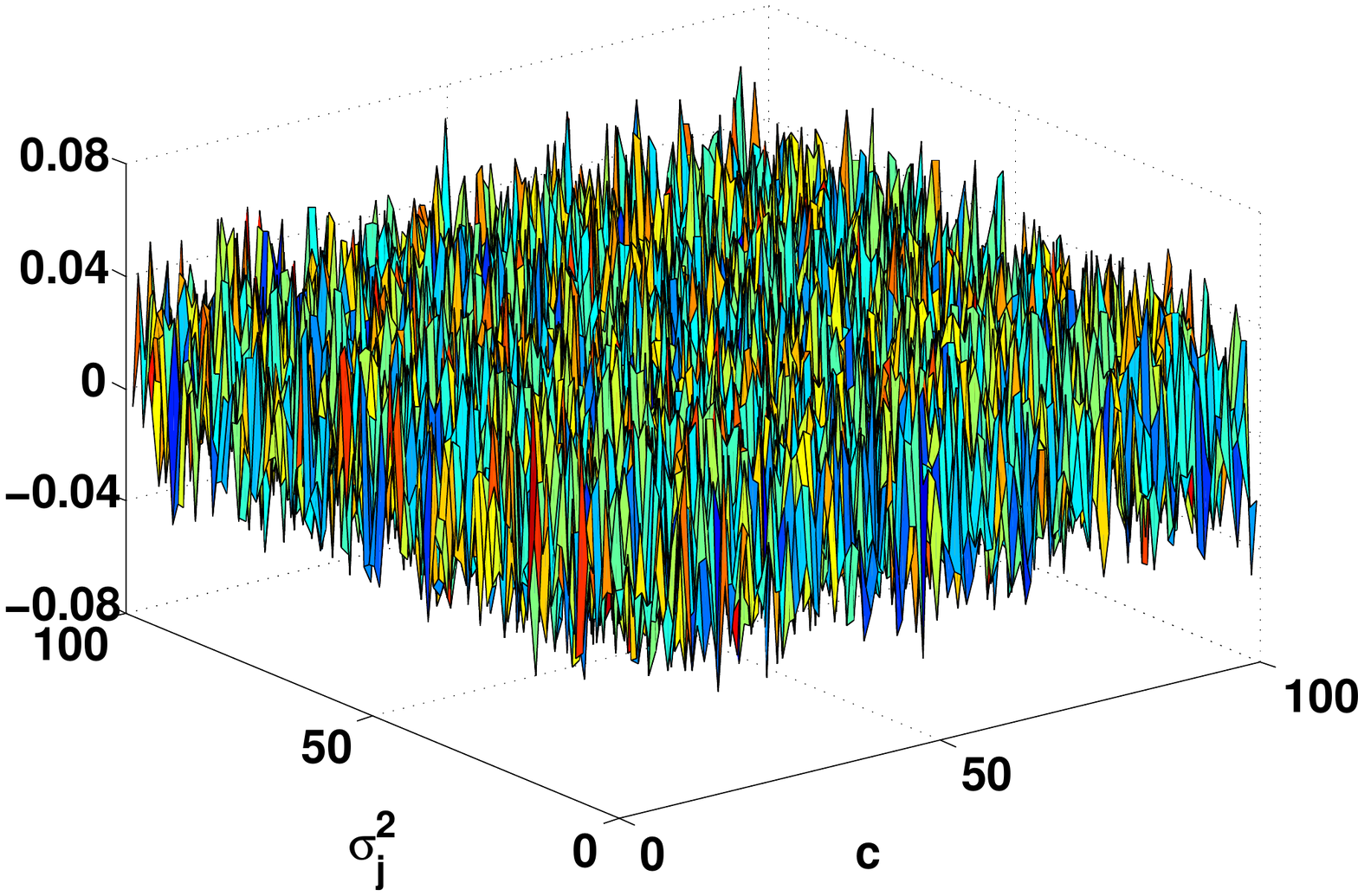}}
\subfigure[]{
\includegraphics[clip, trim=0 2.8in 0 1in, width=1.4in]{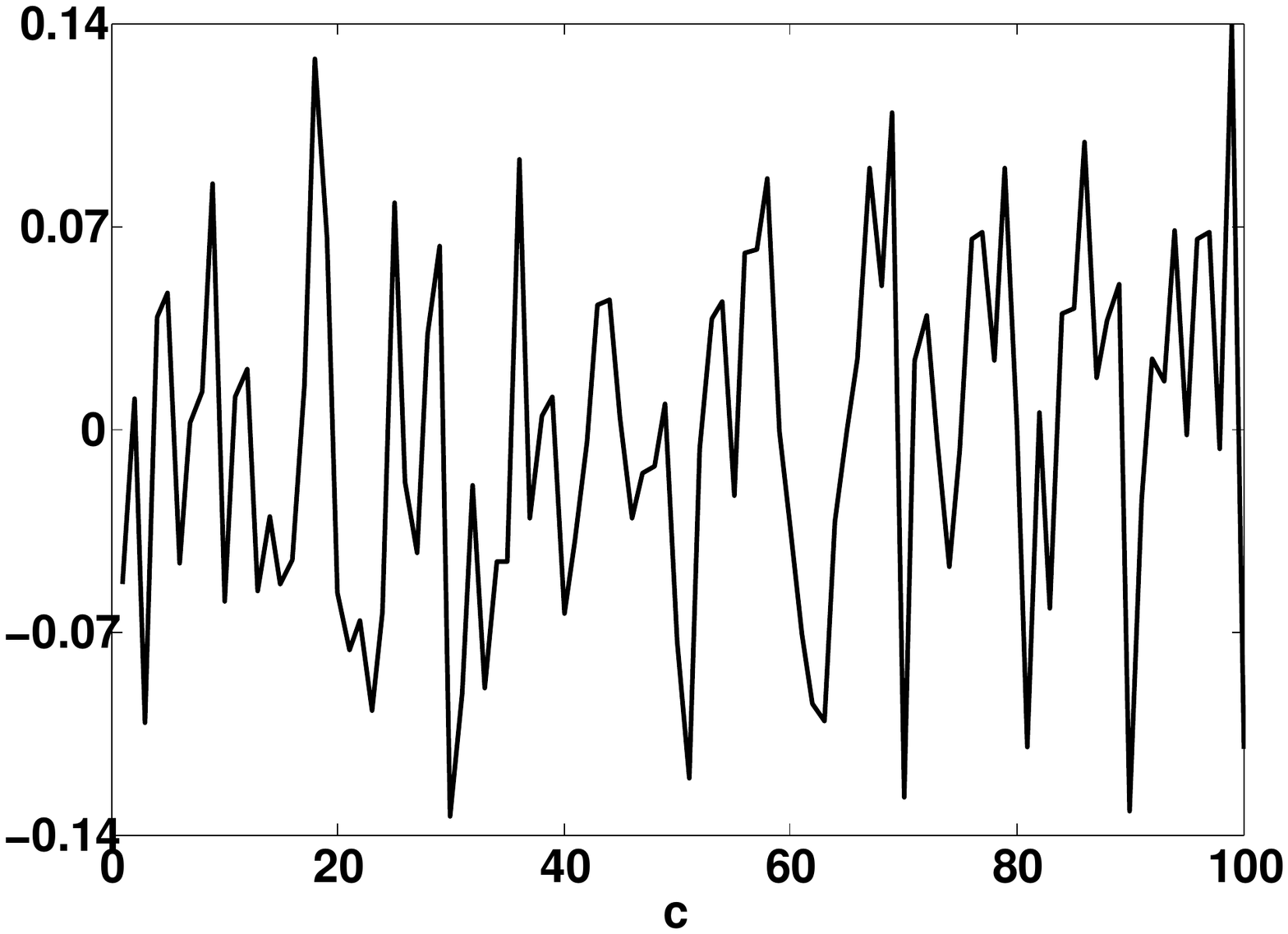}}
\subfigure[]{
\includegraphics[clip, trim=0 2.8in 0 1in, width=1.4in]{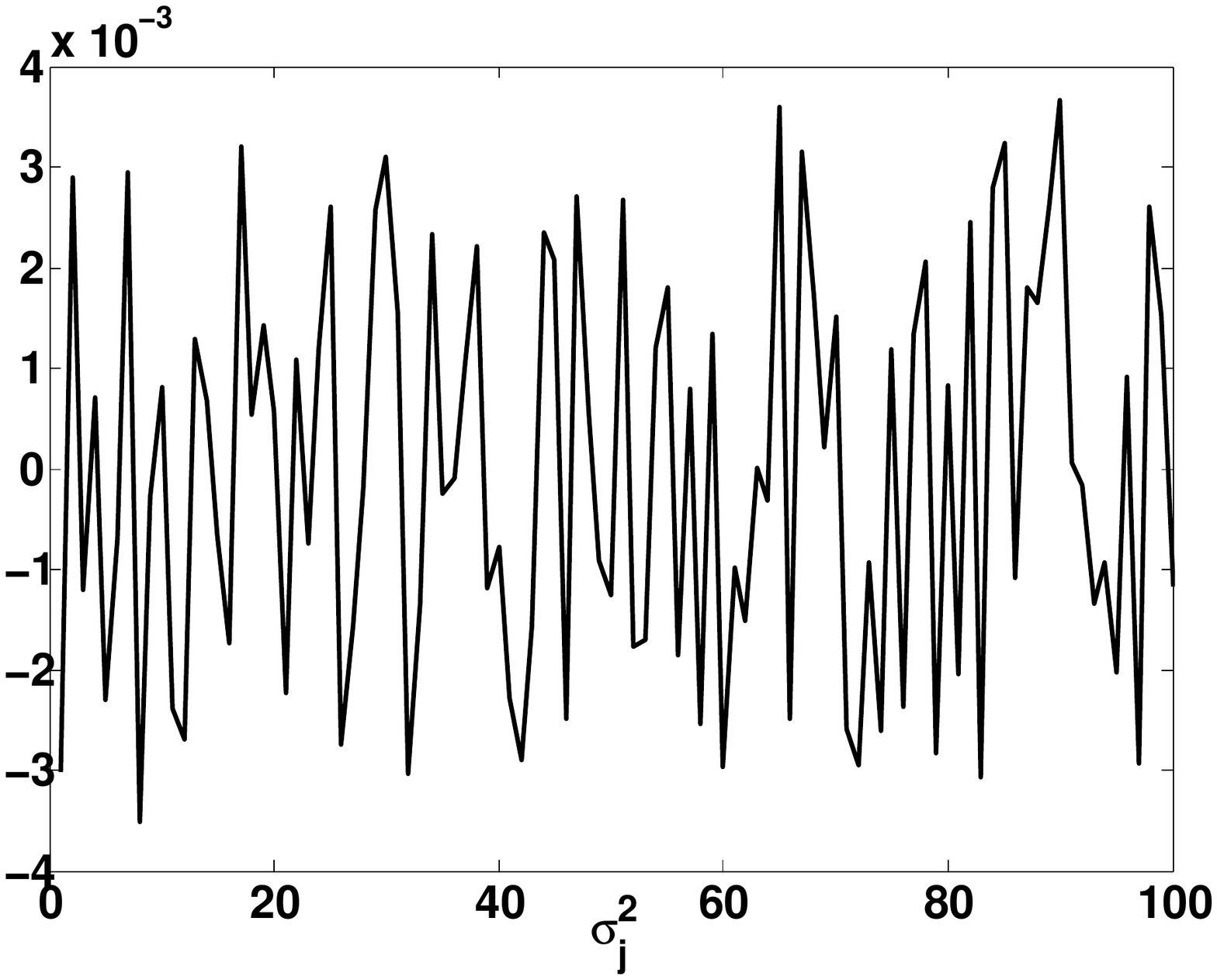}}
\caption[Independence of random variables $\{\mathbf{c}, \sigma^2_1, ..., \sigma^2_{N_C}, \sigma^2_{\mathbf{c}}\}$]{Correlation coefficients $\rho_{\mathbf{c},\sigma^2_j} (j=1, ..., N_C)$, $\rho_{\mathbf{c},\sigma^2_{\mathbf{c}}}$ and $\rho_{\sigma^2_i,\sigma^2_{\mathbf{c}}} (j=1, ..., N_C)$ between the vectors $\mathbf{c}$, $\sigma^2_{\mathbf{c}}$ and $\sigma^2_j$. The noise level $\delta=0.001$. (a) Correlation coefficients $\rho_{\mathbf{c},\sigma^2_j}$. (b) Correlation coefficients $\rho_{\mathbf{c},\sigma^2_{\mathbf{c}}}$. (c) Correlation coefficients $\rho_{\sigma^2_i,\sigma^2_{\mathbf{c}}}$. }
\label{correlation}
\end{figure}

Now, we investigate the robustness of Algorithm 2. We first study the sensitivity of the numerical results with respect to the prior parameters $\{\alpha_j, \beta_j\}^{N_C}_{j=1}$ and $\{\alpha_{\mathbf{c}}, \beta_{\mathbf{c}}\}$, and illustrate the rationale behind these choices in the following examples. If the value of $\alpha_j, \alpha_{\mathbf{c}}$ is not too large, the distributions of $\sigma^2_j$ and $\sigma^2_{\mathbf{c}}$ in (\ref{result}) are dominated by $N_T/2$ and $n/2$ in the numerators $\alpha_j+N_T/2$ and $\alpha_{\mathbf{c}}+n/2$, and can thus only have a marginal impact on the algorithm. Therefore, we only study the influence of parameters $\beta_j, \beta_{\mathbf{c}}$. Since $\{\beta_j\}^{N_C}_{j=1}$ plays the same role in the algorithm, we only present the result for $\beta_1$ here. The results are displayed in Table \ref{priorParameters}, where the data noise level equals 1\%, and parameters $\beta_j$ ($j>1$), $\{\alpha_j\}$ and $\alpha_{\mathbf{c}}$ are all set as 1. The results show that both the estimated rate constant mean map and the conditional independence assumption (\ref{independence}) seem relatively independent of the parameters $\beta_{\mathbf{c}}$ and $\beta_1$ for the magnitude of the relative error L2Err, and the correlation coefficients $|\rho_{\sigma^2_i,\sigma^2_{\mathbf{c}}}|_{\infty}$, $|\rho_{\mathbf{c},\sigma^2_j}|_{\infty}$ and $|\rho_{\mathbf{c},\sigma^2_{\mathbf{c}}}|_{\infty}$ remain almost unchanged as the value of $\beta_{\mathbf{c}}$ and $\beta_1$ varies.

\begin{table}[!htb]
\caption{Sensitivity analysis for parameters $\beta_{\mathbf{c}}$ and $\beta_1$. $\delta=1\%$. $\beta_j=1$ for $j\neq 1$. $\alpha_j=\alpha_{\mathbf{c}}=1$.}
  \centering
   \begin{center}
   \begin{tabular}{cccccc}
   \hline
    $\beta_{\mathbf{c}}$ & $\beta_1$ & L2Err & $|\rho_{\sigma^2_i,\sigma^2_{\mathbf{c}}}|_{\infty}$ & $|\rho_{\mathbf{c},\sigma^2_j}|_{\infty}$ & $|\rho_{\mathbf{c},\sigma^2_{\mathbf{c}}}|_{\infty}$\\
    \hline
    0.0001 & 0.0001 & 0.0661 & 0.0032 & 0.0130 & 0.2247\\
    0.01& 0.01 & 0.0659 & 0.0002 & 0.0304 & 0.0414\\
    1 & 1 & 0.0658 & 0.0029 & 0.0232 & 0.2175\\
    100 & 100 & 0.0661 & 0.0021 & 0.0258 & 0.1502\\
    10000 & 10000 & 0.0660 & 0.0038 & 0.0268 & 0.0931\\
    \hline
    0.0001 & 10000 & 0.0664 & 0.0092 & 0.0268 & 0.1406\\
    0.01 & 100 & 0.0667 & 0.0020 & 0.0836 & 0.2609\\
    0.1 & 10 & 0.0660 & 0.0021 & 0.0135 & 0.3592\\
    \hline
   \label{priorParameters}
  \end{tabular}
  \end{center}
\end{table}

We next consider the influence of the initial guess of the distribution for the rate constant map $q^0(\mathbf{c})$. In Table \ref{priorMap} we investigate the case with the Gaussian prior. The expectation is assumed as the Tikhonov regularized solution $\mathbf{c}^*=\left( \mathbf{K}^T \mathbf{K} + \lambda \mathbf{I} \right)^{-1} \mathbf{K}^T \mathbf{R}$, where the regularization parameter $\lambda$ varies from 0.00001 to 100. The covariance matrix is assumed to be of the type $\kappa \mathbf{I}$. The first four lows in Table \ref{priorMap} show that with a small value of $\lambda$ the relative error is small for $\kappa\leq1$. The last four lows in Table \ref{priorMap} show that the smaller the value of $\lambda$, the better the estimated rate constant mean map. Other types of mean $\mathbf{c}_0$ and covariance matrix $\Sigma_{\mathbf{c}}$ were also tested for this example. Worse results are obtained for other types of mean $\mathbf{c}_0$. However, arbitrary types of covariance matrix $\Sigma_{\mathbf{c}}$ with small values of $\|\Sigma_{\mathbf{c}}\|\leq1$ give the same good results if one chooses the mean $\mathbf{c}_0$ appropriately. Therefore, based on this group of simulations, we suggest using the initial guess of $q^0(\mathbf{c})$ as a Gaussian distribution with mean $\left( \mathbf{K}^T \mathbf{K} + \lambda \mathbf{I} \right)^{-1} \mathbf{K}^T \mathbf{R}$ $(\lambda\leq0.001)$ and covariance matrix $\Sigma_{\mathbf{c}}$ ($\|\Sigma_{\mathbf{c}}\|\leq1$).

\begin{table}[!t]
\caption{Sensitivity analysis for the initial guess of $q^0(\mathbf{c})$ with fixed type of distribution $\mathcal{N}\left( \left( \mathbf{K}^T \mathbf{K} + \lambda \mathbf{I} \right)^{-1} \mathbf{K}^T \mathbf{R} , \kappa \mathbf{I} \right)$. $\delta=1\%$. $\alpha_j=\beta_j=\alpha_{\mathbf{c}}=\beta_{\mathbf{c}}=1$.}
  \centering
   \begin{center}
   \begin{tabular}{cccccc}
   \hline
    $\lambda$ & $\kappa$ & L2Err & $|\rho_{\sigma^2_i,\sigma^2_{\mathbf{c}}}|_{\infty}$ & $|\rho_{\mathbf{c},\sigma^2_j}|_{\infty}$ & $|\rho_{\mathbf{c},\sigma^2_{\mathbf{c}}}|_{\infty}$\\
    \hline
    0.00001 & 0.0001 & 0.0078 & 0.0008 & 0.0281 & 0.0162\\
    0.00001 & 0.01 & 0.0080 & 0.0026 & 0.0143 & 0.1272\\
    0.00001 & 1 & 0.0125 & 0.0016 & 0.0137 & 0.0381\\
    0.00001 & 10 & 0.0531 & 0.0022 & 0.0101 & 0.2360\\
    \hline
    0.001 & 0.01 & 0.0517 & 0.0009 & 0.0061 & 0.2440\\
    0.1 & 0.01 & 0.0621 & 0.0019 & 0.0145 & 0.0828\\
    1 & 0.01 & 0.2179 & 0.0098 & 0.0178 & 0.1738\\
    100 & 0.01 & 0.4512 & 0.0203 & 0.0782 & 0.3491\\
    \hline
   \label{priorMap}
  \end{tabular}
  \end{center}
\end{table}

\begin{table}[!b]
\caption{Sensitivity analysis of the algorithm with different noisy data. $\mathbf{c}_0\sim \mathcal{N}\left( \left( \mathbf{K}^T \mathbf{K} + \lambda \mathbf{I} \right)^{-1} \mathbf{K}^T \mathbf{R} , \mathbf{I} \right)$. $\alpha_j=\beta_j=\alpha_{\mathbf{c}}=\beta_{\mathbf{c}}=1$.}
  \centering
   \begin{center}
   \begin{tabular}{ccccc}
   \hline
    $\delta$ & L2Err & $|\rho_{\sigma^2_i,\sigma^2_{\mathbf{c}}}|_{\infty}$ & $|\rho_{\mathbf{c},\sigma^2_j}|_{\infty}$ & $|\rho_{\mathbf{c},\sigma^2_{\mathbf{c}}}|_{\infty}$\\
    \hline
    0.001 & 0.0098 & 0.0011 & 0.0263 & 0.0459 \\
    0.005 & 0.0111 & 0.0024 & 0.0193 & 0.3285 \\
    0.01 & 0.0318 & 0.0023 & 0.0540 & 0.1941 \\
    0.02 & 0.1520 & 0.0015 & 0.0164 & 0.5070 \\
    0.10 & 0.4365 & 0.0037 & 0.0313 & 0.1497 \\
    0.20 & 0.4253 & 0.0029 & 0.0462 & 0.2979 \\
    0.40 & 0.5049 & 0.0018 & 0.0440 & 0.3387 \\
    \hline
   \label{TabNoise}
  \end{tabular}
  \end{center}
\end{table}

\begin{figure}[!htb]
\centering
\subfigure[]{
\includegraphics[clip, trim=0 3.8in 0 1in, width=2.2in]{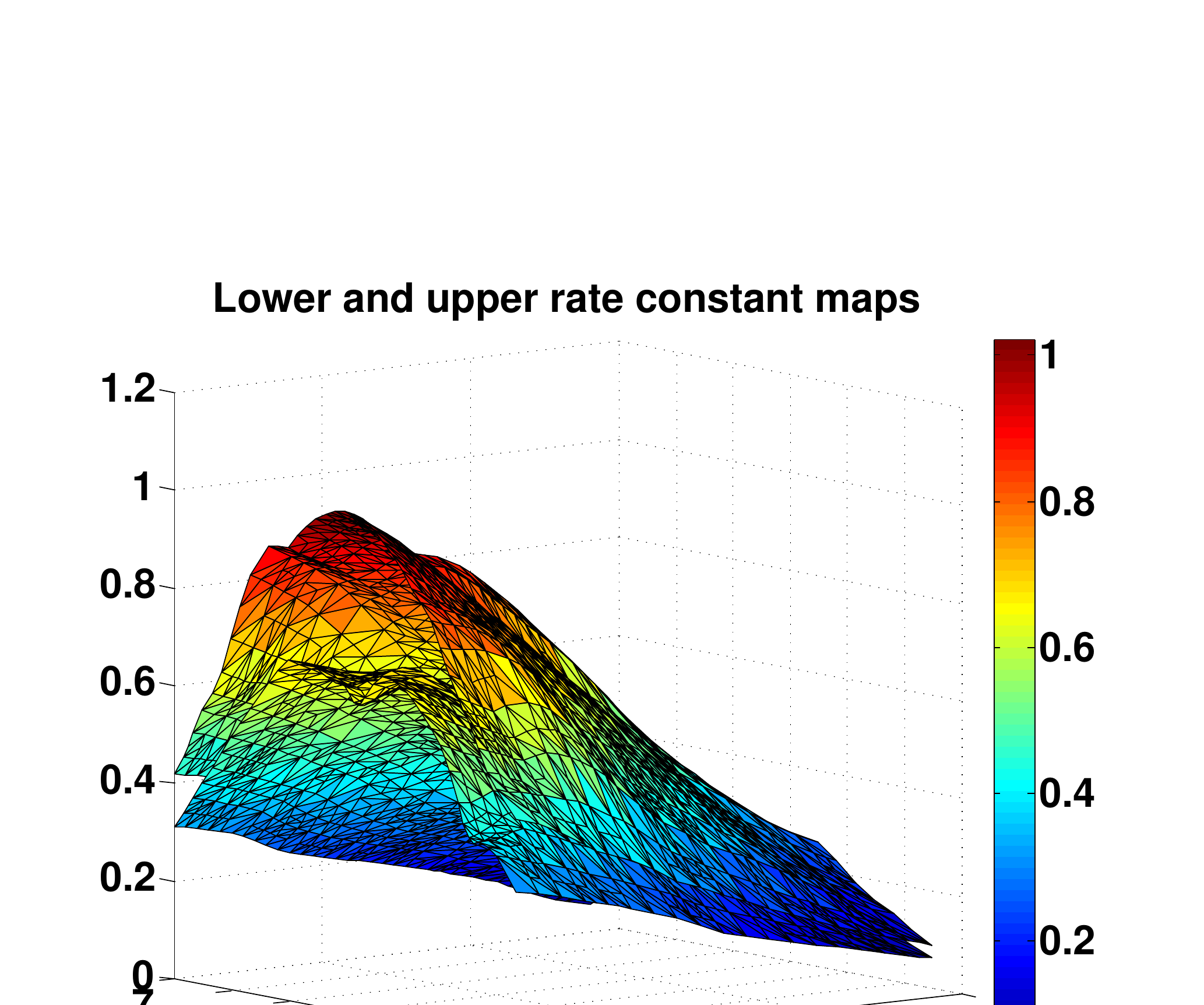}}
\subfigure[]{
\includegraphics[clip, trim=0 2.8in 0 1.8in, width=2.2in]{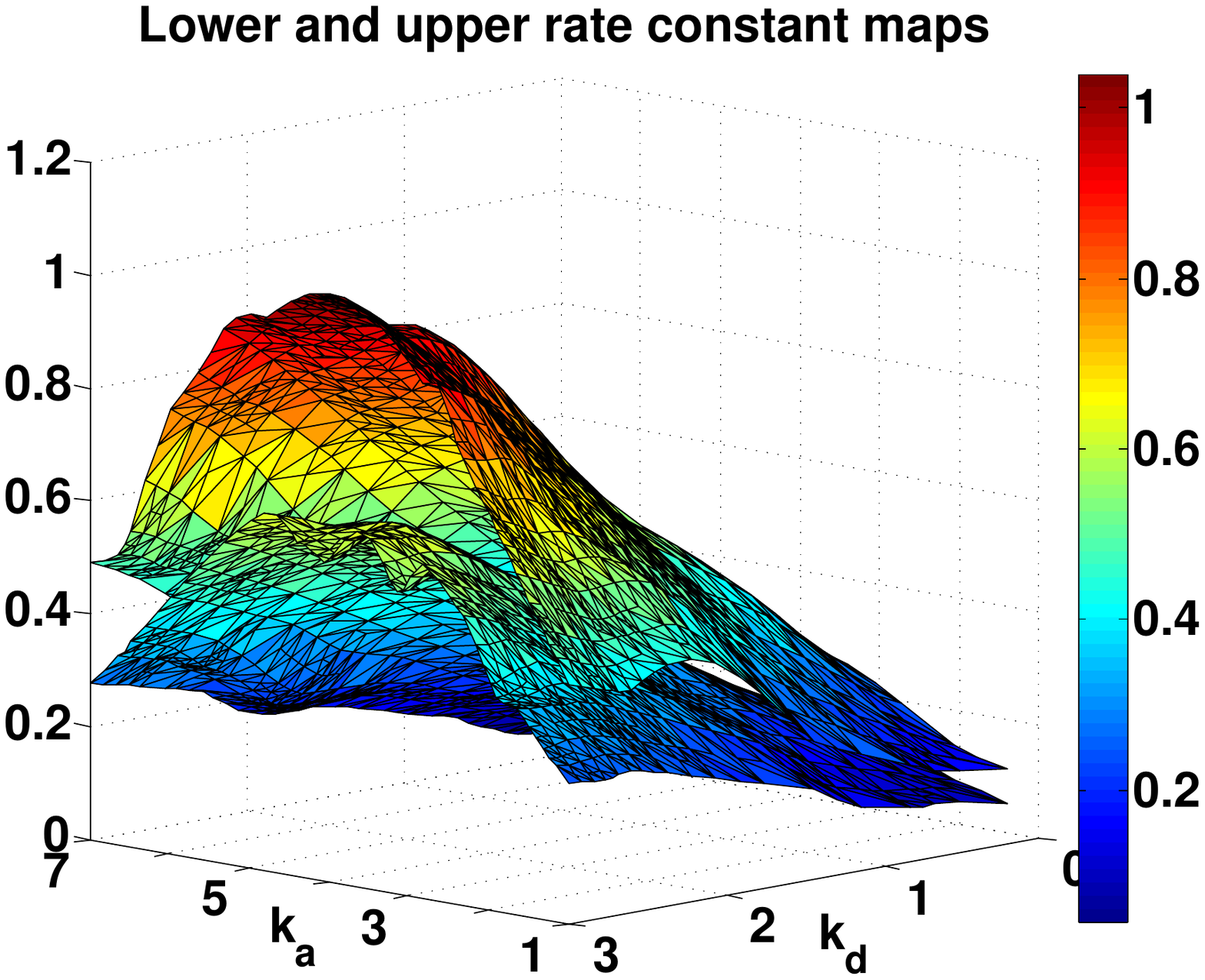}}
\subfigure[]{
\includegraphics[clip, trim=0 2.8in 0 1.8in, width=2.2in]{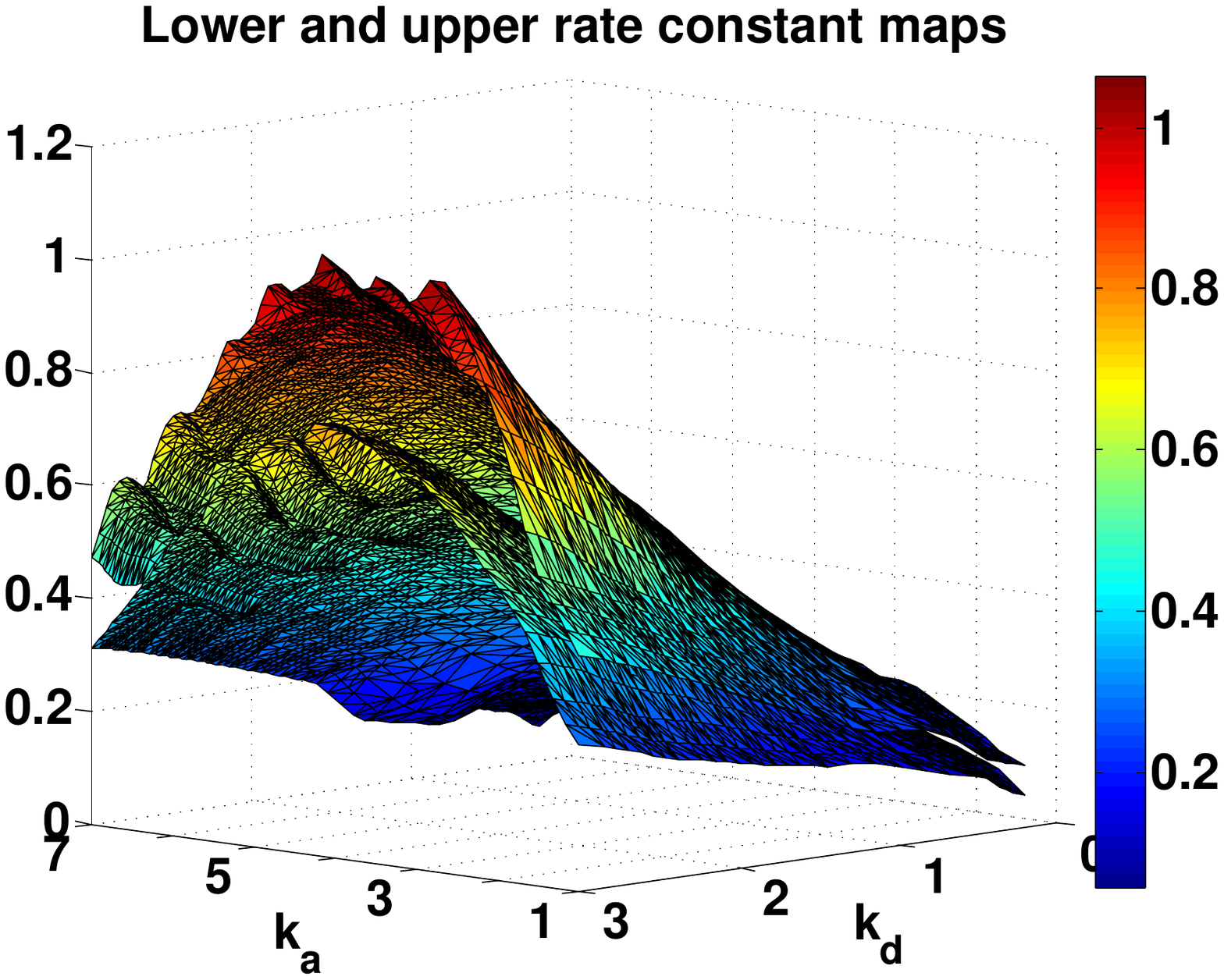}}
\subfigure[]{
\includegraphics[clip, trim=0 2.8in 0 1.8in, width=2.2in]{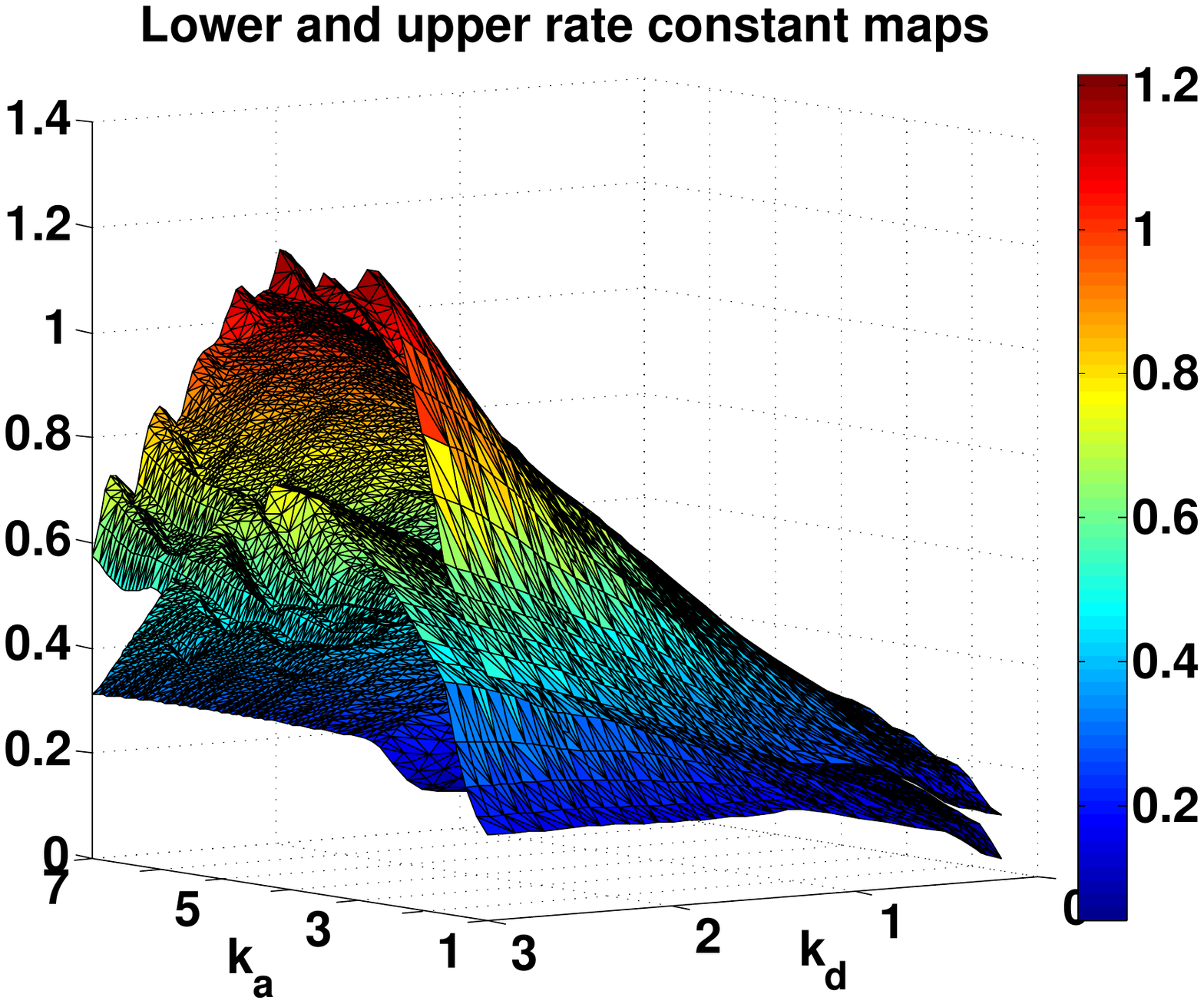}}
\caption[The estimated lower and upper rate constant maps for different noise data]{The estimated lower and upper rate constant maps for different noise data. The initial triangulation: the number of nodes equals 100, and the number of triangles is 162. (a) $\delta=0.001$. The final triangulation: the number of nodes equals 1261, and the number of triangles is 2396. (b) $\delta=0.005$. The final triangulation: the number of nodes equals 1146, and the number of triangles is 2166. (c) $\delta=0.01$. The final triangulation: the number of nodes equals 3930, and the number of triangles is 7652. (d) $\delta=0.02$. The final triangulation: the number of nodes equals 4566, and the number of triangles is 8914.}
\label{FigNoise}
\end{figure}

Now, we investigate the stability of the estimated solutions with respect to the noise level $\delta$.
Specifically, set $\delta=0.001, 0.005, 0.01, 0.02$ and implement Algorithm 2 repeatedly.  The initial guess of expectation is set as the Tikhonov regularized solution $\mathbf{c}^*=\left( \mathbf{K}^T \mathbf{K} + \lambda \mathbf{I} \right)^{-1} \mathbf{K}^T \mathbf{R}$, where the regularization parameter $\lambda$ is chosen by the $L$-curve method~\ctp{engl1996regularization}. The pairs of parameters $\{ (\alpha_j,\beta_j) \}^{N_C}_{j=1}$ and $(\alpha_{\mathbf{c}},\beta_{\mathbf{c}})$ for the inverse Gamma distribution are taken to be $(1,1)$. The relative errors of the estimate rate constant mean map and the correlation between parameters are displayed in Table \ref{TabNoise}. The estimated lower and upper rate constant maps for different noise data are shown in Figure \ref{FigNoise}. These results demonstrate that our algorithm is stable with respect to a small noise error ($\delta\leq0.02$). It should be pointed out that with a large noise data, based on $\delta\geq0.1$, the estimated rate constant mean map is quite poor and the value of $\|\Sigma^{\mathbf{c}}_*\|$ is relatively large, though the algorithm is still convergent with the given system parameters.

\begin{figure}[!t]
\centering
\subfigure[]{
\includegraphics[clip, trim=0 2.8in 0 2.8in, width=2.0in]{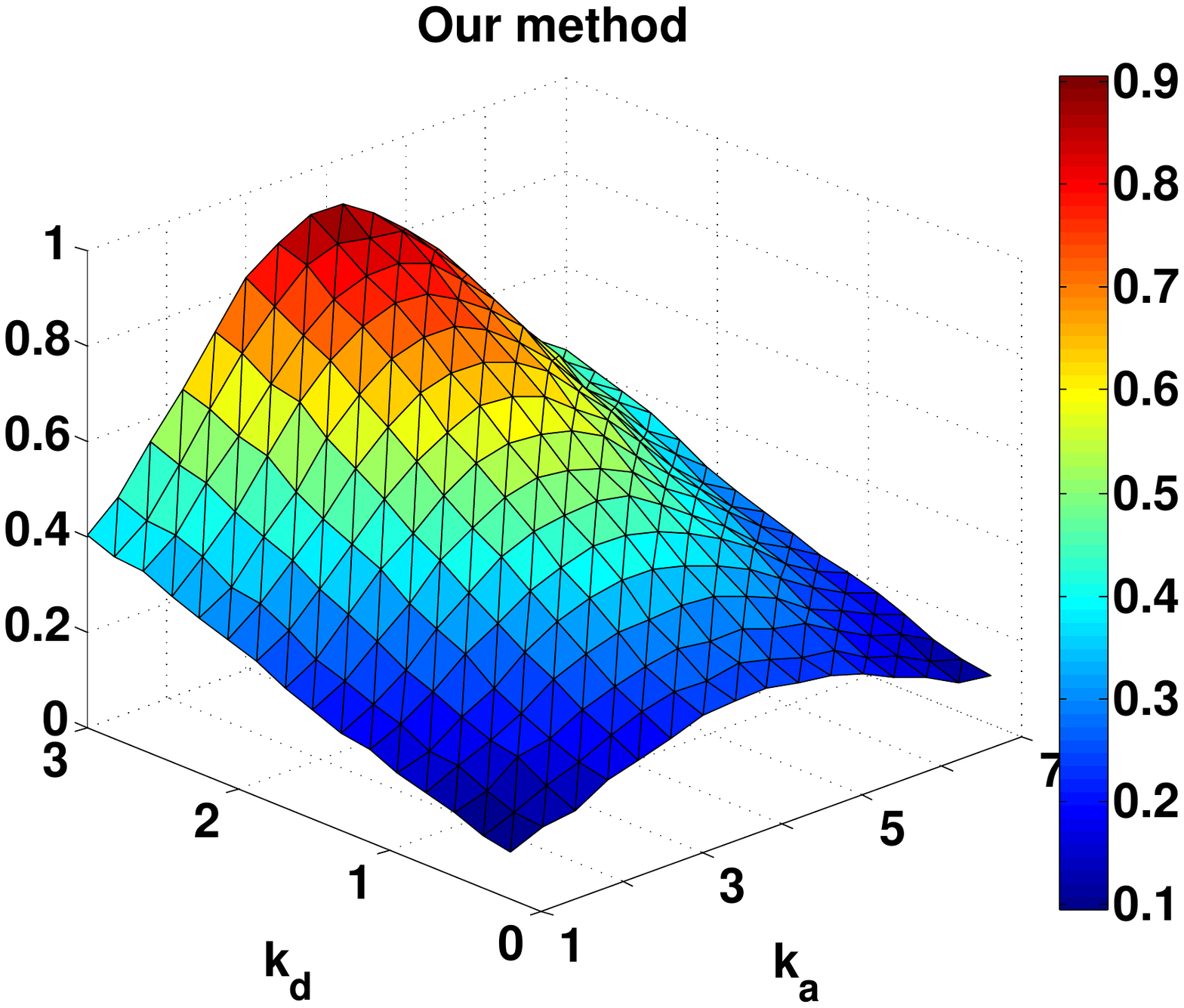}}
\subfigure[]{
\includegraphics[clip, trim=0 2.8in 0 2.8in, width=2.0in]{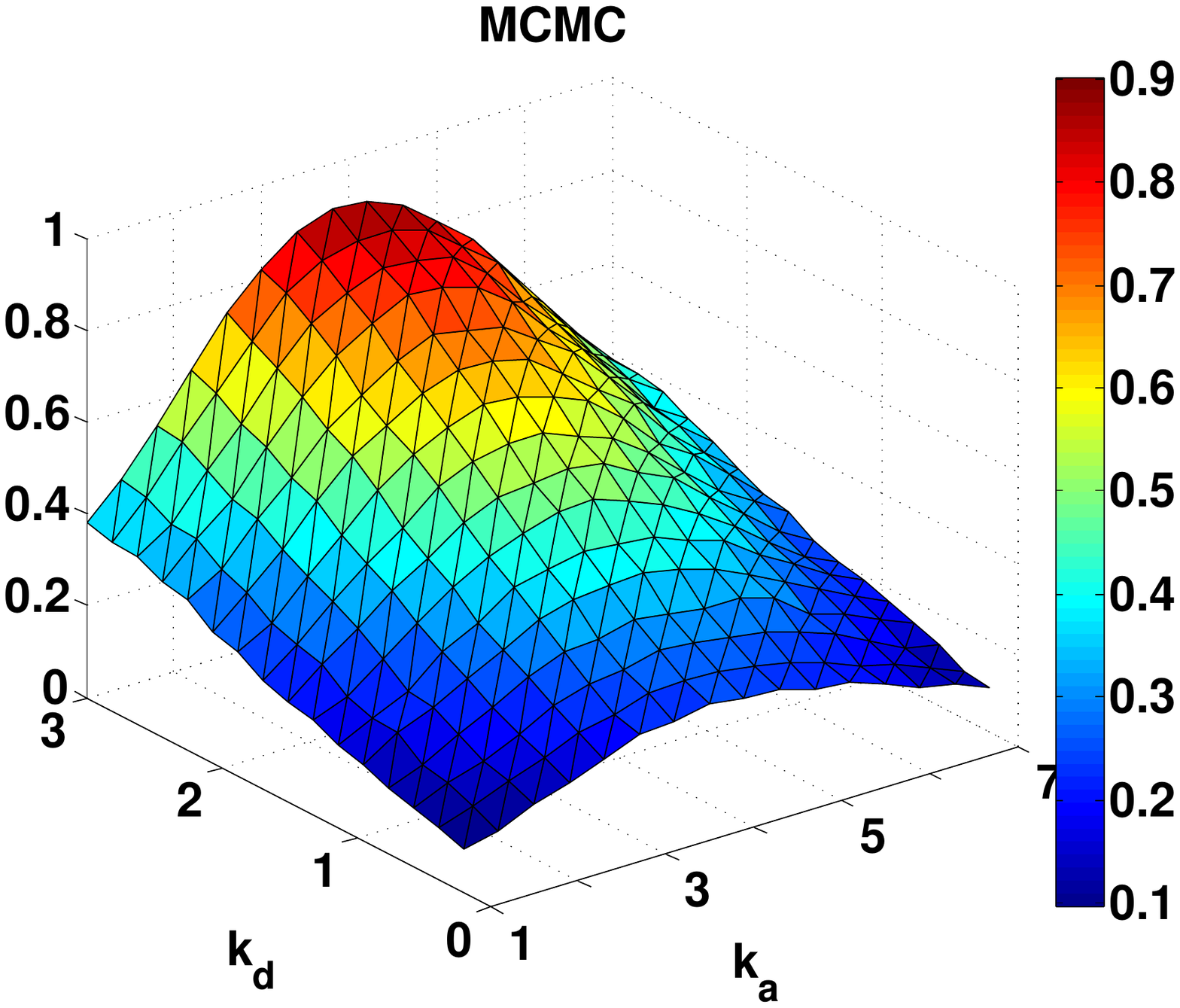}}
\subfigure[]{
\includegraphics[clip, trim=0 2.8in 0 2.8in, width=2.0in]{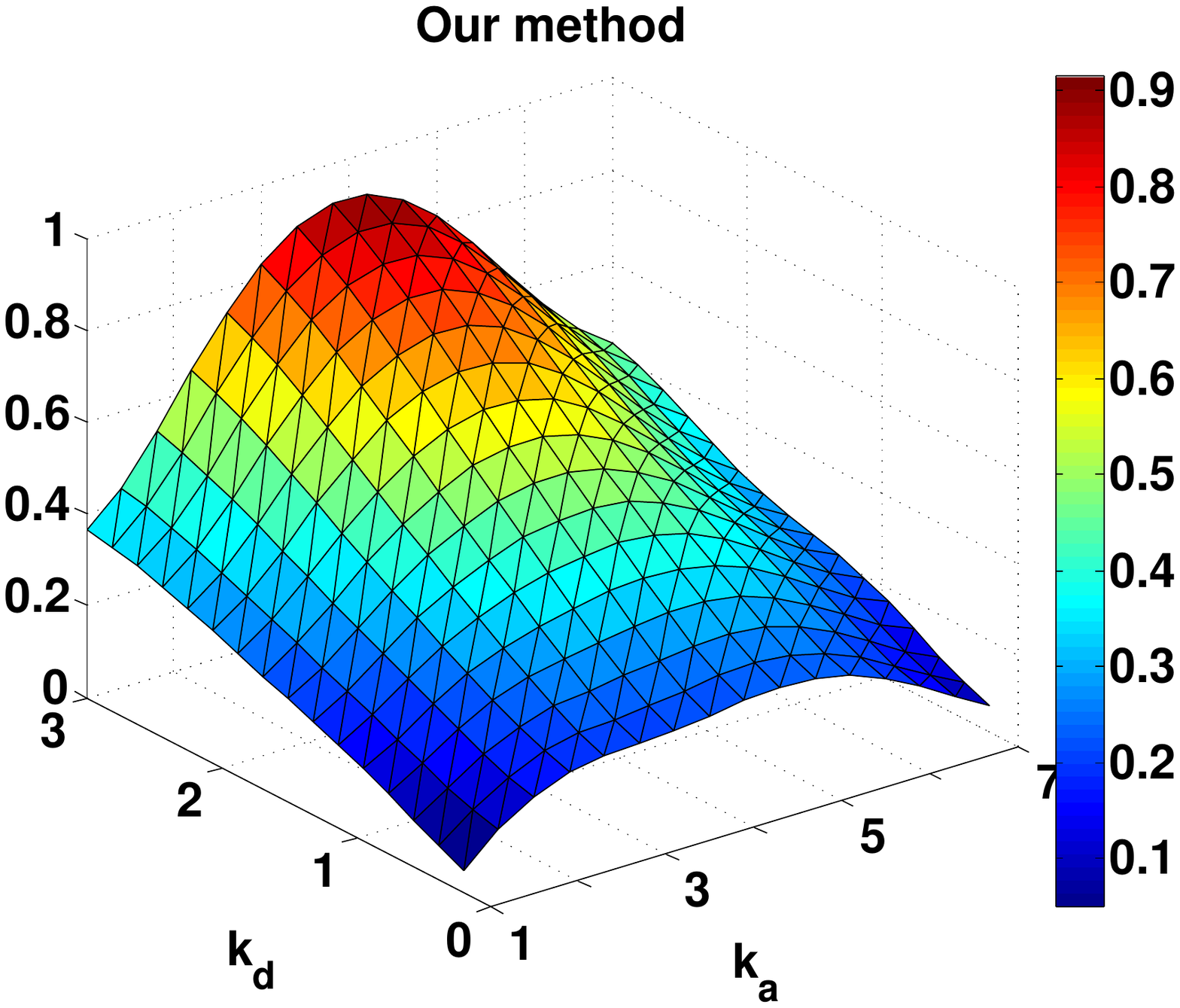}}
\subfigure[]{
\includegraphics[clip, trim=0 2.8in 0 2.8in, width=2.0in]{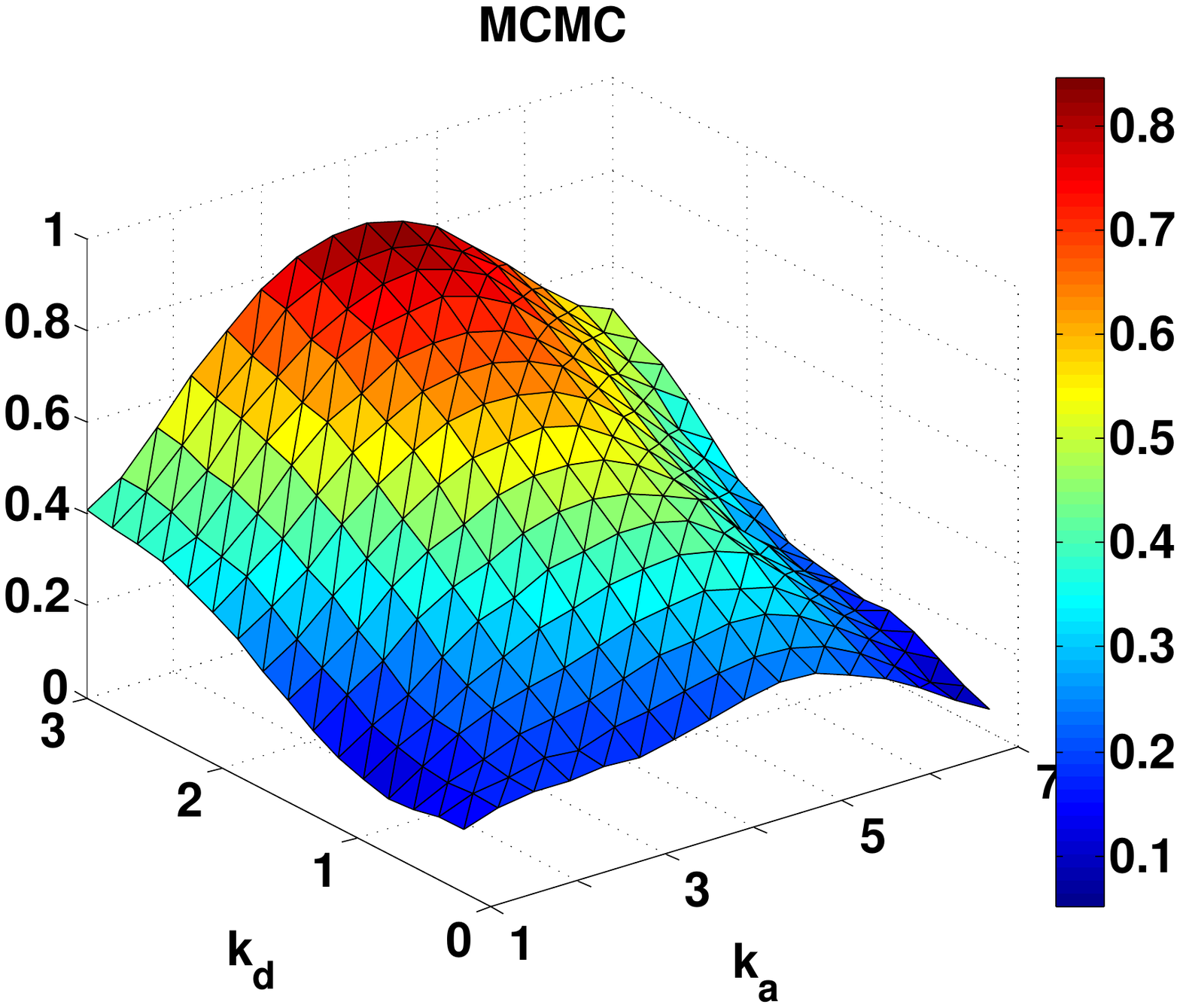}}
\subfigure[]{
\includegraphics[clip, trim=0 2.8in 0 2.8in, width=2.0in]{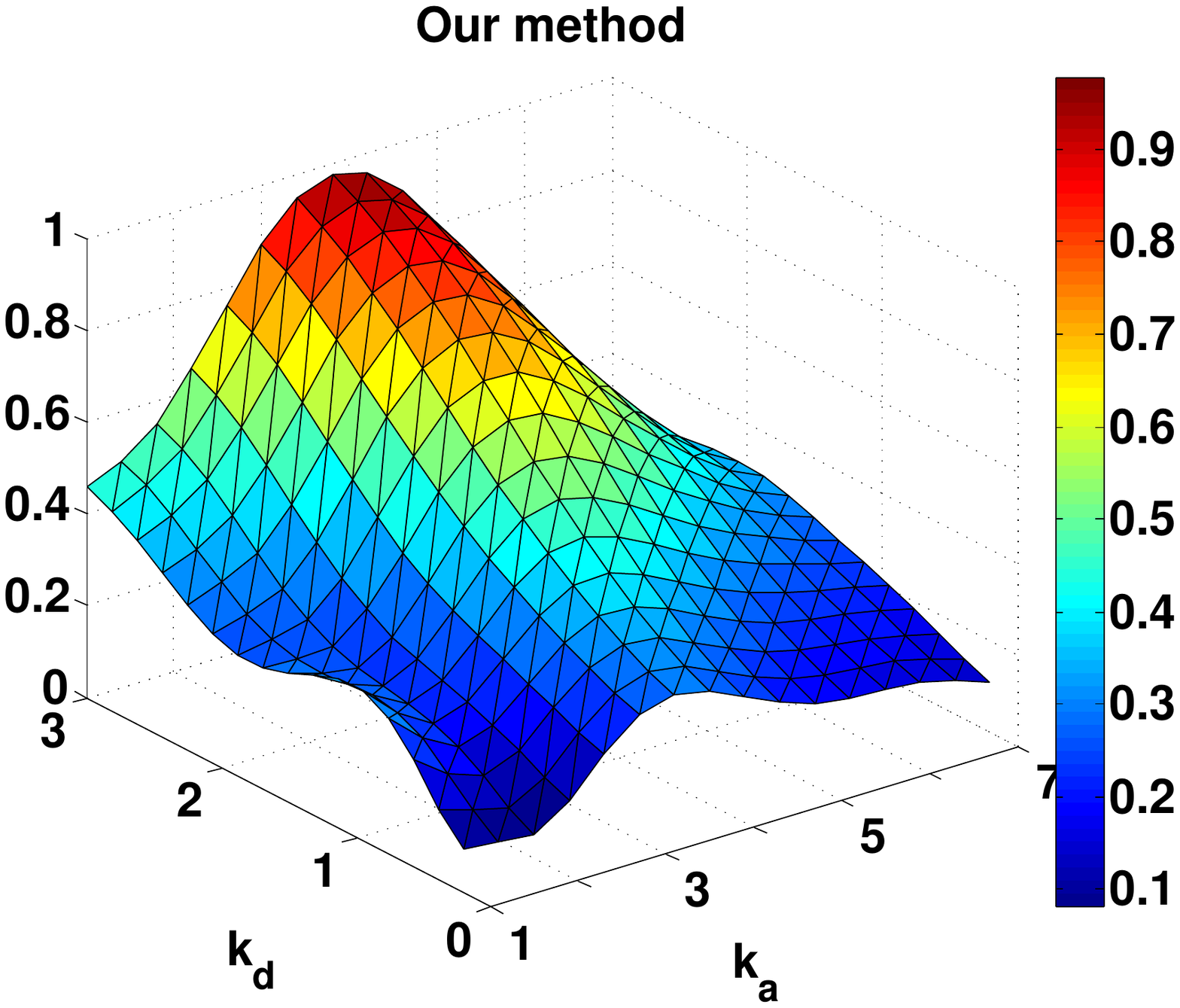}}
\subfigure[]{
\includegraphics[clip, trim=0 2.8in 0 2.8in, width=2.0in]{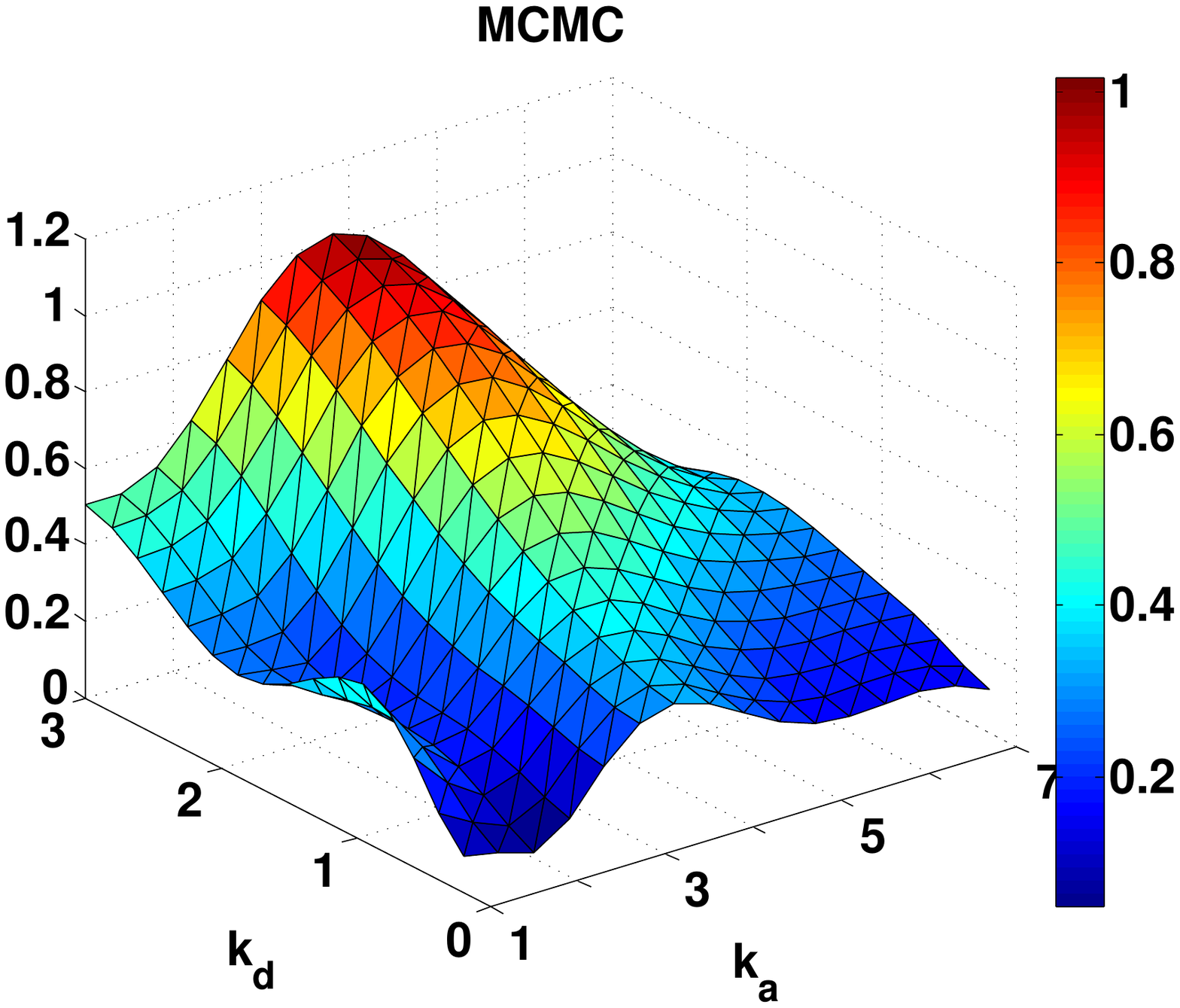}}
\caption[Comparison with MCMC]{A comparison of our method (Algorithm 2) with MCMC. (a) Algorithm 2 for $\delta=0.002$. (b) MCMC for $\delta=0.002$. (c) Algorithm 2 for $\delta=0.01$. (d) MCMC for $\delta=0.002$. (e) Algorithm 2 for $\delta=0.01$. (f) MCMC for $\delta=0.01$. (g) Algorithm 2 for $\delta=0.01$. $\delta=0.05$. (h) MCMC for $\delta=0.05$.}
\label{ComparisonFig}
\end{figure}

Finally, we compare our approach with the MCMC. For simplicity, we fix the grid with 450 triangles and 256 nodes, and employ Algorithm 2 and the MCMC with the same initial guess of parameters ($\alpha_j=\beta_j=\alpha_{\mathbf{c}}=\beta_{\mathbf{c}}=1$ and $\mathbf{c}_0\sim \mathcal{N}\left(\left( \mathbf{K}^T \mathbf{K} + \lambda \mathbf{I} \right)^{-1} \mathbf{K}^T \mathbf{R} , \mathbf{I} \right)$). The resulting posterior distribution (\ref{ppdfEnd}) is sampled using the standard Gibbs sampler, and the length of the Markov chain is 50,000, with the first 5,000 realizations discarded as transient states. The mixing of the Markov chain is monitored by visually inspecting the trace plot and calculating the correlation coefficient. Note that we set a maximal implement time of $T_{max}=300$ min. in all of our simulations. We obtained an output of the result even in the occurrence of the maximal time point of the algorithm. Numerical results with different noise levels are shown in Figure \ref{ComparisonFig} and Table \ref{ComparisonTab}. Furthermore, in order to quantitatively evaluate the difference between our method and MCMC, we compute their relative approximate Wasserstein distance $W$, which is defined as follows: let $\Omega^{VB}_j$ and $\Omega^{MCMC}_j$ be the samples of the two groups (sample sizes are different) of $\mathbf{c}_j$ from our method and MCMC, respectively. Then, the relative approximate Wasserstein distance for $\mathbf{c}_j$ between our method and MCMC is defined as
\begin{equation*}\label{Wasserstein}
W_{j} := \frac{1}{100} \sum^{100}_{i=1} \left| \tilde{c}^{VB,i}_{j} - \tilde{c}^{MCMC,i}_{j} \right| \Big/ \left( \tilde{c}^{MCMC,i}_{j} + 1 \right),
\end{equation*}
where $\tilde{c}^{VB,i}_{j}$ and $\tilde{c}^{MCMC,i}_{j}$ represent the $i$\% quantiles from two sample sets $\Omega^{VB}_j$ and $\Omega^{MCMC}_j$, respectively.  The calculated relative approximate Wasserstein distances for all $\mathbf{c}_j$ ($j=1, ..., 256$) with different data sets are presented in Figure \ref{Wasserstein}(a). The histograms of $\mathbf{c}_{100}$ by the data set with $\delta=0.01$ for both our method and MCMC are provided in Figure \ref{Wasserstein}(b). All of these outputs indicate that the results obtained using our method and the MCMC almost coincided for the small noise level, while for large noisy data the estimated mean map achieved through our method was much better than the results acquired using the MCMC. It may even be possible (e.g., in a case with noise level $\delta=0.4$ in Table \ref{ComparisonTab}) that the MCMC does not converge, while our algorithm converges quickly. Moreover, as one can see in Table \ref{ComparisonTab}, it clearly shows the computational efficiency of our method over the MCMC for all situations. Finally, it should be noted that when the problem size is large, e.g., $n>1000$, for the majority of cases the MCMC cannot offer a satisfactory result within a reasonable time (e.g. two days). Instead, our method provided a convergent result within a short time. For instance, in our real data application, when $n=1812$, the implementation of our method took only 42 minutes.

\begin{table}[!t]
\caption{Relative error and running time of the estimated rate constant mean map for our method and MCMC. $\mathbf{c}_0\sim \mathcal{N}\left( \left( \mathbf{K}^T \mathbf{K} + \lambda \mathbf{I} \right)^{-1} \mathbf{K}^T \mathbf{R} , \mathbf{I} \right)$. $\alpha_j=\beta_j=\alpha_{\mathbf{c}}=\beta_{\mathbf{c}}=1$. }
\begin{center}
\begin{tabular}{cccccc} \hline
\multirow{2}{*}{$\delta$} &
\multicolumn{2}{c}{Relative error} &
\multicolumn{2}{c}{Running time (min.)} \\
\cline{2-5}
& MCMC & AVBA & MCMC & AVBA  \\ \hline
    0.002 & 0.0024 & 0.0022 & 21.18 & 1.75 \\
    0.01 & 0.0112 & 0.0098 & 29.45 & 3.01 \\
    0.05 & 0.0638 & 0.0232 & 25.77 & 3.68 \\
    0.10 & 0.2840 & 0.1632 & 73.12 & 7.23 \\
    0.20 & 0.5114 & 0.2915 & 49.94 & 5.68 \\
    0.40 & 3.8116 & 0.7451 & $T_{max}$ & 7.68 \\ \hline
\end{tabular}
\end{center}
\label{ComparisonTab}
\end{table}

\begin{figure}[!t]
\centering
\subfigure[]{
\includegraphics[clip, trim=0 3.0in 0 2.2in, width=2.4in]{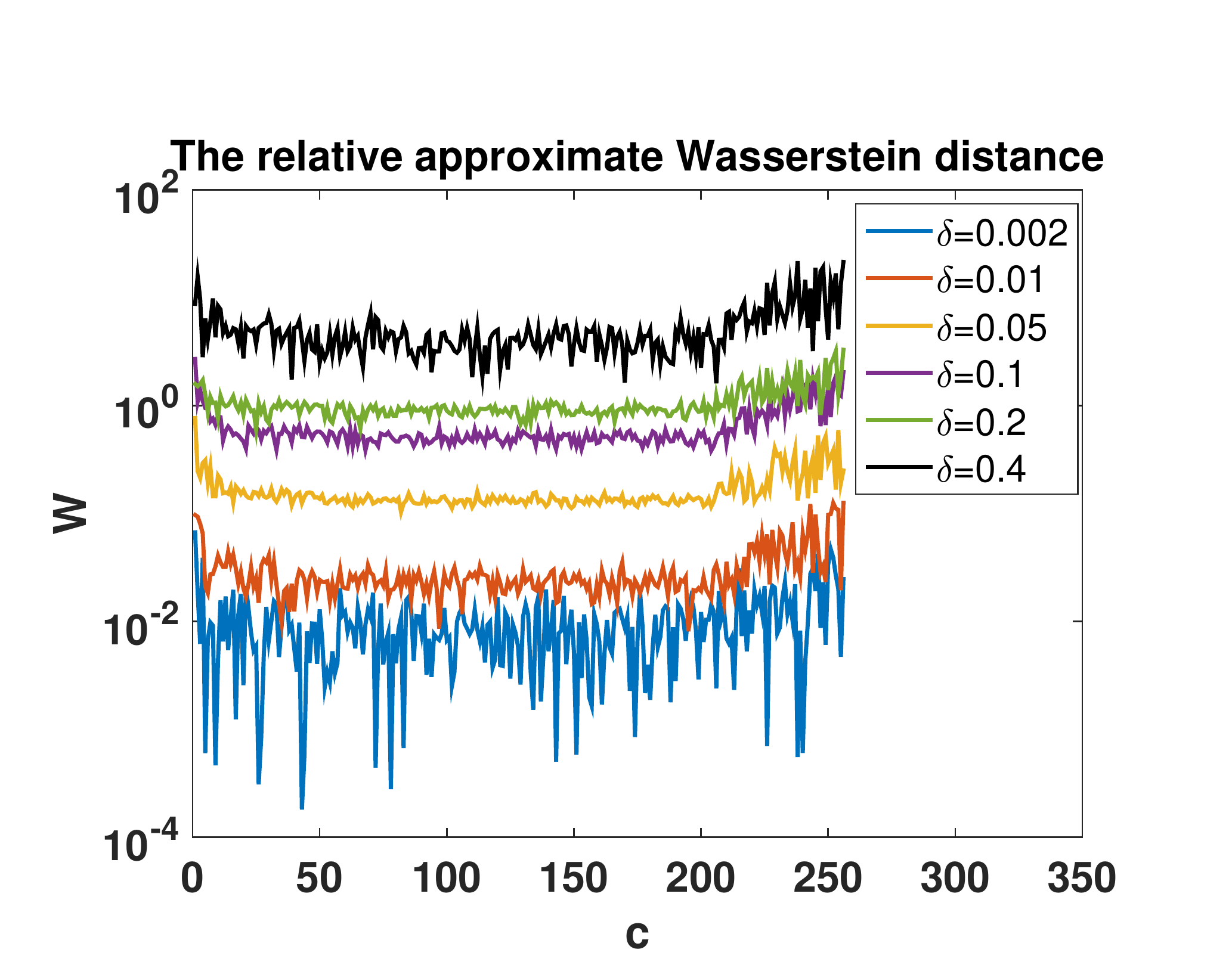}}
\subfigure[]{
\includegraphics[clip, trim=0 3.0in 0 2.2in, width=2.4in]{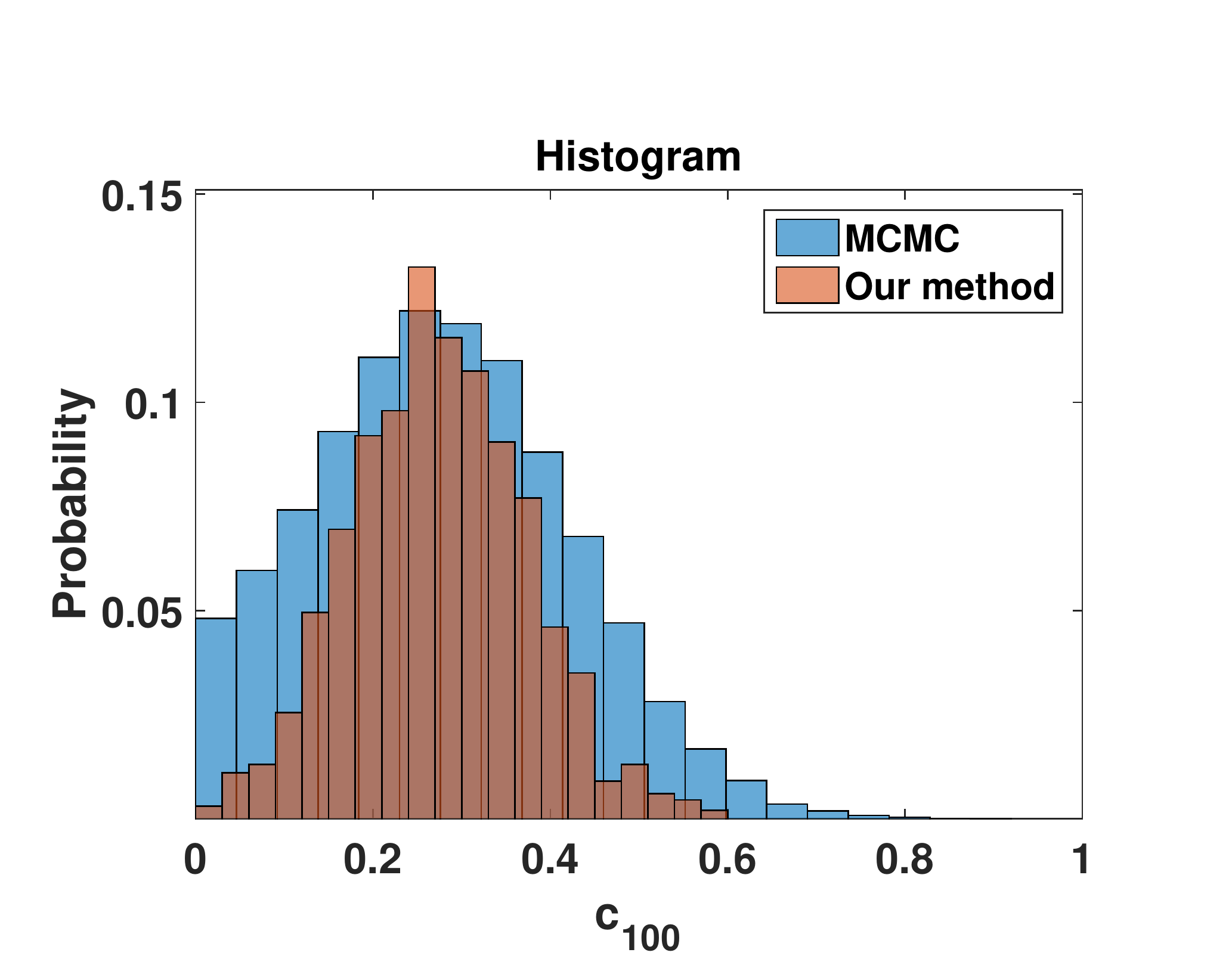}}
\caption[Mean Map]{(a) The relative approximate Wasserstein distance between our method the MCMC of $\mathbf{c}$ for data sets with different noise level $\delta$. (b)  Histograms for $\mathbf{c}_{100}$ by the data set with $\delta=0.01$ for our method (2735 samples) and MCMC (50000 samples). The difference between two histograms in (b) corresponds to a point with the abscissa 100 on the red curve in (a).}
\label{Wasserstein}
\end{figure}

\section{A real data application}\label{realData}

\begin{figure}[!htb]
\centering
\subfigure[]{
\includegraphics[clip, trim=0.5in 2.8in 0.5in 2.8in, width=2.4in]{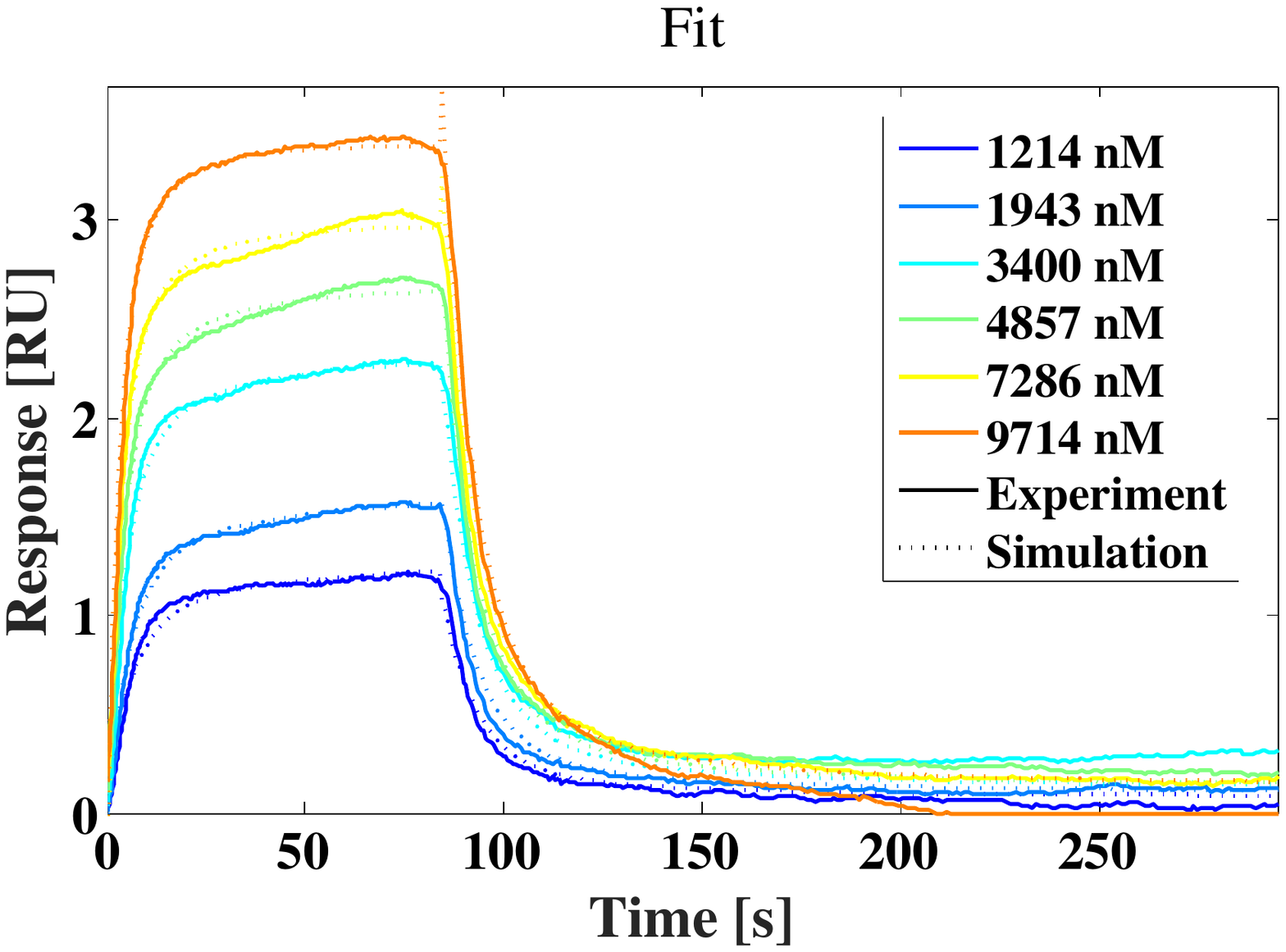}}
\subfigure[]{
\includegraphics[clip, trim=0.5in 2.8in 0.5in 2.8in, width=2.4in]{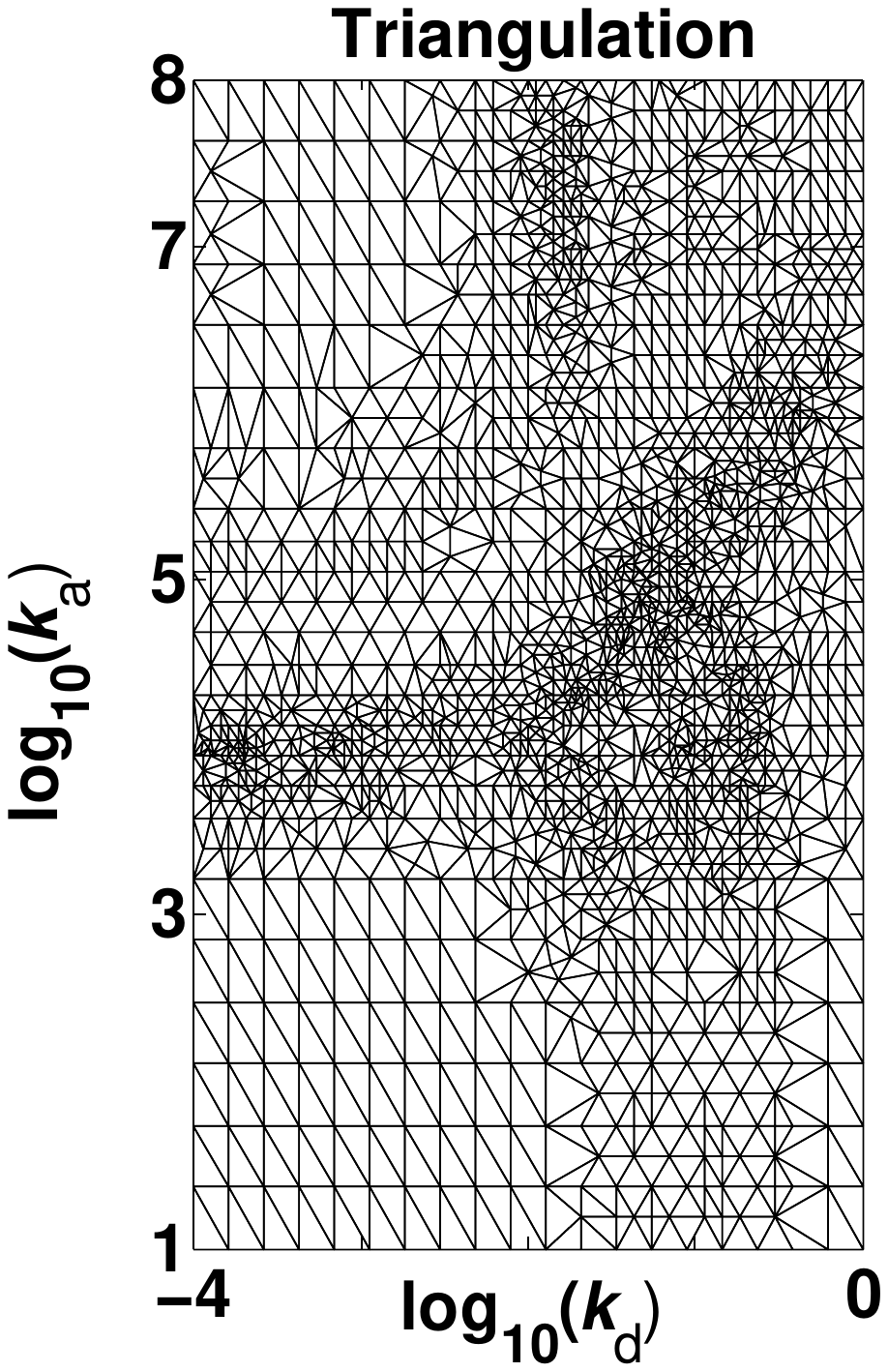}}
\caption[Data fit]{(a) The experimental data (solid line) and simulated response (dashed line) with different concentrations of PTH at the temperature of 20 $^\circ C$. (b) Triangulations at the final iteration. The number of nodes equals 1812, and the number of triangles is 3498. }
\label{DataSimulation}
\end{figure}

In this section, our method (AVBA) is tested on real experimental data~-- parathyroid hormone (PTH). In the experiment, the human PTH1R receptor was immobilized on a LNB-carboxyl biosensor chip using amine coupling according to the manufacturer's instructions. Using the flow rate 25 $\mu L/ min$ at 20.0$^{\circ}$C, we did 35 $\mu L$ injections of the peptide PTH(1-34) at six concentration levels from 1214 nM to 9714 nM; see the solid lines in (a) of Figure \ref{DataSimulation}. The sensorgrams were measured using a QCM biosensor Attana Cell 200  (Attana AB, Stockholm, Sweden) instrument.

\begin{figure}[!htb]
\centering
\subfigure[]{
\includegraphics[clip, trim=0 2.8in 0 2.8in, width=2.4in]{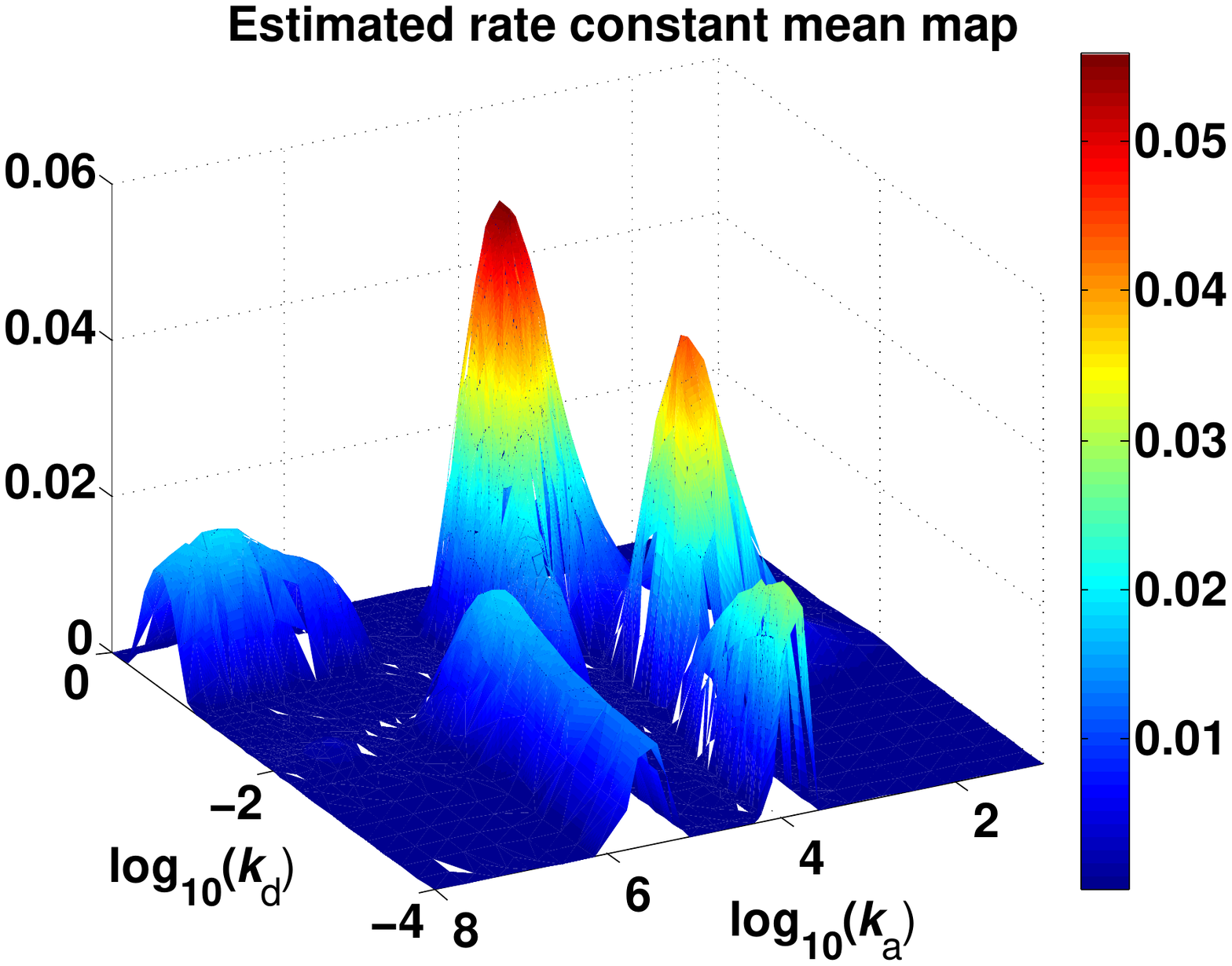}}
\subfigure[]{
\includegraphics[clip, trim=0 2.8in 0 2.8in, width=2.4in]{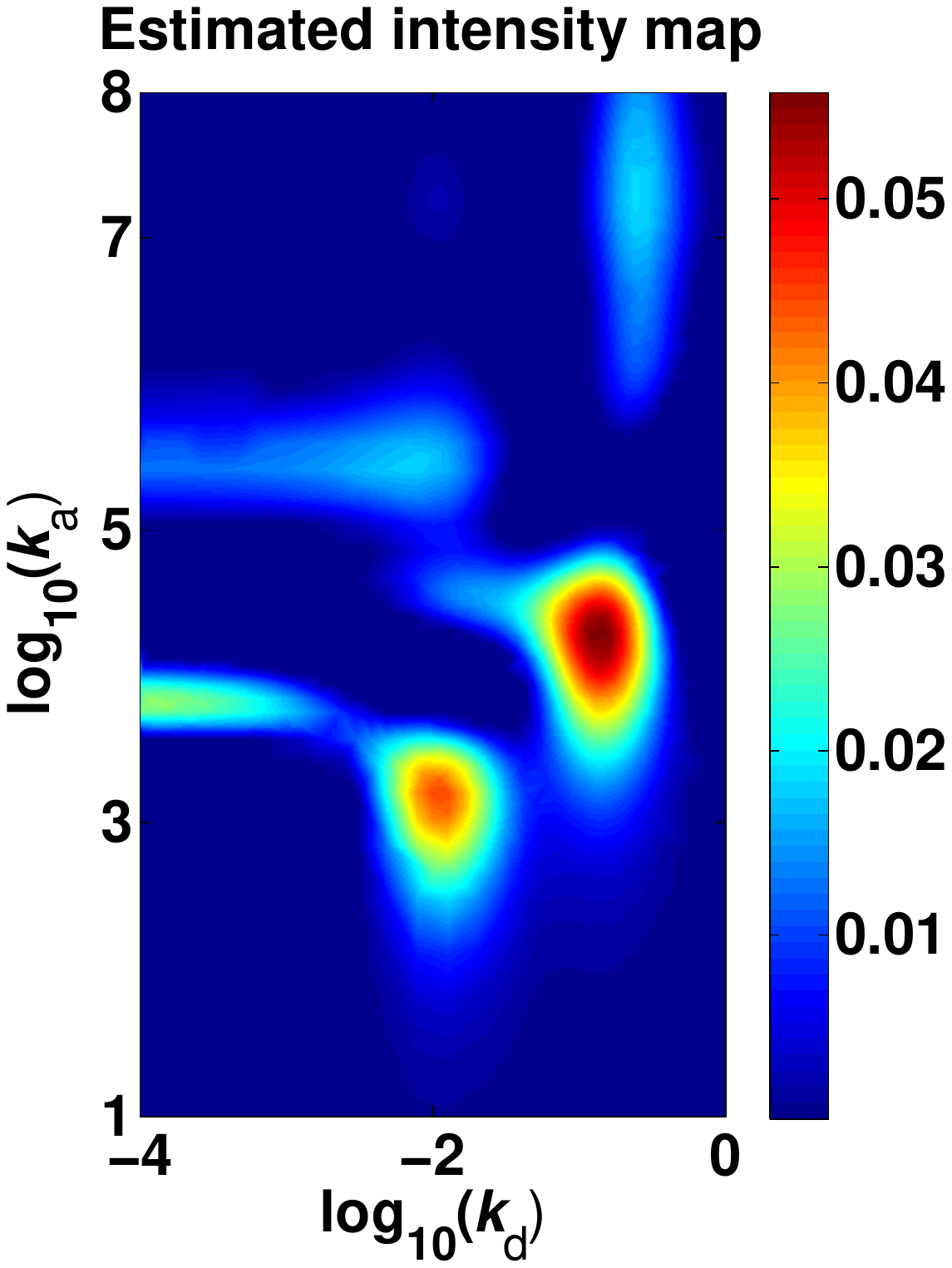}}
\caption[Mean Map]{(a) The estimated rate constant mean map. (b) The intensity map. }
\label{Result1}
\end{figure}

The initial triangulation is uniformly distributed in the log-scale domain $(\log_{10}(k_d), \log_{10}(k_a))\in [-4,0]\times[1,8]$ with $20\times20=400$ node points. Algorithm 2 stopped at the 9-th iteration with 1812 nodes and 3498 triangles~-- see (b) in Figure \ref{DataSimulation}. The estimated rate constant mean map and corresponding intensity map are shown in Figure \ref{Result1}. In (a) of Figure \ref{DataSimulation}, we show a comparison between the experimental data (solid line) and simulated response curves (dashed line), which is obtained by solving the forward problem with the estimated rate constant mean solution. The results with approximately 96\% overlap show that the estimated rate constant mean map can be used as the real rate constant map in a deterministic model.

As discussed in the Introduction, the peaks form the most important piece of information in the rate constant map, as they enable the conclusion of the number of interactions as well as the value of the active association and dissociation constants for the given large molecular system. These peaks can be considered as the geometry representation of the Intrinsic Property of the Interaction (IPoI). Note that the rate constant mean map in Figure \ref{Result1}, as well as the similar maps in \ctp{Svitel-2003}, \ctp{Altschuh-2012}, and \ctp{Gorshkova-2008}, contains too much redundant information about a chemical reaction besides the IPoI. Though we can derive the reaction information from the mean maps (see Figure \ref{Result1}) via ad-hoc work~-- as demonstrated in \ctp{Svitel-2003}, \ctp{Altschuh-2012}, \ctp{Gorshkova-2008} and \ctp{AIDA}~-- to improve the accuracy and efficiency of the estimation, we take an alternative perspective in this work; that is, we provide two new approaches that automatically derive the IPoI of the biosensor system.

\begin{figure}[!htb]
\centering
\subfigure[]{
\includegraphics[clip, trim=2in 2.8in 2in 2.8in, width=1.5in]{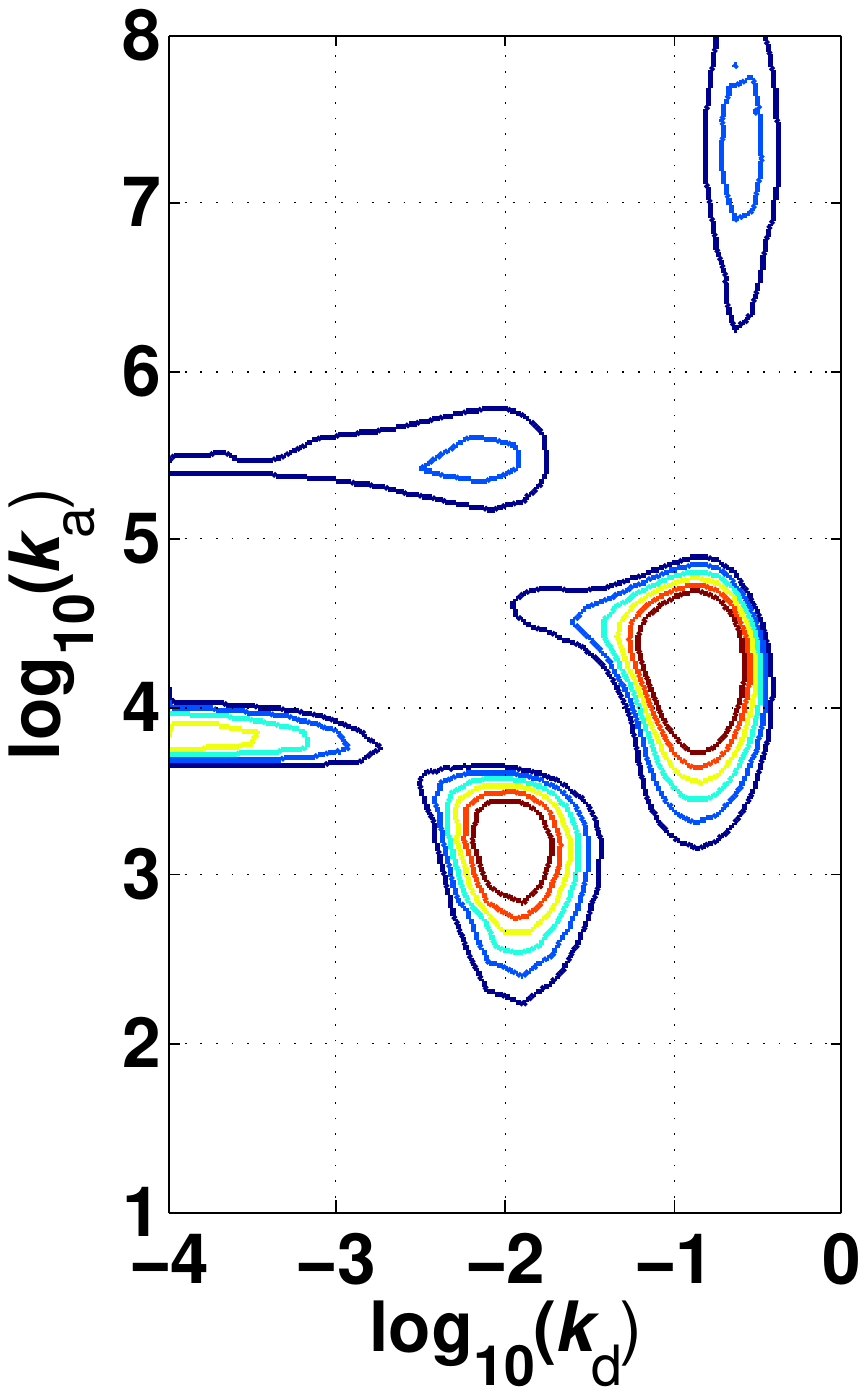}}
\subfigure[]{
\includegraphics[clip, trim=2in 2.8in 2in 2.8in, width=1.5in]{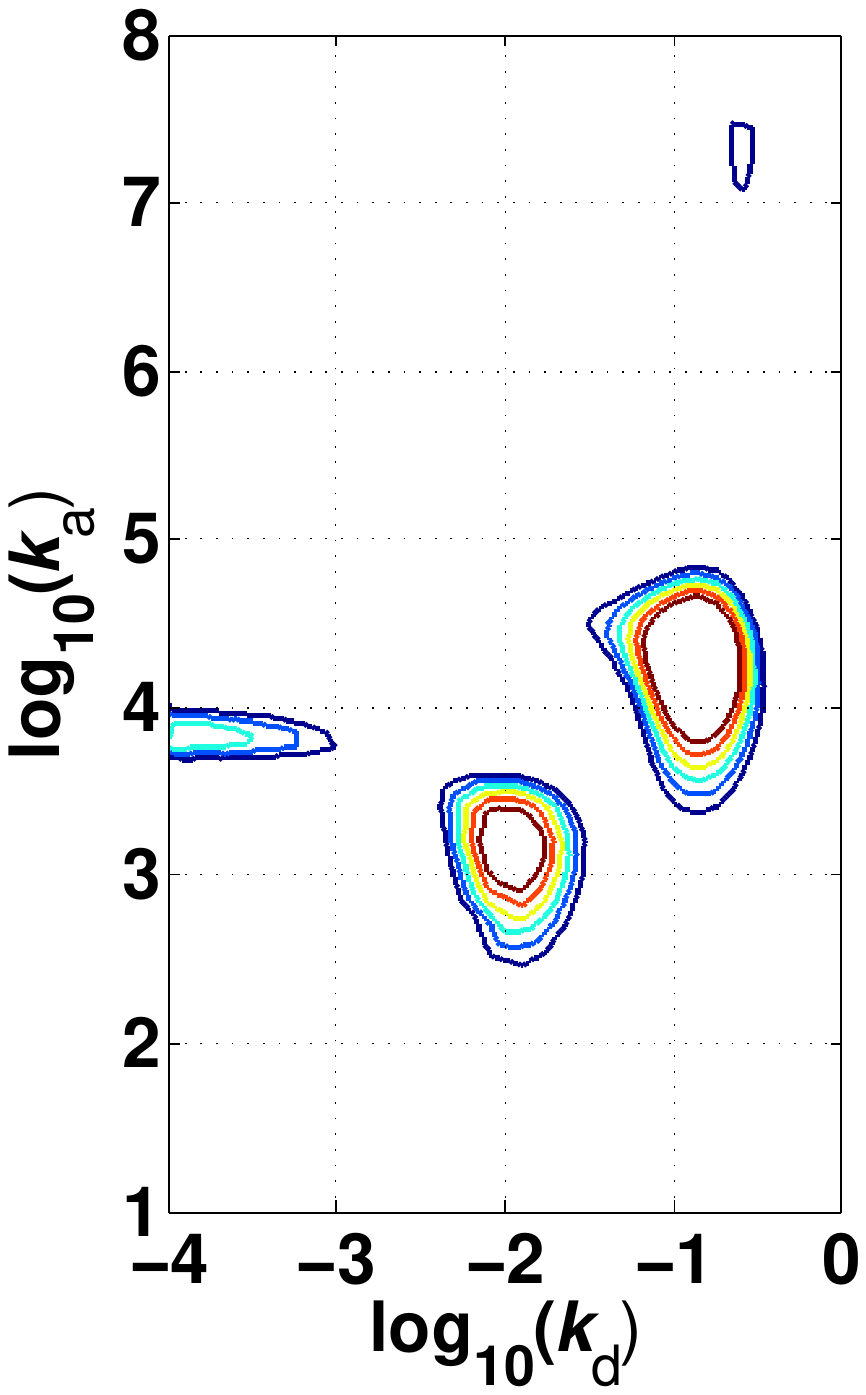}}
\subfigure[]{
\includegraphics[clip, trim=2in 2.8in 2in 2.8in, width=1.5in]{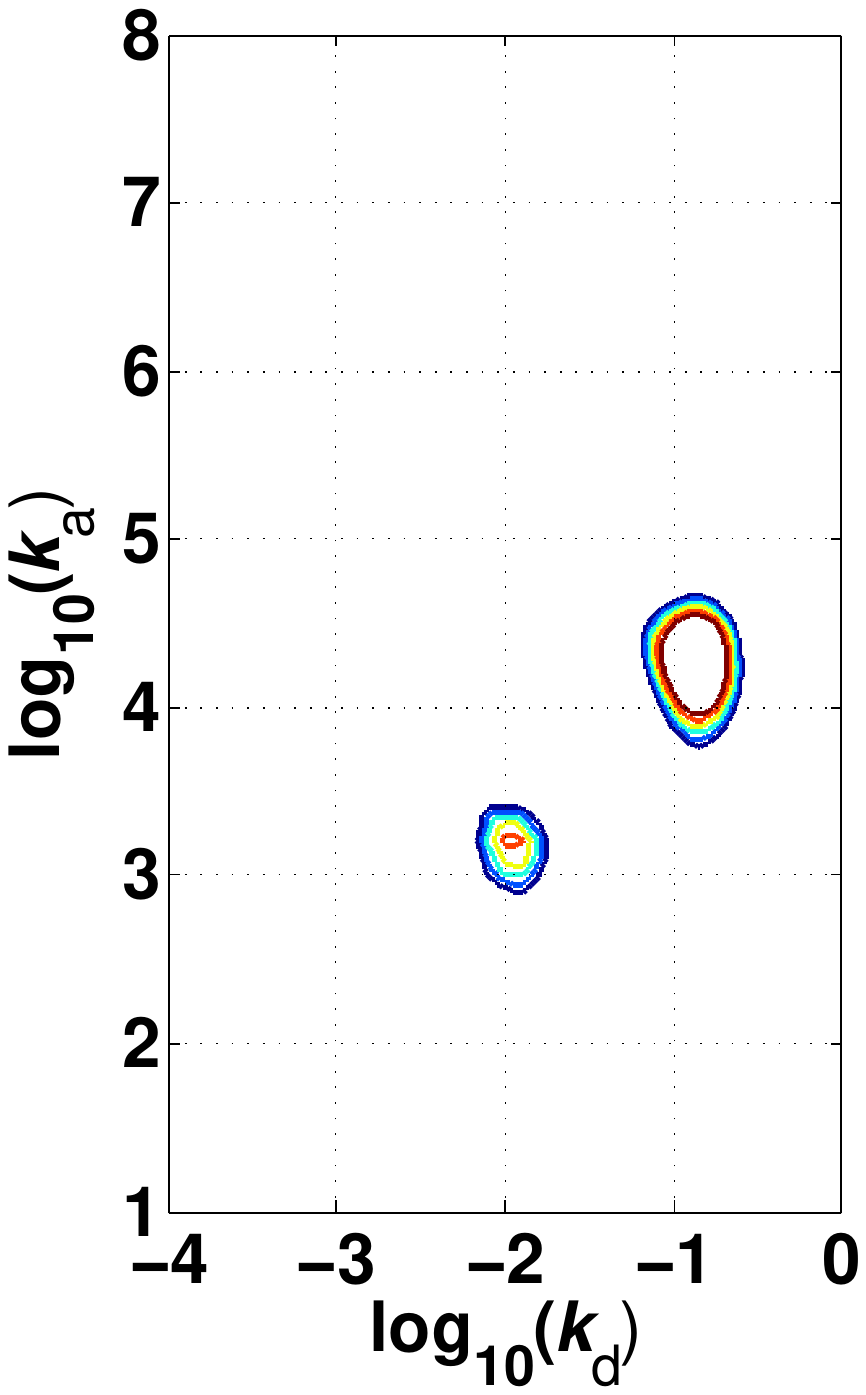}}
\caption[Contour]{TCM with different thresholds $\nu$: $\nu=15$ in (a), $\nu=10$ in (b), $\nu=5$ in (c). }
\label{Contour}
\end{figure}

The first approach is a thresholding method for the mean map, which will be named the $\nu$-Thresholding Contour Method (TCM). This approach cuts off all regions where the value of the rate constant mean map is less than the $1-\nu\%$ of the maximum value of the estimated rate constant. In Figure \ref{Contour}, we display the results of TCM with different thresholds $\nu$. Actually, with $\nu\in[2,8]$, the TCM provides exactly two isolated regions with some contours. We point out that the value of $\nu$ might vary case by case for the different datasets. However, $\nu\in[4,6]$ seems to work well for all tested data. Therefore, we recommend using $\nu=5$ in practice. From (c) of Figure \ref{Contour}, we conclude that two reactions exist for our PTH system, and their positions are $(\log_{10}(k_d), \log_{10}(k_a))=(-2.0,3.2)$ and $(-0.9, 4.4)$.

\begin{figure}[!t]
\centering
\subfigure[]{
\includegraphics[clip, trim=0 4.1in 0 0.6in, width=2.4in]{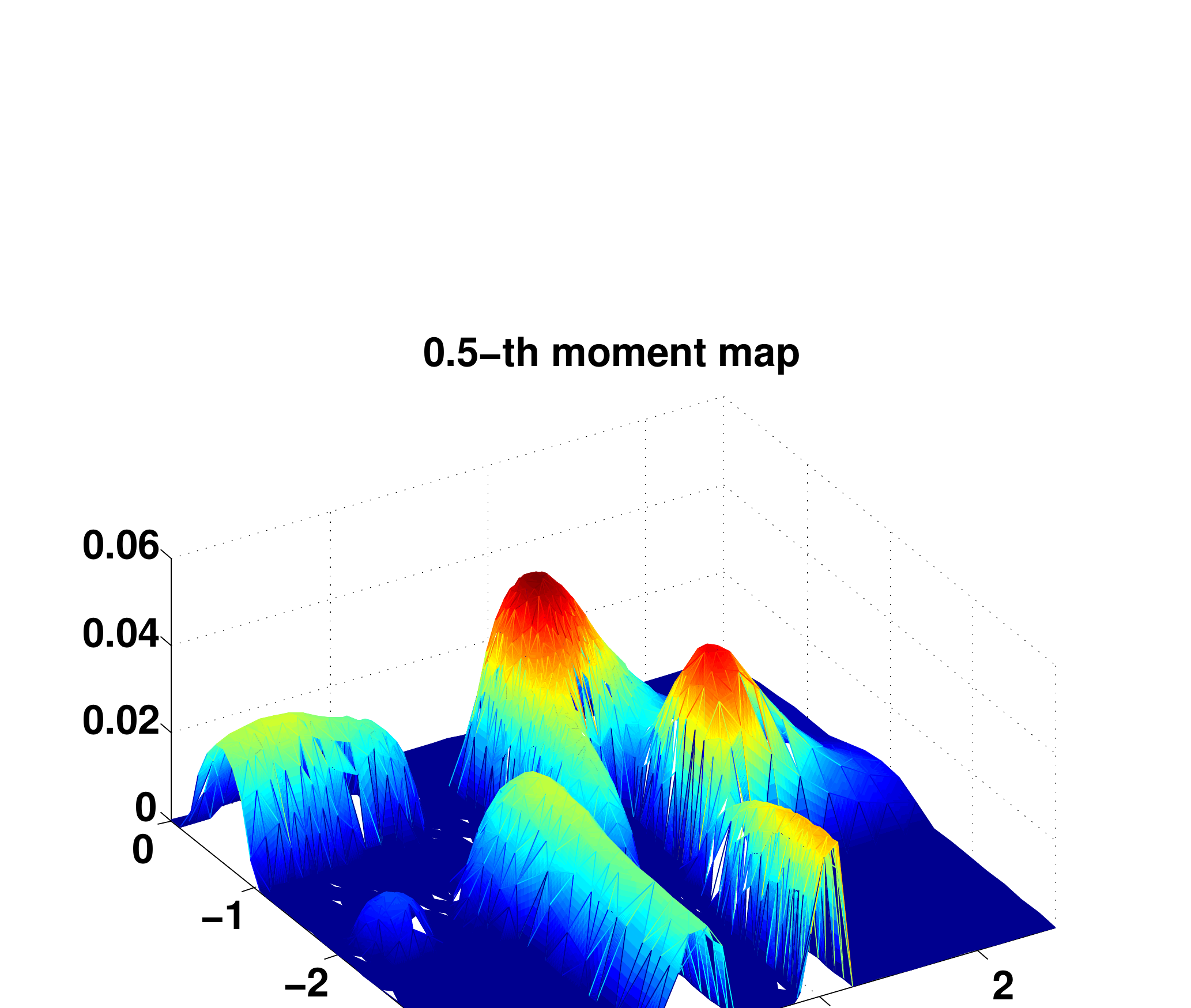}}
\subfigure[]{
\includegraphics[clip, trim=0 3.8in 0 1in, width=2.4in]{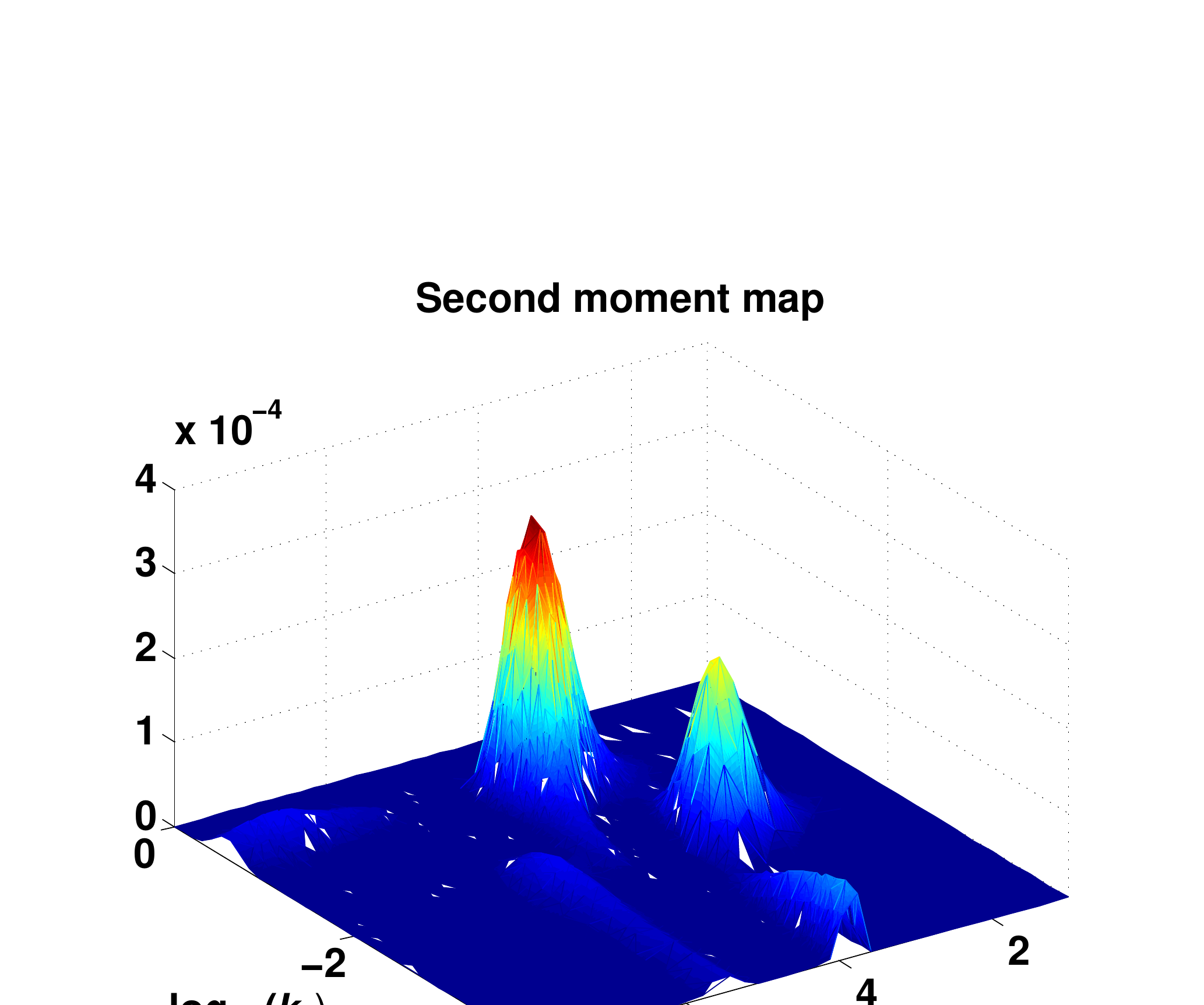}}
\subfigure[]{
\includegraphics[clip, trim=0 2.8in 0 2in, width=2.4in]{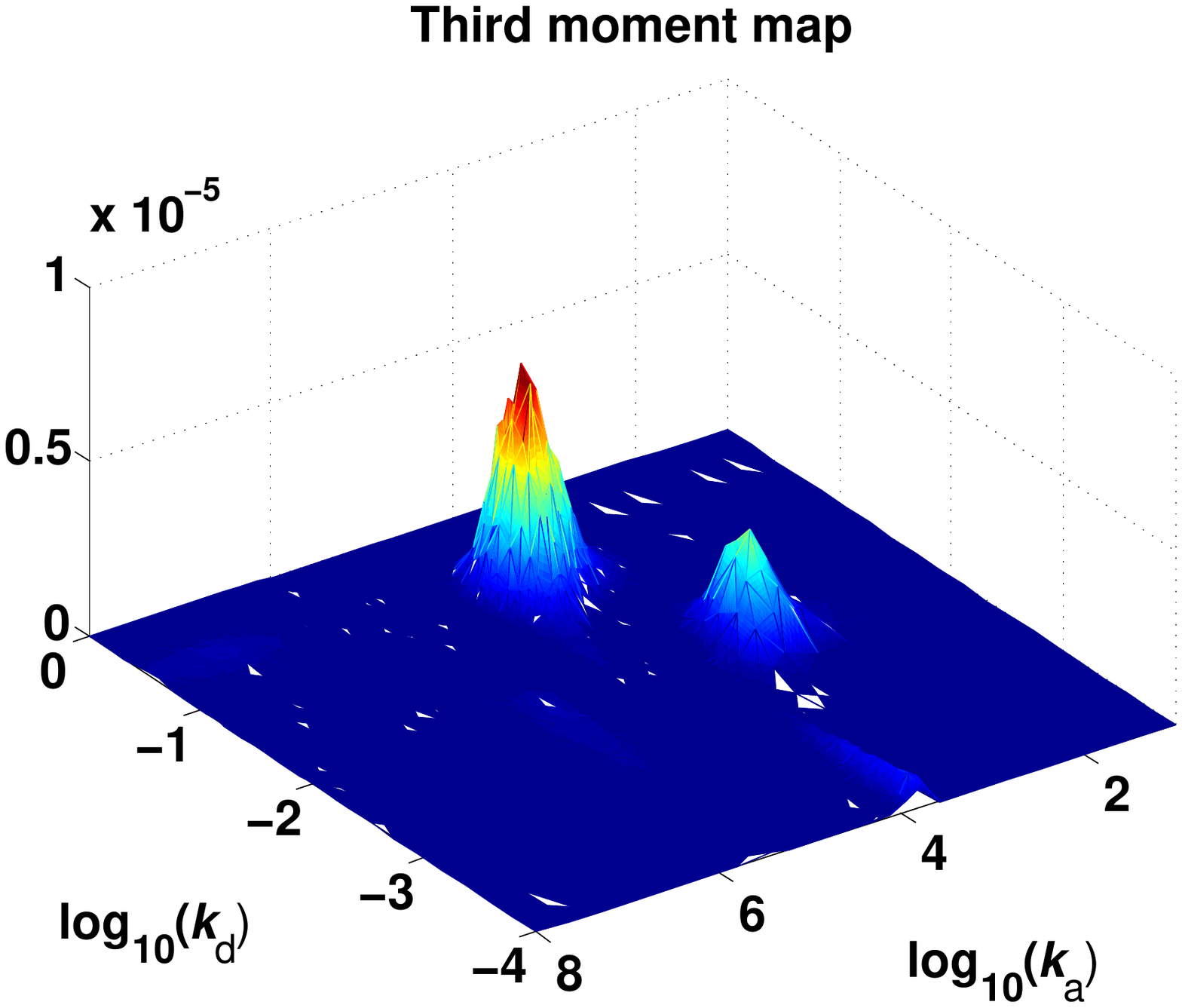}}
\subfigure[]{
\includegraphics[clip, trim=0 3.8in 0 1in, width=2.4in]{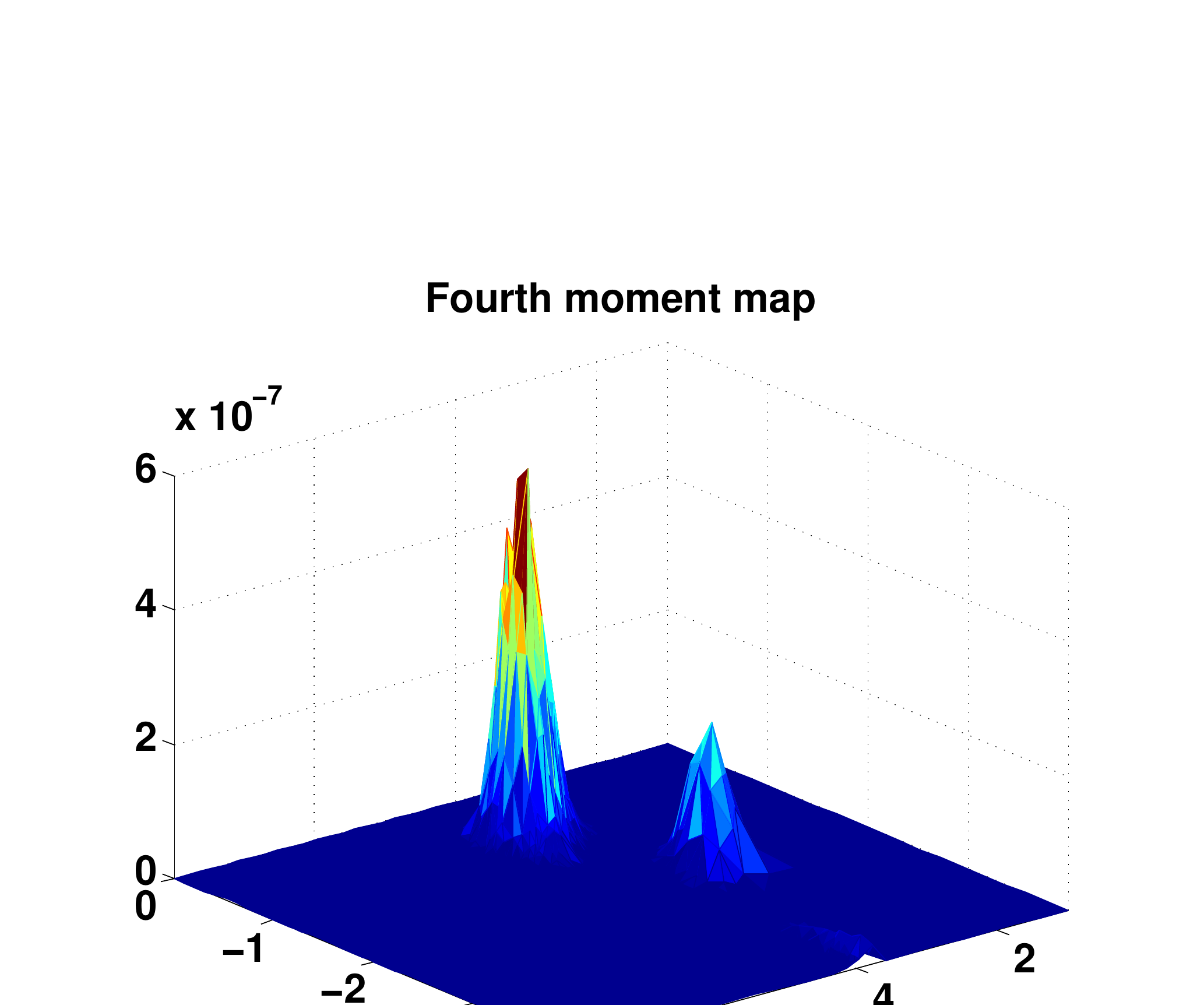}}
\caption[Moment maps]{Moment maps: (a) 0.5-th moment map. (b) Second moment map. (c) Third moment map. (d) Fourth moment map.}
\label{moments}
\end{figure}

\begin{figure}[!htb]
\centering
\subfigure[]{
\includegraphics[clip, trim=1.5in 2.8in 1.5in 2.8in, width=1.8in]{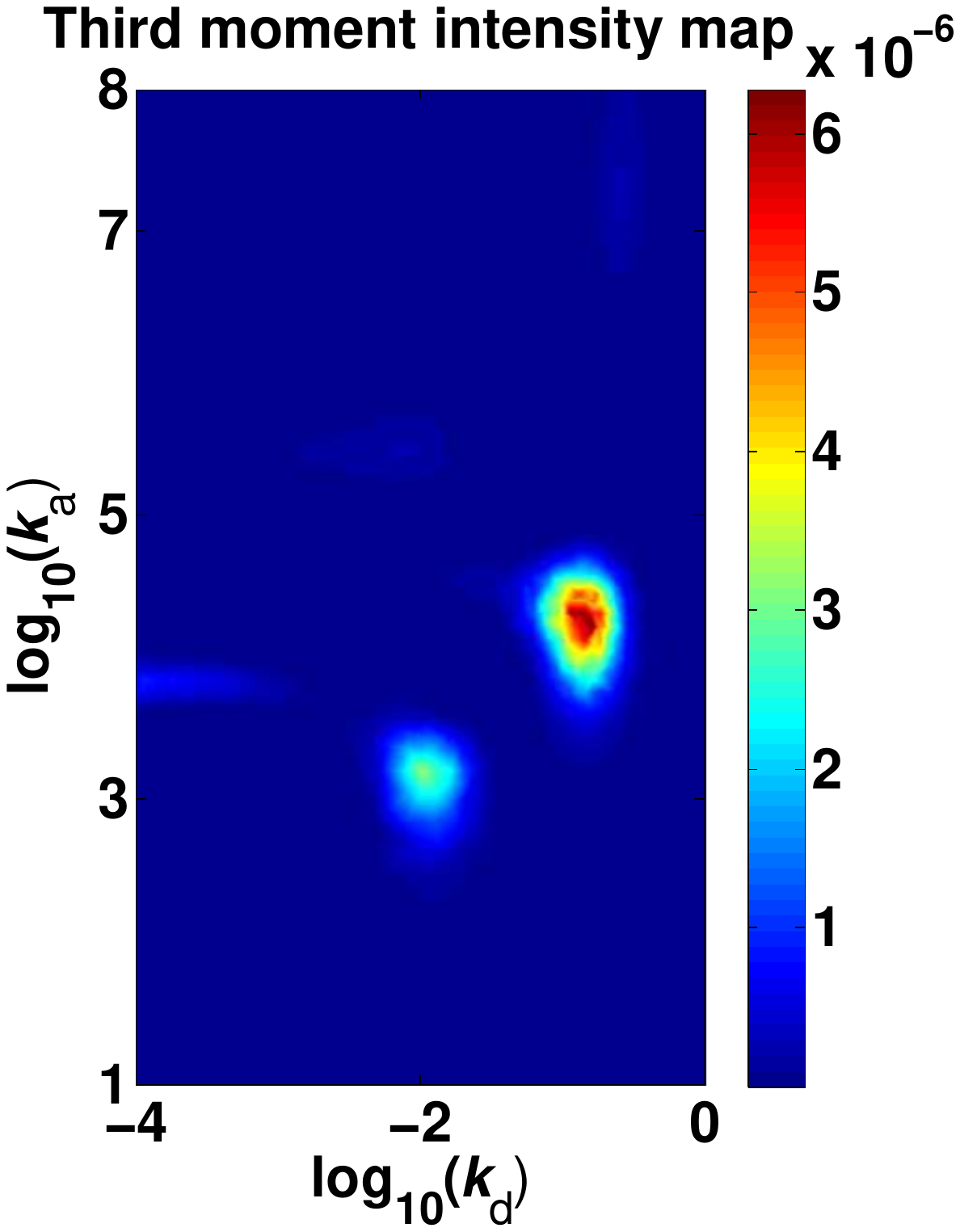}}
\subfigure[]{
\includegraphics[clip, trim=1.5in 2.8in 1.5in 2.8in, width=1.8in]{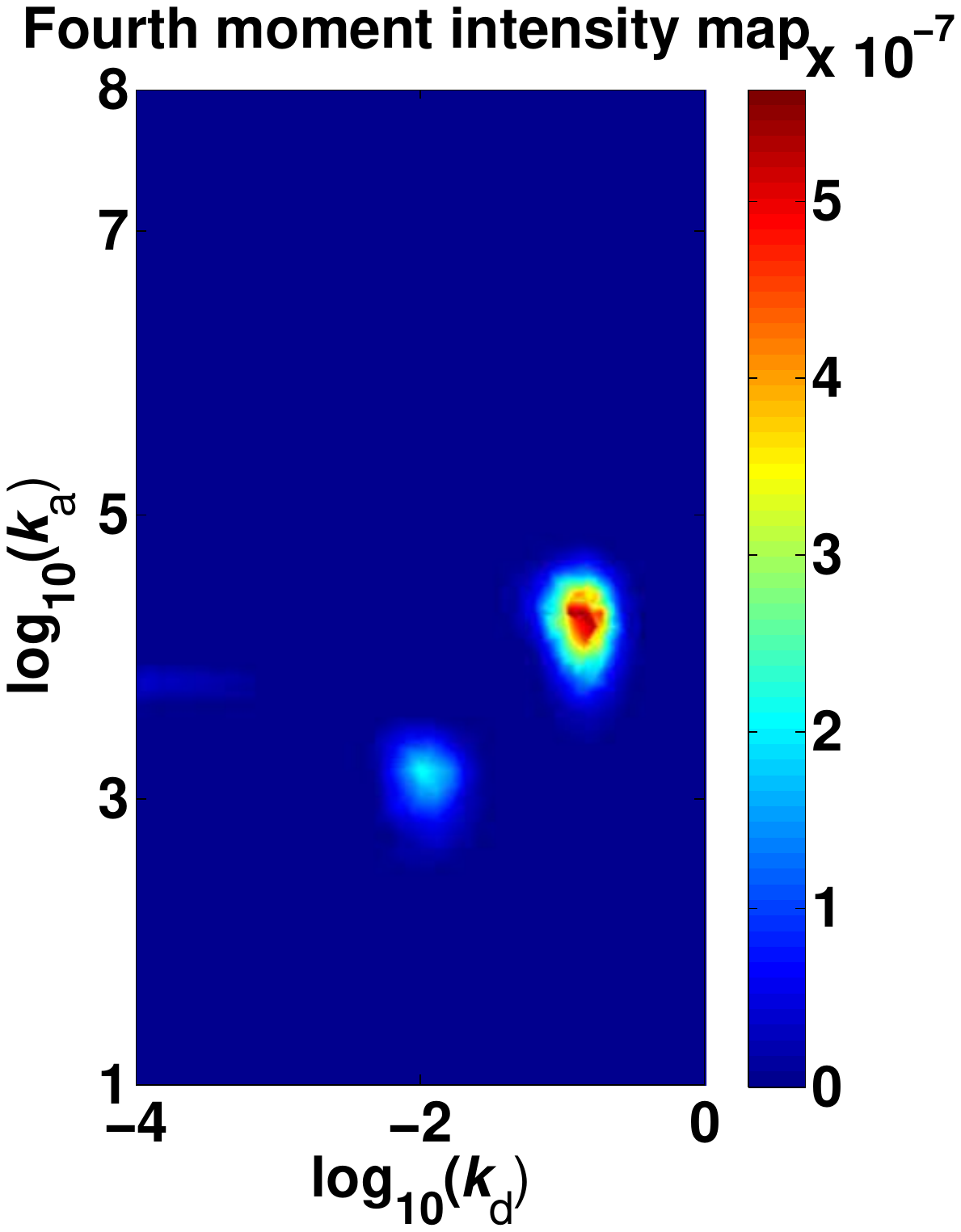}}
\caption[High order moment intensity maps]{High order moment intensity maps: (a) Third moment intensity map. (b) Fourth moment intensity map.}
\label{momentsIntensive}
\end{figure}

Now, consider our second approach to automatically determine the IPoI. In this approach, instead of using the mean map, we consider the other moment maps. Figure \ref{moments} displays the different moment maps from the 0.5-th moment map to the fourth moment map. As can be seen in the evolution of the moment maps in Figure \ref{moments}, the higher the moment, the sharper the map. For the map higher than the third moment, the moment map only has two peaks, and the other flatter local peaks that appear in the lower moment maps vanish. This phenomenon implies that the high moment map filters the redundant information, and highlights the interaction information. From this, we can conclude that the high moment map represents IPoR in a more efficient way. However, the third moment map provides almost the full picture of the IPoR, i.e., moment maps beyond that might not be necessary; see Figure \ref{momentsIntensive}, which shows the similarity between the third and fourth moment intensity maps. The high order moment intensity maps~-- see Figure \ref{momentsIntensive} ((a) for the third moment and (b) for the fourth moment)~-- indicate that there are two interactions in the PTH system. The corresponding association and dissociation rate constants are spotted at $(-2.0,3.2)$ and $(-0.9, 4.4)$ respectively. These results coincide with the conclusion obtained by the TCM.

Finally, we remark that in this real data application, we did not show the lower and upper rate constant maps since they present very similar sharps to the mean map. Hence, the lower and upper rate constant maps do not contain more information about the IPoR of the PTH system. However, for other biosensor systems, we remain inclined to suggest performing this analysis since it is possible that the sharps of the lower and upper rate constant maps behave differently, which will be interesting to study when looking at the advanced properties of IPoR for biosensor systems.

\section{Conclusion}

In this paper, an inverse problem of estimating the interaction number as well as their rate constants in biosensor systems has been studied in detail. Unlike the conventional deterministic mathematical models in the field of biosensor systems, we propose a statistical model based on the Bayesian framework. By engaging a large number of hyper parameters in the model, the full Bayesian solution results in a very high dimensional parameter space, which leads to a serious computational problem in terms of sampling the posterior distribution of the rate constants. To overcome this obstacle, we develop an Adaptive Variational Bayesian Approach (AVBA) for estimating the rate constants in biosensor systems. Numerical examples for both synthetic and real data demonstrate that the AVBA is stable and efficient. By combining the approach with the thresholding contour method or high order moment map method, the AVBA has shown its potential in capturing the intrinsic property of the interaction of a biosensor system; that is, it can provide us with a method to automatically find the number of interactions in a biosensor system, and the value of association and dissociation rate constants corresponding to these interactions.

It should be noted that the AVBA is a single step method that could accurately resolve the two underlying interactions without \emph{a priori} assumptions of the existence of a parallel interactions or the range of expected kinetic parameters, while the classical regression analysis (the parallel reactions model) of the kinetic sets produces correct results depending on the starting dissociation (e.g., $k_d$ in this paper) values. Although more experience with other biosensor systems would be needed to better understand the potential and limitations of the AVBA, its applications in several artificial problems and the PTH system seem very promising. By exploiting the full kinetic data set available, the obtained two-dimensional kinetic and affinity distributions have a higher resolution than the corresponding affinity distributions based on the isotherm analysis alone. We believe this will provide a useful tool for the study of the interactions by affinity biosensors.

\appendix

\section{Rate constant map theory}
\label{Kinetics}

First, let us discuss the kinetics for biosensors. Consider the ``1-to-1'' kinetic model for the binding process. For each interaction we have
\begin{equation*}
[A] +[L] \autoleftrightharpoons{$k_a$}{$k_d$} [AL],
\label{equilibrium}
\end{equation*}
where $k_a$ ($k_d$) denotes the association (dissociation) rate constant, and $[A]$, $[L]$, and $[AL]$ represent the concentrations of the analyte, ligand, and complex, respectively. In this work, the analyte $A$ is injected and flushed over the surface in such a way that the concentration $[A](t)$ can be assumed to be constant during the study. The amount of free ligand will decrease with time according to $[L](t)=[L](0)-[AL](t)$. Suppose that the sensor response $R$ is proportional to the complex concentration $[AL](t)$, i.e. $R(t)=\lambda \cdot [AL](t)$, where $\lambda$ is a constant. Denote by $R_{max}=\lambda \cdot [L](0)$. The rate of complex formation will be
\begin{equation}
\frac{d R(t)}{dt} = k_a \cdot [A](t) \cdot \left( R_{max} - R(t) \right) - k_d \cdot R(t)
\label{rateEq3}
\end{equation}
assuming that the mass transfer kinetics are infinitely fast. Set $C=[A](t)$ ($C$ is a constant, as mentioned above) and $R(t_0)=0$, the solution to (\ref{rateEq3}) is
\begin{equation}
R(t) = R_{max} \cdot \frac{k_a C}{k_d + k_a C} \cdot\left( 1- e^{(k_d + k_a C)(t-t_0)} \right).
\label{rateEq4}
\end{equation}

Now, let us develop a \emph{Rate Constant Map Theory} in the biosensor system. Assume that the binding can be described by an ``m-to-n'' kinetic model, i.e., we have $m$ analyte and $n$ binding sites on the biosensor surface and first order kinetics. Denote by $(k_{a,i}, k_{d,j})$ the pair of association and dissociation constants for the interaction between the $i$th analyte and $j$th binding site. Let $R_{i,j}(t)$ be the response at time $t$ of a complex with association constant $k_{a,i}$ and dissociation constant $k_{d,j}$. Then, according to (\ref{rateEq4}) we have
\begin{equation}\label{Rij}
R_{i,j}(t) = \left\{
\begin{array}{>{\displaystyle}l>{\displaystyle}l}
& 0, \qquad t\leq t_0 + \Delta t, \\
& R^{max}_{i,j}  \frac{k_{a,i} C}{k_{d,j}+ k_{a,i} C} \left( 1- e^{-(k_{d,j}+ k_{a,i} C)(t-t_0)} \right), \\
& \qquad\qquad\qquad\qquad\qquad t_0+ \Delta t < t \leq t_0+ t_{inj}+ \Delta t, \\ &
R^{max}_{i,j} \frac{k_{a,i} C}{k_{d,j}+ k_{a,i} C} \left( 1- e^{-(k_{d,j}+ k_{a,i} C) t_{inj}} \right) e^{-k_{d,j}(t-t_0-t_{inj})} , \\
& \qquad\qquad\qquad\qquad\qquad\qquad\qquad  t> t_0+ t_{inj}+ \Delta t,
\end{array}\right.
\end{equation}
where constant $C$ is the concentration of the analyte, $t_0$ is the time when the injection of the analyte begins, and $t_{inj}$ is the injection time. The adjustment parameter $\Delta t$ is a time delay that accounts for the fact that it usually takes some time for the detector to respond to the injection. Constant $R^{max}_{i,j}$ is the total surface binding capacity, corresponding to association and dissociation constants $k_{a,i}$ and $k_{d,j}$, i.e., the detector response when every binding site on the biosensor surface has formed a complex with the analyte.

We now use the functions $R_{i,j}$ to make an approximation of the measured sensorgrams $R_{obs}$, by assuming that the total measured response, $R_{obs}$, of a system can be written as a linear combination of some individual responses, namely $R_{obs}  = \sum^{m,n}_{i,j=1} R_{i,j}$. If we let $m,n\to+\infty$ in the above equation, we get the integral equation (\ref{IntegralEq}). Finally, we remark that the function $f(k_a,k_d)$, which is the generalization of the total surface binding capacity $\{R^{max}_{i,j}\}$, is known as the (continuous) rate constant map. See \ctp{Svitel-2003} for details.

\section{Proof of Theorem 2}\label{Appendix2}

Since $\sigma^2_{\mathbf{c}}$ obeys the inverse Gamma distribution during the iteration by the definition of $\sigma^{-2}_{\mathbf{c},k}$ in (\ref{minimizer_c}), we have the uniform boundedness for the sequence $\left\{ \sigma^{-2}_{\mathbf{c},k} \right\}$, i.e. (note that $(\mathbf{L} \mathbf{c})^T (\mathbf{L} \mathbf{c})\equiv \|\mathbf{L} \mathbf{c}\|^2\geq0$)
\begin{equation*}
0\leq \sigma^{-2}_{\mathbf{c},k} = \frac{\alpha_{\mathbf{c}} + \frac{n}{2}}{\beta_{\mathbf{c}} + \frac{1}{2} \mathbb{E}_{q^{k}(\mathbf{c})} \left[ (\mathbf{L} \mathbf{c})^T (\mathbf{L} \mathbf{c}) \right]} \leq \frac{\alpha_{\mathbf{c}} + \frac{n}{2}}{\beta_{\mathbf{c}}}.
\end{equation*}
By the same argument, parameters $\{\sigma^2_{k,j}, j=1, ..., N_C\}$ (note that $\sigma^{-2}_{k,j}:=[\Sigma^{-1}_{k}]_{j,j}$) are also uniformly boundedness during the iterations. Therefore, a subsequence exists, denoted by $\{\sigma^{-2}_{\mathbf{c},k_i}, \sigma^2_{k_i,j}\}$ and fixed parameters $\{\sigma^{-2}_{\mathbf{c},*}, \sigma^2_{*,j}\}$ such that
\begin{equation*}
\lim_{i\to\infty} \sigma^{-2}_{\mathbf{c},k_i} = \sigma^{-2}_{\mathbf{c},*}, \quad \lim_{i\to\infty} \sigma^2_{k_i,j} = \sigma^2_{*,j}, ~ j=1, ..., N_C.
\end{equation*}
Since $\mathbf{c}_k$ and $\Sigma^{\mathbf{c}}_{k}$ are solely determined by parameters $\{\sigma^{-2}_{\mathbf{c},k_i}, \sigma^2_{k_i,j}\}$, we can deduce that $\lim_{i\to\infty} q^{k_i} (\mathbf{c}) = q^{*} (\mathbf{c})$. Using the formula (\ref{minimizer_c}), we can conclude $\lim_{i\to\infty} q^{k_i} (\sigma^{2}_{\mathbf{c}}) = q^{*} (\sigma^{2}_{\mathbf{c}})$ and $\lim_{i\to\infty} q^{k_i} (\sigma^2_{j}) = q^{*} (\sigma^2_{j}), j=1, ..., N_C$.

Now, let us show that $q^*(\bullet):=q^*(\mathbf{c}) q^*(\sigma^2_1) \cdot\cdot\cdot q^*(\sigma^2_{N_C}) q^*(\sigma^2_{\mathbf{c}})$ is the stationary point of the KL distance functional $\mathcal{D}$. By the definition of subsequence $\{\sigma^{-2}_{\mathbf{c},k_i+1}, \sigma^2_{k_i+1,j}\}$ in (\ref{minimizer_c}), subsequence $\{\sigma^{-2}_{\mathbf{c},k_i+1}, \sigma^2_{k_i+1,j}, j=1, ..., N_C\}$ as well as subsequence $\{ q^{k_i+1} (\mathbf{c}), q^{k_i+1} (\sigma^{2}_{\mathbf{c}}), q^{k_i+1} (\sigma^2_{j}), j=1, ..., N_C\}$ also converge. Denote by
\begin{equation*}
\lim_{i\to\infty} q^{k_i+1} (\mathbf{c}) = q^{**} (\mathbf{c}), ~ \lim_{i\to\infty} q^{k_i+1} (\sigma^{2}_{\mathbf{c}}) = q^{**} (\sigma^{2}_{\mathbf{c}}), ~ \lim_{i\to\infty} q^{k_i+1} (\sigma^2_{j}) = q^{**} (\sigma^2_{j}).
\end{equation*}
Obviously, the following inequalities hold true
\begin{equation}\label{IneqAlgorithm}
\begin{array}{>{\displaystyle}l>{\displaystyle}l}
& \mathcal{D} \left( q^{**}(\mathbf{c})q^{**}(\sigma^2_1) \cdot\cdot\cdot q^{**}(\sigma^2_{N_C}) q^{**}(\sigma^2_{\mathbf{c}}) | p(\bullet,\mathbf{R}) \right) \leq \cdot\cdot\cdot  \\ &
\qquad\qquad \leq
\mathcal{D}\left( q^{*}(\mathbf{c})q^{**}(\sigma^2_1) \cdot\cdot\cdot q^{**}(\sigma^2_{N_C}) q^{**}(\sigma^2_{\mathbf{c}}) | p(\bullet,\mathbf{R}) \right) \leq \cdot\cdot\cdot \\ &
\qquad\qquad \leq  \mathcal{D}\left( q^{*}(\mathbf{c})q^{*}(\sigma^2_1) \cdot\cdot\cdot q^{*}(\sigma^2_{N_C}) q^{**}(\sigma^2_{\mathbf{c}}) | p(\bullet,\mathbf{R}) \right) \\ &
\qquad\qquad \leq \mathcal{D}\left( q^{*}(\mathbf{c})q^{*}(\sigma^2_1) \cdot\cdot\cdot q^{*}(\sigma^2_{N_C}) q^{*}(\sigma^2_{\mathbf{c}}) | p(\bullet,\mathbf{R}) \right).
\end{array}
\end{equation}

Let us show that $q^*(\bullet)=q^{**}(\bullet):= q^{**}(\mathbf{c}) q^{**}(\sigma^2_1) \cdot\cdot\cdot q^{**}(\sigma^2_{N_C}) q^{**}(\sigma^2_{\mathbf{c}})$. To this end, define $\pi$ as the algorithmic map of scheme (\ref{Algorithm1}), i.e., the solution operator between the two iterations that maps $q^k(\bullet)$ into $q^{k+1}(\bullet)$. The continuity of the functional $\mathcal{D}$ implies the closedness of the mapping  $\pi$. Hence, we have $q^{**}(\bullet)=\pi q^*(\bullet)$. This equation, together with inequalities (\ref{IneqAlgorithm}) and the monotone convergence of $\mathcal{D}$, implies
\begin{equation}\label{EqAlgorithm}
\begin{array}{>{\displaystyle}l>{\displaystyle}l}
& \mathcal{D} \left( q^{**}(\mathbf{c})q^{**}(\sigma^2_1) \cdot\cdot\cdot q^{**}(\sigma^2_{N_C}) q^{**}(\sigma^2_{\mathbf{c}}) | p(\bullet,\mathbf{R}) \right) \\ &
\qquad\qquad = \mathcal{D}\left( q^{*}(\mathbf{c})q^{**}(\sigma^2_1) \cdot\cdot\cdot q^{**}(\sigma^2_{N_C}) q^{**}(\sigma^2_{\mathbf{c}}) | p(\bullet,\mathbf{R}) \right) =  \cdot\cdot\cdot \\ &
\qquad\qquad = \mathcal{D}\left( q^{*}(\mathbf{c})q^{*}(\sigma^2_1) \cdot\cdot\cdot q^{*}(\sigma^2_{N_C}) q^{**}(\sigma^2_{\mathbf{c}}) | p(\bullet,\mathbf{R}) \right) \\ &
\qquad\qquad = \mathcal{D}\left( q^{*}(\mathbf{c})q^{*}(\sigma^2_1) \cdot\cdot\cdot q^{*}(\sigma^2_{N_C}) q^{*}(\sigma^2_{\mathbf{c}}) | p(\bullet,\mathbf{R}) \right).
\end{array}
\end{equation}

On the other hand, it holds true that for any $q(\mathbf{c})$:
\begin{equation}\label{IneqAlgorithm2}
\begin{array}{>{\displaystyle}l>{\displaystyle}l}
& \mathcal{D} \left( q^{**}(\mathbf{c})q^{**}(\sigma^2_1) \cdot\cdot\cdot q^{**}(\sigma^2_{N_C}) q^{**}(\sigma^2_{\mathbf{c}}) | p(\bullet,\mathbf{R}) \right) \\ &
\qquad\qquad
\leq \mathcal{D}\left( q(\mathbf{c}) q^{**}(\sigma^2_1) \cdot\cdot\cdot q^{**}(\sigma^2_{N_C}) q^{**}(\sigma^2_{\mathbf{c}}) | p(\bullet,\mathbf{R}) \right).
\end{array}
\end{equation}

By combining the above inequality, the first inequality in (\ref{EqAlgorithm}), the first identity in (\ref{IneqAlgorithm2}), and the strict multi-convexity of the functional $\mathcal{D}$, we can deduce that $q^{*}(\mathbf{c})=q^{**}(\mathbf{c})$. Furthermore, by repeating the preceding argument, we conclude that $q^{*}(\sigma^2_j)=q^{**}(\sigma^2_j)$ for $j=1, ..., N_C$, and $q^{*}(\sigma^2_{\mathbf{c}})=q^{**}(\sigma^2_{\mathbf{c}})$, which completes the proof of identity $q^{**}(\bullet)= q^*(\bullet)$, and thus $q^*(\bullet)$ is the stationary point of $\mathcal{D}$. Finally, the relation (\ref{result}) can easily be derived by applying $k_i\to\infty$ to the equations in (\ref{minimizer_c}).

\begin{supplement}
\stitle{Supplementary material for ``Estimating the Rate Constant from Biosensor Data via an Adaptive Variational Bayesian Approach"}
\slink[doi]{COMPLETED BY THE TYPESETTER}
\sdatatype{.pdf}
\sdescription{We provide additional material of the proof of Theorem 1, finite element approximation of integral equations, as well as a demonstration of our main algorithm.}
\end{supplement}

\section*{Acknowledgements}

We express our gratitude to the anonymous reviewers whose valuable comments and suggestions led to an improvement of the manuscript. The authors are grateful to Professor M\aa rten Gulliksson for the useful discussions. Cheng Li read and gave helpful comments on the paper. The authors would also like to thank Camilla K\"{a}ck, Marie Andersson and Teodor Aastrup from our KK H{\"O}G partner Attana AB for the PTH experiments.

\bibliographystyle{imsart-nameyear}
\bibliography{Biosensor}

\end{document}